\documentclass[10pt]{amsart}

\usepackage[small,nohug,heads=vee]{diagrams}
\diagramstyle[labelstyle=\scriptstyle]
\usepackage{amsmath,amstext,amsthm,amscd,amsopn,verbatim,amssymb }
\usepackage{hyperref}
\usepackage{graphicx}
\usepackage{multirow}
\usepackage{array}
\usepackage[all]{xy}
\usepackage{tabularx}
\usepackage[stable]{footmisc}

\newtheorem{Lem}{Lemma}
\newtheorem{Def}{Definition}
\newtheorem{Prop}{Proposition}
\newtheorem{Cor}{Corollary}
\newtheorem{Thm}{Theorem}

\newtheorem{Remark}[Thm]{Remark}
\newtheorem{Example}{Example}

\newtheorem{Def2}{Definition-lemma}

\DeclareMathOperator{\Ab}{\bf{Ab}}

\DeclareMathOperator{\Cat}{Cat}

\DeclareMathOperator{\DD}{\bf{D}}

\DeclareMathOperator{\Cyl}{Cyl}
\DeclareMathOperator{\Homm}{\bf{Hom}}

\DeclareMathOperator{\triang}{\bf{triang}}
\DeclareMathOperator{\thick}{\bf{thick}}

\DeclareMathOperator{\Hom2}{Hom}

\DeclareMathOperator{\Endd}{\bf{End}}

\DeclareMathOperator{\Modd}{\bf{Mod}}
\DeclareMathOperator{\prA'}{pr_{A[1]}}
\DeclareMathOperator{\prBs}{pr_{B[1]}}

\DeclareMathOperator{\Ext}{Ext}
\DeclareMathOperator{\T}{T}
\DeclareMathOperator{\ker2}{ker}
\DeclareMathOperator{\im2}{im}

\DeclareMathOperator{\End2}{End}

\DeclareMathOperator{\ker2}{ker}
\DeclareMathOperator{\im2}{im}

\DeclareMathOperator{\GD}{\bf{bG_{\mathbb K\c1}} }

\DeclareMathOperator{\C}{\bf{bG_{\mathbb K}}}

\DeclareMathOperator{\c1}{ [[\hbar ]] }

\begin{document}
\title{ Deformation Quantization of $A_\infty$-Morita equivalences}
\author{Andrea Ferrario}
\email{aferra@math.ethz.ch, andreaferrario011@gmail.com}
\address{Department of Mathematics, ETH Zurich, R\"amistrasse 101
8092 Zurich, Switzerland}

\maketitle
We show that Deformation Quantization of quadratic Poisson structures
preserves the $A_\infty$-Morita equivalence of a given pair of
Koszul dual  $A_\infty$-algebras. 

\tableofcontents

\section{Introduction}
In this paper we consider a finite dimensional vector space $X$ over a field $\mathbb K$ of characteristic 0 and  
the associative algebras with zero differentials $A=\mathrm S(X^*)$ resp. $B=\wedge(X)$ 
i.e. the symmetric algebra over $X^*$, resp. the exterior  algebra over $X$. 
For simplicity we choose $\mathbb K=\mathbb R$,$\mathbb C$. 
In \cite{CFFR} it is shown that it is possible to endow $K=\mathbb K$ with an  $A_\infty$-$A$-$B$-bimodule given by a codifferential $\mathrm d_K$ 
 whose Taylor components   are defined by  certain perturbative expansions in  Feynman diagrams. The expansions are written by considering
configuration spaces of points on the complex  upper half plane and  differential 1-forms called the 4-colors
propagators. This construction  and those  in \cite{CF0},\cite{CF} are
the first partial example of multi-brane generalization of the  results by M.~Kontsevich on Deformation Quantization of Poisson manifolds;
see \cite{Kont}. 
In \cite{CFFR} it is shown that the $A_\infty$-$A$-$B$-bimodule $(K,\mathrm d_K)$ is s.t. the classical Koszul duality between  $A$ and $B$ holds, 
i.e. there exists  isomorphisms
\begin{eqnarray*}
A\simeq \Ext_B(\mathbb K,\mathbb K),~~~B\simeq \Ext_A(\mathbb K,\mathbb K)^{op},
\end{eqnarray*} 
of  algebras: as left $A_\infty$-$A$-module and  right $A_\infty$-$B$-module $K$ is in fact the classical augmentation module.

Our first goal is to prove an $A_\infty$-derived Morita equivalence for the pair $(A,B)$ explicitly,
 i.e. the equivalence of certain triangulated subcategories 
of the derived categories of strictly unital $A_\infty$-right-modules over $A$ and $B$ by using the $A_\infty$-bimodule $(K,\mathrm d_K)$:
 $A$ and $B$ are just associative algebras with zero differential but we consider categories of $A_\infty$-modules over them.

It is natural to introduce a  bigrading on the triple $(A,K,B)$; the first grading is  cohomological;
the second grading is called internal; consequently
we consider only  bigraded  $A_\infty$-structures, i.e. bigraded $A_{\infty}$ modules, bimodules,  morphisms between them  etc.
By definition, the internal grading is preserved by the $A_\infty$-structures and morphisms between them.

The $A_\infty$-Morita equivalence for the pair $(A,B)$ has been already proved in \cite{Zhang}, where a more general result is shown.
In  \cite{Zhang}, (see prop. 1.14, 3.1. and thm. 5.7, 5.8 $loc.$ $cit.$)  the authors prove the  aforementioned 
equivalence by ``returning'' to the differential bigraded level by 
considering the derived categories of  differential bigraded modules over the enveloping algebras $UA$, resp. $UB$ of $A$ resp. $B$.
The  enveloping algebra $UA'$ of any bigraded $A_\infty$-algebra $A'$ is a differential bigraded algebra.
 It is introduced in \cite{Zhang} as the  theory of differential
bigraded algebras is, in general,  simpler than the theory of bigraded $A_\infty$-algebras.
Such an approach has the advantage of using the already well-known results on the  enveloping algebras
and (bar) resolutions of differential bigraded algebras.
On the other way, using this approach one introduces  the iterated use of the Koszul dual functor $E(\cdot)$, 
which associates to any augmented $A_\infty$-algebra  $A'$ its $A_{\infty}$-Koszul dual
$E(A')=\Hom2(UA',\mathbb K)$.
Moreover the  enveloping algebra $UA'$ is  a rather ``big''  bigraded object,
as by definition it is the cobar construction of the bar construction over $A'$.

Our approach is alternative to the one presented in \cite{Zhang};
we use the $A_\infty$-bimodule $K$ to prove the Morita equivalence at the $A_\infty$-level, without using the 
enveloping algebras $UA$, $UB$ and returning
to the differential bigraded level.

The key observation in our construction is that   the left derived derived actions (\cite{Kel2}, \cite{CFFR})
\begin{eqnarray*}
\mathrm L_A: A\rightarrow \underline{\End2}_B(K),~~~\mathrm R_B: B\rightarrow \underline{\End2}_A(K)^{op},
\end{eqnarray*}
 are quasi-isomorphisms of strictly unital $A_\infty$-$A$-$A$-bimodules and strictly unital $A_\infty$-$B$-$B$-bimodules; this is done
in  subsection~\ref{trywalking}. We use this fact to prove the equivalences of categories before and after deformation quantization.

The pair of functors inducing the equivalence is studied in subsection~\ref{stripped}. We define them by using the tensor products
$\bullet~\underline{\otimes}_A\bullet$, $\bullet~\underline{\otimes}_B\bullet$ of $A_\infty$-modules described in subsection~\ref{s-1}.
The main advantage of such ``pure''  $A_\infty$-approach, aside from the explicit use of the bimodule $K$,
is represented by the possibility of quantizing the equivalences: this is the content of section~\ref{quantized-A-infty}.
Let  $\hbar\pi$ be an  $\hbar$-formal quadratic Maurer-Cartan-element of cohomological degree 1  in  $\T_{poly}(X)\c1$, the ring of formal power
 series in $\hbar$ with coefficients in $\T^{\bullet}_{poly}(X)=S(X^{*})\otimes\wedge^{\bullet+1}(X)$. $\T_{poly}(X)\c1$ is a 
differential graded Lie algebra with zero differential and graded Lie braket $[ \cdot,\cdot ]_{\hbar}$ obtained by extending $\mathbb K\c1$-linearly 
the Schouten-Nijenhuis bracket $[ \cdot,\cdot ]$ on $\T_{poly}(X)$.
With such a choice of Poisson bivector the internal grading  on the triple on $(A,K,B)$ is  preserved; i.e.
using the ``2-branes Formality theorem'' contained in
\cite{CFFR} it follows that the quantizations $A_\hbar$, resp. $B_\hbar$ of $A$, resp. $B$ are   associative bigraded algebras
with zero differentials. The quantized bimodule $K_\hbar=(K\c1,\mathrm d_{K_\hbar})$
satisfies the quantized version of the Keller condition, and it is  a left $A_\hbar$-module and a right $B_\hbar$-module with zero differential.
Moreover, it is possible to quantize  straightforwardly the  bar resolutions
$A\underline{\otimes}_A K$, $K\underline{\otimes}_B B$ and the $A_\infty$-bimodules introduced in section~\ref{homee}.

In this ``deformed'' or quantized setting we introduce the categories $\Modd^{\infty}_{tf}(A_\hbar)$, resp. 
$\Modd^{\infty}_{tf}(B_\hbar)$ of strictly unital topologically free right $A_\infty$-$A_\hbar$-modules
(resp. $A_\infty$-$B_\hbar$-modules ).
 An object $N_\hbar$ in $\Modd^{\infty}_{tf}(A_\hbar)$ is a collection
$\{N^i_j\c1\}_{i,j\in\mathbb Z}$
of topologically free $\mathbb K\c1$-modules, endowed with a topological
$A_\infty$-$A_\hbar$-module structure, i.e. a quantized codifferential $\mathrm d_{M_\hbar}=\sum_{i\geq 0}\mathrm d^{(i)}_{M_\hbar}\hbar^i$ s.t.
$\mathrm d_{M_\hbar}\circ\mathrm d_{M_\hbar}=0$.

$\Modd^{\infty}_{tf}(A_\hbar)$ and $\Modd^{\infty}_{tf}(B_\hbar)$ are additive categories but they are not closed under taking
 cohomology. Quasi-isomorphisms in $\Modd^{\infty}_{tf}(A_\hbar)$ (and $\Modd^{\infty}_{tf}(B_\hbar)$) are then defined by considering the bigger abelian
category $\Modd_{bg}(\mathbb K\c1)$ of $all$ bigraded $\mathbb K\c1$-modules.

We define the  homotopy categories
$\mathcal H^{tf}_\infty(A_\hbar)$, $\mathcal H^{tf}_\infty(B_\hbar)$;
these are naturally  triangulated categories. To prove this we define topological cones and cylinders of topological $A_\infty$-morphisms;
all details are contained  in subsection~\ref{trtrtr}. Our approach is quite ``down-to-earth'': we adapt the definitions and results in \cite{GM}
to our topological $A_\infty$-setting.

We finish by introducing  ``the derived categories'' $\DD^{\infty}_{tf}(A_\hbar)$ resp. $\DD^{\infty}_{tf}(B_\hbar)$ as 
the localization at quasi-isomorphisms of $\mathcal H^{tf}_\infty(A_\hbar)$ resp.
 $\mathcal H^{tf}_\infty(B_\hbar)$: they are canonically endowed with a triangulated category structure
induced by the one on the corresponding homotopy category.

With $\triang^{\infty}_{A_{\hbar}}(M_{\hbar})$, $\triang^{\infty}_{B_{\hbar}}(N_{\hbar})$
we denote the full triangulated subcategories in $\DD^{\infty}_{tf}(A_\hbar)$ resp. $\DD^{\infty}_{tf}(B_\hbar)$
generated by $\{M_\hbar [i]\langle j\rangle\}_{i,j\in\mathbb Z}$ and similarly for $N_\hbar$, where $[\cdot]$, resp. $\langle\cdot\rangle$ are the 
shifts w.r.t. the cohomological resp. internal grading. With
$\thick^{\infty}_{A_{\hbar}}(M_{\hbar})$ and $\thick^{\infty}_{B_{\hbar}}(N_{\hbar})$ we denote their thickenings.
We recall their definitions in Appendix C. Let  $\tilde{\otimes}$ be the completed tensor product of bigraded $\mathbb K\c1$-modules w.r.t.
the $\hbar$-adic topology. The completed tensor products $\bullet~\tilde{\underline{\otimes}}_{A_\hbar}\bullet$, $\bullet~
\tilde{\underline{\otimes}}_{B_\hbar}\bullet$ of topological $A_\infty$-modules are defined accordingly.

The main result of these notes is then

\begin{Thm}
Let $X$ be a finite dimensional vector space over $\mathbb K=\mathbb R$,
or  $\mathbb C$ and $(A,K,B)$ be the triple of bigraded
$A_{\infty}$-structures with $A=S(X^{*})$, $B=\wedge(X)$ and $K=\mathbb K$ endowed 
with the $A_\infty$-$A$-$B$-bimodule structure given in \cite{CFFR}.
 By  $\pi_\hbar\in (T_{poly}(X)\c1,0,[\cdot,\cdot]_{\hbar})$ we denote an
    $\hbar$-formal  quadratic Poisson  bivector on $X$ 
and by $(A_{\hbar},K_{\hbar},B_{\hbar})$ the Deformation Quantization of $(A,K,B)$ w.r.t. $\pi_\hbar$.
 The triangulated functor
\begin{eqnarray*}
  \mathcal F_\hbar : \DD^{\infty}_{tf}(A_{\hbar})\rightarrow
\DD^{\infty}_{tf}(    B_{\hbar}), ~~~~~\mathcal
F_\hbar(\bullet)=\bullet~\underline{\tilde{\otimes}}_{A_\hbar} K_\hbar
\end{eqnarray*}
induces the equivalence of triangulated categories
\begin{eqnarray*}
\triang^{\infty}_{A_{\hbar}}(A_{\hbar})\simeq
\triang^{\infty}_{B_{\hbar}}(K_{\hbar}), & \thick^{\infty}_{A_{\hbar}}(A_{\hbar})\simeq \thick^{\infty}_{B_{\hbar}}(K_{\hbar}).
\end{eqnarray*}
Let $(\tilde{K},\mathrm d_{\tilde{K}})$ be the
$A_{\infty}$-$B$-$A$-bimodule
      with $\tilde{K}=K$ and $\mathrm d_{\tilde{K}}$ obtained from $\mathrm
d_K$ exchanging $A$ and $B$; the triangulated functor
\begin{eqnarray*}
  \mathcal F^{''}_\hbar : \DD^{\infty}_{tf}(B_{\hbar})\rightarrow
\DD^{\infty}_{tf}(    A_{\hbar}), ~~~~~\mathcal
F^{''}_\hbar(\bullet)=\bullet~\underline{\tilde{\otimes}}_{B_\hbar}
\tilde{K}_\hbar
\end{eqnarray*}
induces the equivalence of triangulated categories
\begin{eqnarray*}
\triang^{\infty}_{A_{\hbar}}(\tilde{K}_{\hbar})\simeq
\triang^{\infty}_{B_{\hbar}}(B_{\hbar}), &
\thick^{\infty}_{A_{\hbar}}(\tilde{K}_{\hbar}) \simeq \thick^{\infty}_{B_{\hbar}}(B_{\hbar}).
\end{eqnarray*}
\end{Thm}

In other words, Deformation Quantization of $\hbar$-formal  quadratic Poisson
 bivectors  preserves the $A_\infty$-Morita equivalence of the Koszul dual $A_\infty$-algebras $A$ and $B$.

In Appendix A we show the proof of prop.~\ref{End-bim}, while in Appendix B-C we prove thm.~\ref{Thm29} and thm.~\ref{Thm30} in some detail. 
Such proof are conceptually quite easy; using the very definition of the triangulated subcategories $\triang^{\infty}_{A_{\hbar}}(A_{\hbar})$
 $\dots$ $\thick^{\infty}_{B_{\hbar}}(K_{\hbar})$ we just need to check the commutativity of diagrams in which the quasi-isomorphisms 
of $A_\infty$-bimodules of section~\ref{homee} appear. Moreover, the proof of thm.~\ref{Thm30} is analogous to the one of thm.~\ref{Thm29},
 with mild changes.

\section{Acknowledgments}
We gratefully thank  D.~Calaque,  G.~Felder, B.~Keller, D.M.~Lu, C.~Rossi, P.~Shapira, M.~Van den Bergh,  for inspiring  discussions, 
 constructive criticism and useful e-mail exchange.

\section{Notation and Conventions}
Let $\mathbb{K}$ be a field of characteristic $0$. Throughout this work
we fix $\mathbb{K}=\mathbb{R}$ or $\mathbb{C}$.
Let $\C$ be the category of $\mathbb{Z}$-bigraded vector spaces, i.e.
collections $\{M^i_j\}_{i,j\in\mathbb{Z}}$ of
vector spaces over $\mathbb{K}$. The upper grading is also called the
``cohomological grading''. The lower index denotes the ``internal grading''.
The space of morphisms $\Hom2_{\C}(M,N)$  is the
$\mathbb{Z}$-bigraded vector space with ${(r,s)}^{th}$ component
\[
\Hom2^{r,s}_{\C}(M,N)=\prod_{n,m\in\mathbb{Z}}\Hom2_{\mathbb{K}}(M^{n}_{m},N^{n+r}_{m+s}),
\]
for every $r,s\in\mathbb{Z}^2$.

Any $f\in \Hom2^{r,s}_{\C}(M,N)$ is said to be a bigraded morphism of
bidegree $(r,s)$.
The identity morphisms in $\C$ are denoted simply by $1$.
For any object $M$ in $\C$, we denote by $M[n]$ the object in $\C$ such
that $(M[n])^{i}_{j}:=M^{i+n}_j$;
the degree -1 isomorphism $s:M\rightarrow M[1]$, $s(m):=m$ is called the
suspension map; its inverse of
degree 1 $s^{-1}:M[1]\rightarrow M$ is the desuspension. Both are
endofunctors of $\C$, with $(s^{-1})^{\otimes i}\circ s^{\otimes i}=(-1)^{\frac{i(i-1)}{2}}1$.
We use the short notation $sm$ for $s(m)\in M[1]$. The cohomological
degree of bihomogeneous elements of $M$ is
denoted by $|\cdot|$; in particular $|sm|=|m|-1$, for every $sm\in M[1]$.

Similarly, the object $M\langle j\rangle$ in $\C$ is s.t. ${M\langle
j\rangle}^n_m:=M^n_{m+j}$, for any $j\in\mathbb Z$.
It follows that $\Hom2^{r,s}_{\C}(M,N)=\Hom2^{0,0}_{\C}(M,N[r]\langle s\rangle)$.

The tensor product $M\otimes N$ of any two objects in $\C$ is the object
in $\C$ with bihomogeneous components
\begin{eqnarray*}
(M\otimes N)^{n}_{m} = \bigoplus_{p+q=n \atop r+s=m} M^{p}_{r} \otimes
N^{q}_{s}, 
\end{eqnarray*}
for every $n,m\in\mathbb{Z}$ with  $\otimes=\otimes_{\mathbb K}$.
Throughout this work we will use the shorthand conventions $m_1,\ldots,m_n\overset{!}=m_1\otimes\ldots\otimes m_n$,
and $(m_1|\cdots|m_n)\overset{!}=s(m_1)\otimes\cdots \otimes s(m_n)$, for any $m_1$, $\dots$ $m_n\in M\in\C$. So, in particular,
$(m_1,m_2|m_3)=m_1\otimes s(m_2)\otimes s(m_3)$ and
$(m_1|m_2,m_3)=s(m_1)\otimes s(m_2)\otimes m_3$. In what follows we assume that the Koszul sign rule holds.

\section{$A_\infty$-structures}

In this section we introduce $A_{\infty}$-structures from a purely
algebraic point of view.
We recall the concept of $A_{\infty}$-algebra, $A_{\infty}$-module,
$A_{\infty}$-bimodule and their morphisms.
We focus our attention on unital $A_{\infty}$-structures, augmented
$A_{\infty}$-algebras.
The tensor product of $A_{\infty}$-modules is also considered; it
contains the bar resolution of
a module over a given unital algebra as special case.
$A_{\infty}$-algebras have been introduced by Stasheff \cite{Sta}
in the sixties in algebraic topology; in the nineties they have been
further popularized by Kontsevich's \cite{Kont2}
in his Homological Mirror Symmetry conjecture.
The material here presented is  standard; we refer to
\cite{Kenji,Kel, GJ,Tradler} for all details, in particular
the definitions of coalgebras, coderivations, comodules etc...
For the interested reader, we just note that such definitions 
 can be deduced by taking the ``limit'' $\hbar=0$
in the formul\ae~ appearing in section~\ref{quantized-A-infty}.
Tensoring of $A_{\infty}$-bimodules has been introduced
explicitly in \cite{ART}, extending the case of right
$A_{\infty}$-modules contained in \cite{Kenji}.
In what follows we will consider only bigraded $A_\infty$-structures;
the rule of thumb is that the maps defining the $A_\infty$-structures themselves preserve the internal grading.
In this sense, there is not substantial difference between the graded and bigraded  case.

\subsubsection{$A_{\infty}$-algebras}\label{A-algebras}
Let $A$ be  an object of $\C$.
The coassociative counital tensor coalgebra on $A$ is the triple
\[
\mathcal B(A):=(\T^c(A[1]),\Delta, \epsilon),
\]
where $\T^c(A[1])=\oplus_{k\geq 0} A[1]^{\otimes k}$,
$\Delta:\T^c(A[1])\rightarrow \T^c(A[1])\otimes \T^c(A[1]) $
is the coassociative coproduct $\Delta(a_1|\dots|a_n)=1\otimes
(a_1|\dots|a_n)+(a_1|\dots|a_n)\otimes 1+\sum_{n'=1}^{n-1}
(a_1|\dots|a_{n'})\otimes (a_{n'+1}|\dots|a_n)$ and the counit
$\epsilon$ denotes the projection onto $\mathbb K$; by definition
$(\epsilon\otimes 1)\circ\Delta=(1\otimes \epsilon)\circ\Delta=1$.

\begin{Def}[J.~Stasheff,~\cite{Sta}]
An $A_{\infty}$-algebra is a pair $(A,\mathrm{d}_A)$, where $A$ is an
object of $\C$ and
${\mathrm d}_A$ is a bidegree $(1,0)$ coderivation on $\mathcal B(A)$ s.t.
\[
\mathrm{d}_A\circ \mathrm{d}_A=0.
\]
\end{Def}
By the lifting property of coderivations on $\mathcal B(A)$, such
$\mathrm d_A$ is uniquely determined by its Taylor components, i.e. the
family of morphisms $\bar{\mathrm d}_A^n:=\prA'\circ\mathrm
d_A|_{A[1]^{\otimes n}}$, $n\geq 0$, denoting by
$\prA'$ the projection $\prA':\T^c(A[1])\rightarrow A[1]$.

Then $\mathrm d_A\circ \mathrm d_A=0$ is equivalent to
\begin{eqnarray}
\sum_{s_1=0}^{k}\sum_{j=1}^{k-s_1+1}(-1)^{\epsilon}\bar{{\mathrm
d}}_A^{k-s_1+1}(a_1|\dots|a_{j-1},
\bar{{\mathrm
d}}_A^{s_1}(a_j|\dots|a_{s_1+j-1})|a_{s_1+j}|\dots|a_k)=0,\label{A1}
\end{eqnarray}
for every $k\geq 0$ and $(a_1,\dots,a_k)\in \mathcal{B}(A)$. The Koszul
sign is simply $\epsilon=\sum_{i=1}^{j-1}(|a_i|-1)$.
Equivalently, we can consider the bidegree $(2-n,0)$ maps $m_n$ defined through
\begin{eqnarray}
\bar{\mathrm d}^0_A&=-&s\circ m_0, \nonumber \\
\bar{{\mathrm d}}_A^{n}&=-&s\circ m_{n}\circ (s^{-1})^{\otimes
n},~~n\geq 1, \label{desus}
\end{eqnarray}

An $A_{\infty}$-algebra $(A,{\mathrm d}_A)$ is said to be flat if
${\mathrm d}_A^0=0$.
In this case $m_1$ is a differential and $m_2$ is associative up to
homotopy.
It reduces to an associative product on the cohomology $H(A)$ with
respect to $m_1$.
If a flat $A_{\infty}$-algebra is s.t. $m_3=m_4=\dots=0$, then it is a differential bigraded algebra.
If $(A,\mathrm d_A)$ is not flat, then it is called curved, with
curvature $\bar{\mathrm d}^0_A$
(or $\bar{\mathrm d}^0_A(1)$; we use both notations ).
In presence of non trivial curvature, $\bar{\mathrm{d}}^1_A$ is not a
differential. Any graded associative algebra $A$ s.t. $\bar{\mathrm{d}}^0_A(1)$ is a
degree $2$ element in the center of $A$ is a curved $A_{\infty}$-algebra.
Curvature appears naturally in Deformation Quantization: see for example
\cite{CFR}. Curved $A_{\infty}$-algebras are also related to
models in theoretical physics \cite{CT1}. 
With a little abuse of notation we introduce the following
\begin{Def}
Let $(A,\mathrm{d}_A)$ and $(B,\mathrm{d}_B)$ be $A_{\infty}$-algebras.
A morphism $F: A\rightarrow B$ of $A_{\infty}$-algebras is a morphism
$F\in \Hom2^{0,0}_{\C}(\mathcal{B}(A), \mathcal{B}(B) )$
of coassociative counital coalgebras s.t.
\[
F\circ \mathrm d_A=\mathrm d_B\circ F.
\]
\end{Def}
$F: \T^c(A[1])\rightarrow \T^c(B[1])$ is uniquely
determined by the family of morphisms $F_n: A[1]^{\otimes n}\rightarrow
B[1]$ s.t.
$\prBs\circ F|_{A[1]^{\otimes n}}=F_n$ and $F(1)=1$. The morphisms
$F_n$ are called the Taylor components of $F$.
$F\circ\mathrm{d}_A=\mathrm{d}_B\circ F$ is equivalent to a tower
of quadratic relations involving the Taylor components $F_{\bullet}$,
$\bar{\mathrm{d}}^{\bullet}_A$ and
$\bar{\mathrm{d}}^{\bullet}_B$ of $F$, $A$ and $B$,
respectively. If $(A,\mathrm{d}_A)$, resp. $(B,\mathrm{d}_B)$, are
curved $A_{\infty}$-algebras with curvature $\mathrm{d}^0_A$, resp.
$\mathrm{d}^0_B$, then, by definition of $F$: $F_1(\bar{\mathrm{d}}^0_A(1))= \bar{\mathrm{d}}^0_B(1)$.

It is useful to introduce the degree $1-n$ desuspended morphisms $f_n :
A^{\otimes n}\rightarrow B$ in $\C$, through
\begin{eqnarray}
F_n=s\circ f_n\circ (s^{-1})^{\otimes n}, \label{desus2}
\end{eqnarray}
for every $n\geq 0$.
A morphism $F: A\rightarrow B$ of $A_{\infty}$-algebras is said to be
strict if $F_n=0$ for $n\geq 2$. If $A$ and $B$ are flat,
$F$ is a quasi-isomorphism if $F_1$ is a quasi-isomorphism in $\C$.

\subsubsection{Units and augmentations in flat $A_{\infty}$-algebras}
Let $(A,{\mathrm d}_A)$ be an $A_{\infty}$-algebra;
the maps $m_n$, $n\geq 0$ and $f_m$, $m\geq 1$, have been defined in
(\ref{desus}), (\ref{desus2}).
\begin{Def}
An $A_{\infty}$-algebra $(A,{\mathrm d}_A)$ is said to be strictly
unital if it contains an element $1_A\in A^0_0$ s.t.
\begin{eqnarray*}
m_2(a ,1_A)=m_2(1_A, a)=a,
\end{eqnarray*}
for any $a\in A$ and $m_n(a_1,\dots,a_n)=0$ for $n \geq 3$ if $a_i=1$
for some $i=1,\dots,n$.
\end{Def}
We note that,  if $A$ is stricly unital, then $\bar{\mathrm d}^1_A(s1_A)=0$, also in presence of curvature on $A$.

A morphism $F: A_1\rightarrow A_2$ of strictly unital
$A_{\infty}$-algebras is said to be
strictly~unital if
\begin{eqnarray*}
f_1(1_{A_1})=1_{A_2},
\end{eqnarray*}
and $f_m(a_1,\dots,a_m)=0$ for $m \geq 2$ if $a_i=1_{A_1}$ for some
$i=1,\dots,m$.
In particular, it follows that $\bar{\mathrm d}^1_B(F_1(1_A))=0$.
\begin{Lem}
Any strictly unital flat $A_{\infty}$-algebra $A$ with unit $1_A$ comes
equipped with a strict strictly unital morphism
$\eta : K \rightarrow A$, sending the unity $1$ of the ground field
$\mathbb{K}$ to $1_A$.
\end{Lem}
This allows us to introduce the following
\begin{Def}
A strictly unital flat $A_{\infty}$-algebra $(A,d_A)$ with unit $1_A$ is
augmented if there exists a
strictly unital $A_{\infty}$-algebra morphism $\epsilon: A\rightarrow
K$, s.t. $\epsilon \circ \eta =1$.
\footnote{For any $A_{\infty}$-algebra $B$,
the identity morphism $1: B\rightarrow B$ is the
strict $A_{\infty}$-morphism with non trivial Taylor component
$\bar{1}^1(b)=b$, for every $b\in B$.}
\end{Def}
We note that the morphism $\epsilon\circ \eta$ is strict as $\epsilon$
is strictly unital.
If $A$ is an augmented $A_{\infty}$-algebra with augmentation
$\epsilon$, then we call $\ker2 \epsilon_1$ the augmentation ideal of $A$.


\subsubsection{$A_{\infty}$-modules and $A_{\infty}$-bimodules}
In this subsection $(A,{\mathrm d}_{A})$ and $(B,{\mathrm d}_{B})$ are
$A_{\infty}$-algebras.

\begin{Def}
A left $A_{\infty}$-$A$-module is pair $(M,\mathrm{d}_M)$, where $M$ is
an object in $\C$ and
$\mathrm{d}_M\in \Hom2^{1,0}_{\C}(\mathcal{L}(M),\mathcal{L}(M))$ is a
codifferential on
$\mathcal{L}(M):=\T(A[1])\otimes M[1]$ s.t.
\[
\mathrm{d}_M\circ\mathrm d_M=0.
\]
\end{Def}
As in the case of morphisms and coderivations on the tensor coalgebra
$\T^c(V)$, the codifferential $\mathrm{d}_M$ is uniquely
determined by its Taylor components
${\bar{\mathrm d}_M}^{s}: A[1]^{\otimes s}\otimes M[1]\rightarrow M[1]$,
$s\geq 0$, $via$
\begin{eqnarray*}
\mathrm{d}^k_M&=&\sum_{s_1=0}^k\sum_{j=1}^{k-s_1+1}1^{\otimes
j-1}\otimes \bar{\mathrm d}^{s_1}_A\otimes 1^{\otimes k-s_1-j+1}+
\sum_{s=0}^k1^{\otimes k-s}\otimes \bar{\mathrm d}^s_M,
\end{eqnarray*}
where the $\bar{\mathrm d}^{s_1}_A$ denote the Taylor components of the
coderivation ${\mathrm d}_A$
defining the $A_{\infty}$-algebra structure on $A$.

Let $(M,\mathrm d_M)$ be a left $A_{\infty}$-$A$-module.
$\mathrm{d}_M\circ\mathrm{d}_M=0$ is equivalent to
\begin{eqnarray}
&&\sum_{s_1=0}^k\sum_{j=1}^{k-s_1+1}(-1)^{\epsilon_1} {\bar{\mathrm
d}_M}^{k-s_1+1}(a_1|\dots|a_{j-1},\bar{\mathrm d}^{s_1}_A(a_j|
\dots|a_{s_1+j-1})|a_{s_1+j}|\dots|a_k|m )\nonumber+ \\
&&\sum_{s_2=0}^k(-1)^{\epsilon_2}\bar{\mathrm
d}^{k-s_2}_M(a_1|\dots|a_{k-s_2},\bar{\mathrm
d}^{s_2}_M(a_{k-s_2+1}|\dots|a_{k}|m))=0,
\label{modsin}
\end{eqnarray}
with $\epsilon_1=\sum_{i=1}^{j-1}(|a_i|-1)$,
$\epsilon_2=\sum_{i=1}^{k-s_2}(|a_i|-1)$.

\begin{Remark}
With obvious changes it is possible to define right
$A_{\infty}$-$A$-modules on the right $\mathcal{B}(A)$-counital
comodule $\mathcal{R}(M)=M[1]\otimes \T(A[1])$.
\end{Remark}
If $A$ is curved then $\bar{\mathrm{d}}^1_M(\mathrm{d}^0_A(1),sm)+
\bar{\mathrm{d}}^0_M(\mathrm{d}^0_M(sm))=0$, i.e. in presence of non trivial curvature $\mathrm{d}^0_A(1)$,
$\bar{\mathrm{d}}^0_M$ is not a differential on $M[1]$.
\begin{Def}
A morphism $F: M\rightarrow N$ of left $A_{\infty}$-modules
$(M,\mathrm{d}_M)$, $ (N,\mathrm{d}_N)$ is a morphism
$F\in \Hom2^{0,0}_{\C}(\mathcal{L}(M), \mathcal{L}(N) )$ of
left-$\mathcal{B}(A)$-counital-comodules s.t.

\begin{eqnarray}
F\circ\mathrm{d}_M=\mathrm{d}_N\circ F.\label{qwe}
\end{eqnarray}

\end{Def}

Any morphism $F: M\rightarrow N$ of left $A_{\infty}$-modules is
uniquely determined by its Taylor components
$F_n: A[1]^{\otimes n}\otimes M[1]\rightarrow N[1]$. Eq. (\ref{qwe}) is
equivalent to a tower of quadratic
relations involving the Taylor components $F_n$, ${\bar{\mathrm
d}_M}^{\bullet}$, ${\bar{\mathrm d}_N}^{\bullet}$
and ${\bar{\mathrm d}_A}^{\bullet}$; if $A$ is curved then
\[
F_0(\bar{\mathrm{d}}^0_M(sm))+F_1(\bar{\mathrm{d}}^0_A(1),sm)=\bar{\mathrm{d}}^0_N(F_0(sm));
\]
i.e. in presence of non trivial curvature $\mathrm{d}^0_A(1)$, $F_0:
M[1]\rightarrow N[1]$ does not commute with $\mathrm{d}^0_M$ and $\mathrm{d}^0_N$ (which are not differentials).

\begin{Def}
A morphism $F: M\rightarrow N$ of left-$A_{\infty}$-$A$-modules
is said to be strict if $F_n=0$ for $n\geq 1$. If $A$ is flat, $F$ is a
quasi-isomorphism if $F_0$ is a quasi-isomorphism.
\end{Def}

\begin{Def}
An $A_{\infty}$-$A$-$B$-bimodule is a pair $(M,\mathrm{d}_M)$, where $M$
is an object in $\C$ and $\mathrm{d}_M\in
\Hom2^{1,0}_{\C}(\mathcal{B}(M),\mathcal{B}(M) )$ is a codifferential on
$\mathcal{B}(M)=\T(A[1])\otimes M[1]\otimes \T(B[1])$ s.t.
\[
\mathrm{d}_M\circ\mathrm d_M=0.
\]
\end{Def}
Once again, it is possible to show that the codifferential
$\mathrm{d}_M$ is uniquely determined by the Taylor components
$\bar{\mathrm{d}}^{k,l}_M:=A[1]^{\otimes k}\otimes M[1]\otimes
B[1]^{\otimes l}\rightarrow M[1]$, $k,l\geq 0$, with
\begin{eqnarray*}
\mathrm{d}_M^{k,l}&=&\sum_{s_1=0}^k\sum_{j=1}^{k-s_1+1}1^{\otimes
j-1}\otimes \bar{{\mathrm d}}^{s_1}_A\otimes 1^{\otimes
k-j-s_1+1+l+1}+
 \sum_{s_2=0}^l\sum_{j=1}^{l-s_2+1}1^{\otimes k+1}\otimes 1^{\otimes
j-1} \otimes \bar{{\mathrm d}}^{s_2}_B\otimes 1^{\otimes
l-j-s_2+1}+\nonumber \\
&&\sum_{s_3=0}^k\sum_{s_4=0}^{l}1^{\otimes k-s_3}\otimes {\bar{\mathrm
d}}^{s_3,s_4}_M\otimes 1^{\otimes l-s_4}.\label{Bim1}
\end{eqnarray*}
Then $\mathrm d_M\circ \mathrm d_M=0$ is equivalent to a tower of quadratic relations similar to \eqref{modsin}, with due differences.
In presence of non trivial curvatures on $A$ and/or $B$, then
$\bar{\mathrm{d}}^{0,0}_M$ is not a differential on $M[1]$.

\begin{Lem}[\cite{ART}]\label{t`}
Let $(A,\mathrm{d}_A)$, $(B,\mathrm{d}_B)$ be $A_{\infty}$-algebras and
$(M,\mathrm{d}_M)$ be an $A_{\infty}$-$A$-$B$-bimodule.
\begin{itemize}
\item If $B$ is flat, then the family $\bar{\mathrm{d}}^{k,0}_M:
A[1]^{\otimes k}\otimes M[1]\rightarrow M[1]$
defines a left-$A_{\infty}$-$A$-module structure on $M$.
\item If $A$ is flat, then the family $\bar{\mathrm{d}}^{0,l}_M:
M[1]\otimes B[1]^{\otimes l}\rightarrow M[1]$, $l\geq 0$,
defines a right-$A_{\infty}$-$B$-module structure on $M$.
\end{itemize}
\end{Lem}

\begin{Remark}
Every $A_{\infty}$-algebra $(A,{\mathrm d}_A)$ is an
$A_{\infty}$-$A$-$A$-bimodule
with $A_{\infty}$-bimodule structure given by the Taylor components
$\bar{\mathrm d}_A^{k,l}: A[1]^{\otimes k}\otimes A[1] \otimes
A[1]^{\otimes l}\rightarrow A[1]$, with
\begin{eqnarray*}
\bar{\mathrm d}_A^{k,l}:=\bar{\mathrm d}_A^{k+l+1}. 
\end{eqnarray*}

\end{Remark}

\begin{Def} Let $(M,{\mathrm d}_M)$ and $(N,{\mathrm d}_N)$ be two
$A_{\infty}$-$A$-$B$-bimodules, with $\mathcal{B}(M)=\T(A[1])\otimes M[1]\otimes \T(B[1])$, and similarly for 
$\mathcal{B}(N)$.
A morphism of $A_{\infty}$-$A$-$B$-bimodules is a morphism $F\in\Hom2^{0,0}_{\C}(\mathcal{B}(M),\mathcal{B}(N))$
of $\mathcal{B}(A)$-$\mathcal{B}(B)$-codifferential-
counital bicomodules s.t.
\begin{eqnarray*}
F\circ {\mathrm d}_M={\mathrm d}_N\circ F. 
\end{eqnarray*}
\end{Def}
Any $A_{\infty}$-$A$-$B$-bimodule morphism $F$ is uniquely determined by
its Taylor components
$\bar{F}^{k,l}:A[1]^{\otimes k}\otimes M[1]\otimes B[1]^{\otimes
l}\rightarrow N[1]$, $k,l\geq 0$. Explicitly
\begin{eqnarray*}
F^{k,l}=\sum_{s_3=0}^k\sum_{s_4=0}^{l}1^{\otimes k-s_3}\otimes
\bar{F}^{s_3,s_4}\otimes1^{\otimes l-s_4},
\end{eqnarray*}
where $F^{k,l}:=F|_{A[1]^{\otimes k}\otimes M[1]\otimes B[1]^{\otimes l}}$.
If $A$, resp. $B$ are curved with curvature $\bar{\mathrm d}^0_A$, resp.
$\bar{\mathrm d}^0_B$, then $F^{0,0}$ does not commute with $\bar{\mathrm d}^{0,0}_M$ and $\bar{\mathrm d}^{0,0}_N$
(which are not differentials).

\subsubsection{Units in $A_{\infty}$-modules}
Let $(A,\mathrm{d}_A)$ be a strictly unital $A_{\infty}$-algebra with
unit $1_A$ and $(M,\mathrm{d}_M)$ a left $A_{\infty}$-$A$ module.
We introduce the desuspended maps
\[
\bar{\mathrm d}^{l}_M=-s\circ d_l^M \circ (s^{-1})^{\otimes l}, ~~l\geq 0.
\]

\begin{Def}
The module $(M,\mathrm{d}_M)$ is strictly unital if
\begin{eqnarray*}
d_1^M(1_A,m)=m,
\end{eqnarray*}
for every $m\in M$ and $d_n^M(a_1,\dots,a_n,m)=0$ for $n\geq 2$ with
$a_i=1_A$ for some $i=1,\dots,n$.
\end{Def}
Similar considerations hold for right $A_{\infty}$-modules.
A strictly unital morphism of strictly unital $A_{\infty}$-modules is an
$A_{\infty}$-morphism $F$ s.t.
\begin{eqnarray*}
F^n(a_1|\dots|a_n|m)=0, ~~~n\geq 2
\end{eqnarray*}
with $a_i=1_A$ for some $i=1,\dots,n$ and $F^1(1_A|m)=-sm$.

Similar definitions hold for unital $A_{\infty}$-bimodules over strictly
unital $A_{\infty}$-algebras.

\subsubsection{Homotopies of strictly unital $A_\infty$-modules}
 Let $A$ be a strictly unital $A_\infty$-algebra and $(M,\mathrm d_M)$, $(N,\mathrm d_N)$ be strictly unital
 $A_\infty$-$A$-modules.
Let $f,g:M\rightarrow N$ be morphisms of
 $A_\infty$-$A$-modules; we say that $M$ and $N$
 are  $A_\infty$-homotopy equivalent (alternatively:
$A_\infty$-homotopic) if there exists an
$A_\infty$-homotopy between them, i.e. a bidegree (-1,0) morphism $H\in \Hom2^{-1,0}_{\C}(M[1]\otimes \T(A[1])
,N[1]\otimes\T(A[1]))$ of counital $\T(A[1])$-comodules, s.t.
\[
f_\hbar-g_\hbar=\mathrm d_{N_\hbar}\circ H_\hbar+H_\hbar\circ \mathrm
d_{M_\hbar}.
\]

\subsubsection{The tensor product of $A_\infty$-bimodules}\label{s-1}
We consider now three $A_\infty$-algebras $(A,\mathrm d_A)$, $(B,\mathrm
d_B)$ and
$(C, \mathrm d_C)$. Furthermore, we introduce an
$A_\infty$-$A$-$B$-bimodule $(K_1,\mathrm d_{K_1})$ and an
$A_\infty$-$B$-$C$-bimodule
$(K_2,\mathrm d_{K_2})$.
\begin{Def}
The tensor product $K_1\underline{\otimes}_B K_2$ of $K_1$ and $K_2$
over $B$ is the object 
\[
K_1\underline{\otimes}_B K_2=K_1\otimes
\T(B[1])\otimes K_2
\]
in $\C$ .
\end{Def}

\begin{Prop}[\cite{ART}] \label{prop1}
$K_1\underline\otimes_B K_2$ is endowed with an
$A_\infty$-$A$-$C$-bimodule structure given by the codifferential $\mathrm
d_{K_1\underline\otimes_B K_2}$ with Taylor components
 $\bar{\mathrm d}^{m,n}_{K_1\underline\otimes_B
K_2}$  given by
{\footnotesize
\begin{eqnarray}
&&\bar{\mathrm d}_{K_1\underline\otimes_B
K_2}^{m,n}\!(a_1|\cdots|a_m|k_1\otimes
(b_1|\cdots|b_q)\otimes k_2|c_1|\cdots|c_n)=0,\quad m,n>0\nonumber\\
&&\bar{\mathrm d}_{K_1\underline\otimes_B
K_2}^{m,0}\!(a_1|\cdots|a_m|k_1\otimes
(b_1|\cdots|b_q)\otimes k_2)=\nonumber\\
&&\sum_{l=0}^q s\left(s^{-1}(\bar{\mathrm
d}^{m,l}_{K_1}(a_1|\cdots|a_m|k_1|b_1|\cdots|b_l))\otimes
(b_{l+1}|\cdots|b_q)\otimes k_2\right),\quad m>0\nonumber\\
&&\bar{\mathrm d}_{K_1\underline\otimes_B K_2}^{0,n}\!(k_1\otimes
(b_1|\cdots|b_q)\otimes
k_2|c_1|\cdots|c_n)=\nonumber\\
&&(-1)^{|k_1|+\sum_{j=1}^q(|b_j|-1)}\sum_{l=0}^q
s\left(k_1\otimes (b_1|\cdots|b_l)\otimes s^{-1}(\bar{\mathrm
d}^{q-l,n}_{K_2}(b_{l+1}|\cdots|b_q|k_2|c_1|\cdots|c_n)\right),\quad
n>0,\nonumber\\
&&\bar{\mathrm d}_{K_1\underline\otimes_B K_2}^{0,0}\!\left(s(k_1\otimes
(b_1|\cdots|b_q)\otimes k_2)\right)=\nonumber\\
&&\sum_{l=0}^q s\left(s^{-1}(\mathrm
d_{K_2}^{0,l}(k_1|b_1|\cdots|b_l)\otimes (b_{l+1}|\cdots|b_q)\otimes
k_2\right)+\nonumber\\
&&\sum_{0\leq l\leq q\atop 0\leq p\leq
q-l} (-1)^{(|k_1|-1)+\sum_{j=1}^l (|b_j|-1)}s\!\left(k_1\otimes
(b_1|\cdots|\bar{\mathrm d}^p_B(b_{l+1}|\cdots|b_{l+p})|\cdots|b_q)\otimes
k_2\right)+\nonumber\\
&&(-1)^{|k_1|+\sum_{j=1}^q(|b_j|-1)}\sum_{l=0}^q
s\!\left(k_1\otimes (b_1|\cdots|b_l)\otimes s^{-1}(\bar{\mathrm
d}^{q-l,0}_{K_2}(b_{l+1}|\cdots|b_q|k_2)\right).\label{eq-tayl-tens}
\end{eqnarray} }
\end{Prop}

\begin{Cor}
Let $K_1$ be an $A_{\infty}$-$A$-$B$-bimodule, $K_2$ an
$A_{\infty}$-$B$-$C$-bimodule and $K_3$ an $A_{\infty}$-$C$-$D$-bimodule.
The tensor product of $A_{\infty}$-bimodules is associative, i.e. there
exists a strict $A_{\infty}$-$A$-$D$-bimodule morphism
\begin{eqnarray*}
\Theta: (K_1\underline{\otimes}_B K_2)\underline{\otimes}_C K_3
\rightarrow K_1\underline{\otimes}_B (K_2\underline{\otimes}_C K_3)
\end{eqnarray*}
which induces an isomorphism of objects in $\C$.
\end{Cor}

\subsubsection{The $A_\infty$-bar constructions of an
$A_\infty$-bimodule}\label{ss-1-1}
We consider two $A_\infty$-algebras $(A,\mathrm{d}_A)$,
$(B,\mathrm{d}_B)$ and an
$A_\infty$-$A$-$B$-bimodule $(M,\mathrm{d}_M)$.
We recall that $A$ can be canonically endowed with an
$A_{\infty}$-$A$-$A$-bimodule structure; see Remark 11.
Same holds for $B$, with due changes.

\begin{Def}\label{11}
The $A_{\infty}$-$A$-$B$-bimodule $(A\underline{\otimes}_A
M,\mathrm{d}_{A\underline{\otimes}_A M})$ is called the $A_{\infty}$-bar
construction of $(M,\mathrm{d}_M)$ as left $A_{\infty}$-$A$-module.
Similarly, the $A_{\infty}$-$A$-$B$-bimodule
$(M\underline{\otimes}_B B,\mathrm{d}_{M\underline{\otimes}_B B})$ is
called the $A_{\infty}$-bar
construction of $(M,\mathrm{d}_M)$ as right $A_{\infty}$-$B$-module.
\end{Def}
By definition, both $A\underline{\otimes}_A M$ and
$M\underline{\otimes}_B B$ are $A_{\infty}$-$A$-$B$-bimodules.
Let $A$ and $B$ be unital algebras and $M$ an $A$-$B$-bimodule.
Then $A\underline{\otimes }_A M$ is the bar resolution of $M$ as left
$A$-module.
Similarly, $M\underline{\otimes}_B B$ is the bar resolution of $M$ as
right $B$-module.
\begin{Prop}[\cite{ART}]\label{p-bar}
Let $(A,\mathrm{d}_A)$, $(B,\mathrm{d}_B)$ be $A_\infty$-algebras and
$(M,\mathrm{d}_M)$ be an $A_\infty$-$A$-$B$-bimodule.
There exists a natural morphism
\[
\mu: A\underline{\otimes}_A M \rightarrow M,
\]
of $A_\infty$-$A$-$B$-bimodules.
If $A$, $B$ are both flat, and $A$, $M$ are left unital as
$A_{\infty}$-$A$-module, then the
morphism $\mu$ is a quasi-isomorphism.
\end{Prop}

\subsubsection{On the $A_\infty$-bar construction: a remark }

We continue our analysis of the $A_\infty$-bar constructions and the
 morphisms
\[
\mu_A :A\underline{\otimes}_A K\rightarrow K, ~~~\mu_B:
K\underline{\otimes}_B B\rightarrow K
\]
of strictly unital $A_\infty$-$A$-$B$-bimodules introduced in the above
subsection. In the following lemma
we restrict to the case of augmented associative algebras with zero
differentials as they will appear later on.

\begin{Lem}\label{Lem2}
Let $(A,\mathrm d_A)$ and $(B,\mathrm d_B)$ be augmented
associative algebras with zero differential
and $(K,\mathrm d_K)$ be a strictly unital $A_\infty$-$A$-$B$-bimodule.

\begin{itemize}
\item There exists  strictly unital quasi-isomorphisms
\begin{eqnarray*}
K\rightarrow A\underline{\otimes}_A K, & K\rightarrow K\underline{\otimes}_B B,
\end{eqnarray*}
of $A_\infty$-$A$-$B$-bimodules.
\end{itemize}
\end{Lem}
\begin{proof}
We denote by
\[
A_{+}:=\ker2\epsilon_A,~~~ B_{+}:=\ker2\epsilon_B
\]
the augmentation ideals in $A$, resp. $B$, denoting by $\epsilon_A$
resp. $\epsilon_B$ the augmentation maps on $A$, resp. $B$.
We recall that the augmentation maps are morphisms of algebras. So the
augmentation ideals are subalgebras.

We prove the first statement. The second is similar. Let
\[
A\underline{\otimes}_{A_{+}} K=\bigoplus_{n\geq 0} A\otimes A_{+}[1]^{\otimes n}\otimes K,
\]
be the normalized bar resolution of $K$.
$A\underline{\otimes}_{A_{+}} K$ is a strictly unital
$A_\infty$-$A$-$B$-bimodule. There exists a
strict quasi-isomorphism
\[
\mathcal I: A\underline{\otimes}_{A_{+}} K\rightarrow
A\underline{\otimes}_A K
\]
of $A_\infty$-$A$-$B$-bimodules; it is the natural inclusion.
The quasi-isomorphism $K\rightarrow A\underline{\otimes}_A K$ is the
composition
\[
K\stackrel{\bar{\Phi}}{\rightarrow} A\underline{\otimes}_{A_{+}}
K\stackrel{\mathcal I}{\rightarrow} A\underline{\otimes}_A K
\]
where the (bidegree $(0,0)$) morphism $\bar{\Phi}$ is given as follows. Its
$(n,m)$-th
Taylor component $\bar{\Phi}_{n,m}: A[1]^{\otimes n}\otimes K[1]\otimes
B[1]^{\otimes m}
\rightarrow (A\underline{\otimes}_{A_{+}} K)[1]$ is simply
\[
\bar{\Phi}_{n,m}=s\circ \Phi_{n,m}\circ {(s^{-1})}^{n+m+1}
\]
with
\[
\Phi_{n,m}(a_1,\dots,a_n,k,b_1,\dots,b_m)=0~~\mbox{if}~~m\geq 1,
\]
and
\[
\Phi_{n,0}(a_1,\dots,a_n,k)=\left\{ \begin{array}{ccc}
(-1)^{\sum_{i=1}^n(|a_i|-1)}(1,a_1,\dots,a_n,k) & \mbox{if} & a_i\in
A_{+},~\mbox{for all}~ i=1,\dots, n. \\
0 & &\mbox{otherwise} \\
\end{array}\right.
\]
Note that $\Phi_{n,0}$ is of bidegree $(-n,0)$; $\Phi$ is strictly unital by
construction.
To check that
\begin{eqnarray}
\bar{\Phi}\circ \mathrm d_K=\mathrm d_{A\underline{\otimes}_{A_{+}} K} \circ
\bar{\Phi};\label{ice}
\end{eqnarray}
is straightforward. We need to consider \eqref{ice} on all the possible strings of elements  
$ (a_1|\dots|a_m|k|b_1|\dots|b_n)\in\T(A[1])\otimes K[1]\otimes \T(B[1])$, $n,m\geq 0$
paying attention whether $ (a_1|\dots|a_n)\in A_{+}[1]^{\otimes n}$ or $sa_i\in \mathbb K[1]$, for some $i$.  
As $\Phi_{0,0}(1)=1\otimes 1$, then $\bar{\Phi}$ is a quasi-isomorphism.
\end{proof}

\begin{Cor}Let $A$, $B$ and $K$ be as above.
\begin{itemize}
\item $K$ and $A\underline{\otimes}_A K$ are homotopy equivalent as strictly unital $A_\infty$-$A$-$B$-bimodules.
\item $K$ and $K\underline{\otimes}_B B$ are homotopy equivalent as strictly unital $A_\infty$-$A$-$B$-bimodules.
\end{itemize}

\end{Cor}
\begin{proof}
We prove the first statement; the second is analogous. We want to show that there exists a strictly unital $A_\infty$-homotopy 
$\bar{H}: K\rightarrow A\underline{\otimes}_A K$ of $A_\infty$-$A$-$B$-bimodules, s.t.
\begin{eqnarray*}
\bar{\Phi}\circ\mu_A &=& 1+\mathrm d_{A\underline{\otimes}_A K}\circ \bar{H} + \bar{H}\circ \mathrm d_{A\underline{\otimes}_A K},\\
\mu_A\circ\bar{\Phi} &=& 1,
\end{eqnarray*}
denoting by $\mu_A$ the $A_\infty$-morphism appearing in prop.~\ref{p-bar} and by $\Phi$ the one appearing in  lem.~\ref{Lem2}.
The bidegree $(-1,0)$ Taylor components $\bar{H}_{m,n}: A[1]^{\otimes m}\otimes  (A\underline{\otimes}_A K)[1]\otimes B[1]^{\otimes n}\rightarrow  
 (A\underline{\otimes}_A K)[1]$ are given by
$\bar{H}_{n,m}=0$ if $m\geq 1$, and
\[
\bar{H}_{n,0}(a_1|\dots|a_n|(a,{a'}_1|\dots|{a'}_q,k))=\left\{ \begin{array}{ccc}
s(1,a_1|\dots|a_n|a|a'_1|\dots|a'_q,k) & \mbox{if} & a_i\in
A_{+},~\mbox{for all}~ i=1,\dots, n. \\
0 & &\mbox{otherwise} \\
\end{array}\right.
\]
$\mu_A\circ\bar{\Phi} = 1$ easily follows  as $K$ is strictly unital. The equality involving $H$ is long to prove, but straightforward.
By definition, the identity $1$ is a strict and strictly unital $A_\infty$-morphism.
\end{proof}

 \section{The triple $(A,K,B)$}\label{section-AKB}

Let $X$ be a finite dimensional vector space over the field
$\mathbb{K}=\mathbb{R},\mathbb{C}$. In \cite{CFFR} it is shown that, choosing a pair $(U,V)$ of
subspaces in $X$, then it is possible to introduce a pair $(A,B)$ of $A_\infty$-algebras associated to the subspaces
themselves $and$ an $A_\infty$-bimodule $K$ associated to the intersection $U\cap V$.
Choosing $(U,V)=(X,\{0\})$ we arrive at the pair of $A_\infty$-algebras
\[
A=S(X^*),~~~B=\wedge(X).
\]
$A$ and $B$ are objects in $\C$; let us discuss their bigrading. We put
\[
A=\bigoplus_{i\geq 0 }A_i,~~ A_i=A_i^0,
\]
where $A_i$ denotes the vector space of homogeneous polynomials of
degree $i$. It follows that $A_0=A_0^0=\mathbb{K}$.
$A$ is concentrated in cohomological degree $0$. The
$A_{\infty}$-structure on $A$ is encoded in a
codifferential $\mathrm d_A$ whose only non trivial Taylor component is
$\mathrm{\bar{d}}^2_A: A[1]^{\otimes 2}\rightarrow A[1]$.
For the exterior algebra $B$ we put
\[
B=\bigoplus_{i\geq 0} B^i,~~ B^i=B_{-i}^{i},
\]
with $B_{-i}^{i}:=\wedge^{i} X$. A bihomogeneous element $b\in B^i$
has bidegree $(i,-i)$.
Also in this case $B_0=B_0^0=\mathbb{K}$. The $A_{\infty}$-structure on
$B$ is encoded in a
codifferential $\mathrm d_B$ whose only non trivial Taylor component is
$\mathrm{\bar{d}}_B^2: B[1]^{\otimes 2}\rightarrow B[1]$.
In summary, the generators of $B$ are bihomogeneous of bidegree
$(1,-1)$;
the dual generators in $A$ are bihomogeneous with bidegree $(0,1)$.
Both $A$ and $B$ are augmented $A_{\infty}$-algebras with augmentation
ideals
$A_{+}=\bigoplus_{i\geq 1}A^0_i$ and $B_{+}=\bigoplus_{i\geq 1} B^i_{-i}$.
Moreover 
\begin{Prop}[\cite{CFFR}]
Let $X$ be a finite dimemnsional vector field over $\mathbb K$, $A=S(X^*)$, and $B=\wedge(X)$. There 
exists a one-dimensional strictly unital $A_\infty$-$A$-$B$-bimodule $K$ which, as a left $A$-module and as a right $B$-module, is
the augmentation module.
\end{Prop}
The $A_{\infty}$-$A$-$B$-bimodule structure on $K$ is specified by a
codifferential $\mathrm d_K$,
with Taylor components $\bar{\mathrm d}^{k,l}_K: A[1]^{\otimes k}\otimes
K[1]\otimes B[1]^{\otimes l}\rightarrow K[1]$.  We remind that, by definition,
$\mathrm d_K$ (and so $\bar{\mathrm d}^{k,l}_K$,
for every $k,l\geq 0$) is of cohomological degree $1$. The explicit
construction in terms of Feynman diagrams implies that
$\mathrm{\bar{d}}^{k,l}_K (a_1|\dots|a_k|1|b_1|\dots|b_l)$ is non
vanishing iff
\begin{eqnarray}
\sum_{i=1}^k \deg a_i=\sum_{i=1}^l|b_i|=k+l-1, \label{d2}
\end{eqnarray}
where $\deg a_i$ denotes the internal degree of the homogeneous polynomial
$a_i\in A$ and $|b_i|$
the cohomological grading of $b_i\in B^{|b_i|}$. But (\ref{d2})
implies that $\bar{\mathrm d}^{k,l}_K$ is
of degree $0$ w.r.t the internal grading on $A$, $B$ and $K$, for every
$k,l\geq 0$: we recall that suspension
and desuspension do not shift the internal degree.

The explicit construction of the codifferential $\mathrm d_K$ implies
that $K$ is a strictly unital $A_{\infty}$-$A$-$B$-bimodule.

\subsubsection{On the Keller condition for $(A,K,B)$}\label{trywalking}

We return to a more general setting.
\begin{Def}
Let $(A,\mathrm d_A)$ and $(B,\mathrm d_B)$ be flat
$A_{\infty}$-algebras and $(K,\mathrm d_K)$ be a right $A_{\infty}$-
$B$-module. We set $\mathcal R(K):=K[1]\otimes\T(B[1])$.
$(\underline{\End2}_{B}(K), \mathrm
d_{\underline{\End2}_{B}(K)})$ is the flat $A_{\infty}$-algebra defined
as follows. As bigraded object
\[
\underline{\End2}_{B}(K):=\Hom2_{\C}(\mathcal R(K),K[1]);
\]
the codifferential $\mathrm d_{\underline{\End2}_{B}(K)}$ has non
trivial Taylor
components
\begin{eqnarray}
\bar{\mathrm
d}^1_{\underline{\End2}_{B}(K)}&=&-s\circ\partial\circ
s^{-1},~~\partial(\varphi)=(-1)^{|\varphi|+1}(\varphi\circ \mathrm
d_K)+\mathrm d_K\circ \varphi,\nonumber\\
\bar{\mathrm
d}^2_{\underline{\End2}_{B}(K)}(\varphi|\psi)&=&(-1)^{|\varphi|}s(\varphi\circ\
\psi).\label{End_Tayl}
\end{eqnarray}
\end{Def}
We can define $(\underline{\End2}_{A}(K), \mathrm
d_{\underline{\End2}_{A}(K)})$ almost $verbatim$.

\begin{Prop}[\cite{Kont3},\cite{CFFR}]\label{prop4}
Let $(A,\mathrm d_A)$ and $(B,\mathrm d_B)$ be flat
$A_{\infty}$-algebras and $(K,\mathrm d_K)$ be a right $A_{\infty}$-
$B$-module. $K$ is an $A_{\infty}$-$A$-$B$-bimodule\footnote{We define a
codifferential ${\mathrm D}_K$ s.t.
$\bar{\mathrm D}^{0,l}_K=\bar{\mathrm d}^{l}_K$, for every $l\geq 0$.}
if and only if there exists a morphism
\begin{eqnarray*}
\mathrm L_A: A \rightarrow \underline{\End2}_{B}(K)
\end{eqnarray*}
of $A_{\infty}$-algebras.
\end{Prop}
\begin{proof} A detailed proof can be found in \cite{CFFR}; here we sketch it. 
Let $(K,\mathrm d_K)$ be endowed with an $A_{\infty}$-$A$-$B$-bimodule
structure $\mathrm D_K$ s.t.
$\bar{\mathrm D}^{0,l}_K=\bar{\mathrm d}^{l}_K$. The maps
\begin{eqnarray}
\mathrm L_A(a_1|\dots|a_k)\in \underline{\End2}_B(K)[1], ~~ \mathrm
L_A(a_1|\dots|a_k):=s\circ \mathcal L_A(a_1|\dots|a_k) \label{tcomp}
\end{eqnarray}
with $\mathcal L_A(a_1|\dots|a_k)$ of bidegree $(1,0)$ given by

\begin{equation}\label{ttt}
\mathcal L_A(a_1|\dots|a_k)(1|b_1|\dots|b_q):=\bar{\mathrm D}^{k,l}_K( a_1
|\dots | a_k |1|b_1|\dots |b_q),
\end{equation}
are the Taylor components of an $A_{\infty}$-algebra morphism
$\mathrm L_A: A \rightarrow \underline{\End2}_{B}(K)$, for every
$(a_1|\dots|a_k)\in A[1]^{\otimes k}$,
$(1|b_1|\dots |b_q)\in K[1]\otimes B[1]^{\otimes q}$ and $k\geq 1$,
$q\geq 0$.
Viceversa, let $\mathrm L_A: A \rightarrow \underline{\End2}_{B}(K)$ be an
$A_{\infty}$-algebra morphism with Taylor
components as in \eqref{tcomp} . Then the maps
$\bar{\mathrm D}^{k,l}_K$ in \eqref{ttt} are the Taylor components
of a codifferential $\mathrm D_K$ on $\T(A[1])\otimes K[1]\otimes
\T(B[1])$, extending the given right $A_{\infty}$-$B$-module structure
on $K$.
\end{proof}

We call $\mathrm L_A$ in prop.~\ref{prop4} the derived left $A$-action. A
 similar statement can be proved in the case of the derived right
$B$-action,
i.e. the $A_{\infty}$-algebra morphism $\mathrm R_B :B^{op}\rightarrow
\underline{\End2}_{A}(K)$ with obvious Taylor components. The
$A_{\infty}$-algebra $B^{op}$ has $A_{\infty}$-structure
canonically induced by the one on $B$, but the signs are not trivial. We
refer to \cite{Zhang} for all details. 
\begin{Def}[\cite{Kel2}]
Let $(A,\mathrm d_A)$ and $(B,\mathrm d_B)$ be flat
$A_{\infty}$-algebras and $(K,\mathrm d_K)$ be an
$A_{\infty}$-$A$-$B$-bimodule.
The triple $(A,K,B)$ satisfies the Keller condition if the
derived actions
\[
\mathrm L_A: A\rightarrow \underline{\End2}_B(K) ,
\]
and
\[
\mathrm R_B: B^{op}\rightarrow \underline{\End2}_A(K) ,
\]
are quasi-isomorphism of $A_{\infty}$-algebras.
\end{Def}

\subsubsection{The Keller condition for the triple $(A,K,B)$ }
Let $(A,K,B)$ be the triple of bigraded $A_{\infty}$-structures given in
section~\ref{section-AKB}. 
The bigrading on the triple $(A,K,B)$ is such that
\begin{eqnarray*}\underline{\End2}^{ i,j}_{B}(K)=\left\{\begin{array}{cc}
0 &  i+j<0\\
\Hom2^{0,0}_{\C}(K[1]\otimes
B[1]^{\otimes i+j},K[1][i]\langle j\rangle) & i+j\geq 0
\end{array}\right.
\end{eqnarray*} 
Note that $\underline{\End2}^{ 0,0}_{B}(K)\cong \mathbb{K}$ and
$\mathrm{\bar{d}}_K^{0,1}\in\underline{\End2}^{ 2,0}_{B}(K)$.
Similar considerations hold for $\underline{\End2}_{A}(K)$.
The derived left action $\mathrm L_A$ preserves the internal grading, by
definition. Moreover, for every $k\geq 1$ and $(a_1|\dots|a_k)\in A[1]^{\otimes
k}$, then $\mathrm L_A(a_1|\dots|a_k)$
is an element of $\underline{\End2}^{ n,m}_{B}(K)$, with
$(n,m):=(-k+1,\sum_{i=1}^k \deg a_i)$.

For any $l\geq 0$ and $(1|b_1|\dots|b_l)\in (K[1]\otimes B[1]^{\otimes
l})^a_b$, with $(a,b)=(-1+\sum_{i=1}^l |b_i|-l, -\sum_{i=1}^l |b_i|)$, we have
\[
\mathrm L_A(a_1|\dots|a_k)(1|b_1|\dots|b_l):=\mathrm
d^{k,l}_{K}(a_1|\dots|a_k|1|b_1|\dots|b_l)\in (K[1])^{n+a}_{m+b}.
\]
This implies that
\[
n+a+1=0\Rightarrow \sum_{i=1}^l|b_i|=k+l-1,~~~~~~~m+b=0\Rightarrow \sum_{i=1}^k \deg a_i=\sum_{i=1}^l|b_i|.
\]
In other words, the wordlength $l$ is uniquely determined by the constraint
$l=1-k+\sum_{i=1}^k \deg a_i$, for any choice of $(b_1|\dots|b_l)\in (B[1])^{\otimes l}$ as above.
This analysis applies to $\mathrm R_B$, with due changes. In \cite{CFFR} it is shown the important
\begin{Prop}
The triple $(A,K,B)$ given is section~\ref{section-AKB} is s.t.
the derived left $A$-action $\mathrm
L_A$ and the derived right $B$-action
$\mathrm R_B$ are quasi-isomorphism of strictly unital
$A_{\infty}$-algebras.
\end{Prop}
As in the proof of proposition~\ref{prop4} we introduce the notation
\[
\mathcal L_A(a_1|\dots|a_n)\in \underline{\End2}_{B}^{r,m}(K),
~~\mathcal L_A(a_1|\dots|a_n)(1|b_1|\dots|b_q)=\mathrm
d_K^{n,q}(a_1|\dots|a_n|1|b_1|\dots|b_q),
\]
i.e. $\mathrm L_A(a_1|\dots|a_n)=s\circ \mathcal L_A(a_1|\dots|a_n)$ and
$r=\sum_{i=1}^n (|a_i|-1)+1$,
$m=\sum_{i=1}^n \deg a_i$.
We note that $\mathcal L_A(a_1|\dots|a_n)$ is of cohomological degree +1.
 $A=S(X^{*})$ is canonically endowed with a strictly unital $A_\infty$-$A$-$A$-bimodule structure
$\mathrm {\tilde{d}}_A$ whose non trivial
Taylor components are $\bar{\mathrm {\tilde{d}}}^{(1,0)}_A=\bar{\mathrm
{\tilde{d}}}^{(0,1)}_A=\bar{\mathrm {\bar{d}}}^2_A$.

\begin{Prop}\label{End-bim}
There exists a strictly unital $A_{\infty}$-$A$-$A$-bimodule
structure $\mathrm
d_{\underline{\End2}_{B}(K)}$ on $\underline{\End2}_{B}(K)$ such that
the derived action $\mathrm L_A$ descends to a quasi-isomorphisms of
strictly unital
$A_{\infty}$-$A$-$A$-bimodules.
$\mathrm d_{\underline{\End2}_{B}(K)}$ has Taylor components
\begin{eqnarray*}
&&\bar{\mathrm d}_{\underline{\End2}_{B}(K)}^{0,0}= -s\circ\partial_{
\underline{\End2}_{B}(K) } \circ s^{-1},\\
&&\bar{\mathrm
d}_{\underline{\End2}_{B}(K)}^{n,0}(a_1|\dots|a_n|\varphi)=s\circ \mathrm
D_{\underline{\End2}_{B}(K)}^{n,0}(a_1|\dots|a_n|\varphi),~~~(n\geq 1)\\
&&\bar{\mathrm d}_{\underline{\End2}_{B}(K)}^{0,m}(\varphi|a_1|\dots|a_m)=
s\circ \mathrm
D_{\underline{\End2}_{B}(K)}^{0,m}(\varphi|a_1|\dots|a_m),~~~(m\geq 1)\\
\end{eqnarray*}
with
\begin{eqnarray*}
&&\partial_{\underline{\End2}_{B}(K)}(\varphi)=(-1)^{|\varphi|+1}
\varphi\circ \mathrm d_{K}+\mathrm
d_{K}\circ\varphi,\\
&&\mathrm D_{\underline{\End2}_{B}(K)}^{n,0}(a_1|\dots|a_n|\varphi)=
(-1)^{\sum_{i=1}^n(|a_i|-1)-1}\mathcal L_A(a_1|\dots|a_n)\circ\varphi,\\
&&\mathrm D_{\underline{\End2}_{B}(K)}^{0,m}(\varphi|a_1|\dots|a_m)=
(-1)^{|\varphi|}\varphi\circ
\mathcal L_A(a_1|\dots|a_m),
\end{eqnarray*}
and $\bar{\mathrm d}_{\underline{\End2}_{B}(K)}^{n,m}=0$, otherwise.\\

\end{Prop}
\begin{proof} See Appendix A.
\end{proof}
It can also be verified that the derived right-$B$ action $\mathrm
R_B$ descends to a quasi-isomorphism
of $A_{\infty}$-$B^{op}$-$B^{op}$-bimodules.


\section{$A_\infty$-Morita theory}\label{homee}

\subsubsection{On thm. 5.7. in \cite{Zhang}}

In this section we study the $A_{\infty}$-Morita theory for the triple
$(A,K,B)$.
Our approach to the Morita equivalence is purely $A_{\infty}$; all we need is the
$A_{\infty}$-$A$-$B$-bimodule structure on $K$ we described in the
previous section
to prove the equivalence of certain triangulated subcategories of
$A_{\infty}$-modules in the derived categories $\DD^{\infty}(A)$ and $\DD^{\infty}(B)$ of $A$ and $B$. 
The functors giving the equivalences are defined through the $A_\infty$-tensor product of $A_\infty$ modules 
and bimodules.
The formalism is quite simple, using the associativity of the $A_\infty$-tensor product.
The main advantage in using such ``pure'' $A_{\infty}$-approach is
represented by the fact that the computations
which follow are all explicit; the quasi-isomorphisms of
$A_{\infty}$-bimodules which are the core of
the equivalences are induced by the Keller condition on $(A,K,B)$.

\subsubsection{On some bigraded $A_{\infty}$-modules}

Let $M$ be an $A_{\infty}$-$A$-$B$-bimodule and $N$ be an
$A_{\infty}$-$B$-$C$-bimodule, where $A$, $B$ and $C$ are
$A_{\infty}$-algebras. We have already introduced the
$A_{\infty}$-$A$-$B$-bimodule $(\mathcal{B}_B(M),\mathrm
d_{\mathcal{B}_B(M)})$,
where $\mathcal{B}_B(M):=M\underline{\otimes}_B B$, calling it the
$A_{\infty}$-bar
construction of $M$ as right $A_{\infty}$-$B$-module. It is an
$A_{\infty}$-right-$B$-module. If $B$ is a differential bigraded
algebra, then $\mathcal{B}_B(M)$ is a right-$B$-module. Note that $A$
and $B$ are not necessarily
augmented.

Similarly, $(_B \mathcal{B}(N), \mathrm d_{_B \mathcal{B}(N)})$,
with $_B
\mathcal{B}(N):=B\underline{\otimes}_B N$, is the $A_{\infty}$-bar
construction
of $N$ as left $A_{\infty}$-$B$-module. It is an
$A_{\infty}$-left-$B$-module. If $B$ is a differential bigraded
algebra, then $\mathcal{B}_B(M)$ is a left-$B$-module.\\
The following lemma is almost tautological, but it is helpful to fix
notation.

\begin{Lem}\label{grillo!}
Let $A$, $C$ be flat $A_{\infty}$-algebras, $B$ be a unital
associative algebra and
$(N,\mathrm d_N)$ be $A_{\infty}$-$B$-$C$-bimodule. If $(M,\mathrm
d_M)$ is an
$A_{\infty}$-$A$-$B$-bimodule such that
$\mathrm{\bar{d}}^{k,l}_M= 0$ if $(k,l)\neq (0,0),(0,1),(k,0)$ and it is
unital
as right $B$-module, then there exists a
strict $A_{\infty}$-$A$-$C$-bimodule isomorphism
\begin{eqnarray}
M\underline{\otimes}_B N \equiv M \otimes_{B} {_B\mathcal{B}(N)}. \label{q1}
\end{eqnarray}
\end{Lem}
The $A_{\infty}$-$A$-$C$-bimodule $M \otimes_{B}
{_B\mathcal{B}(N)}$ in lem.~\ref{grillo!} is given as follows. As bigraded object we have
\begin{eqnarray*}
(M \otimes_{B} {_B\mathcal{B}(N)})^i_j:=\bigoplus_{i_1+i_2=i, \atop
j_1+j_2=j}
M^{i_1}_{j_1}\otimes {_B\mathcal{B}(N)}^{i_2}_{j_2}/ Q^i_j,
\end{eqnarray*}
where $Q^i_j=\bigoplus_{i_1+i_2=i, j_1+j_2=j}Q\cap (M^{i_1}_{j_1}\otimes
{_B\mathcal{B}(N)}^{i_2}_{j_2})$ and $Q$ denotes the submodule in
$M\otimes {_B\mathcal{B}(N)}$ generated by elements of the form $m\cdot
b\otimes \tilde{B}-m\otimes b\cdot \tilde{B}$, with $m\in
M$, $b\in B$
and $\tilde{B}\in {_B\mathcal{B}(N)}$.
$M \otimes_{B} {_B\mathcal{B}(N)}$ is endowed with an
$A_{\infty}$-$A$-$C$-bimodule structure given by a codifferential
$\mathrm{d}_{M \otimes_{B} {_B\mathcal{B}(N)}}$ with Taylor components
\begin{eqnarray}
&&\mathrm{\bar{d}}^{0,0}_{M \otimes_{B} {_B\mathcal{B}(N)}}=-s\circ
\mathrm{D}^{0,0}\circ s,~~~
\mathrm{\bar{d}}^{n,0}_{M \otimes_{B}
{_B\mathcal{B}(N)}},~~~\mathrm{\bar{d}}^{0,m}_{M \otimes_{B}
{_B\mathcal{B}(N)}},
\end{eqnarray}
with $n,m\geq 1$, s.t.
\begin{eqnarray*}
&&\mathcal{D}^{0,0}(m\otimes_B \tilde{B})=
s^{-1}(\mathrm{\bar{d}}_M^{0,0}(sm))\otimes_B
\tilde{B}+(-1)^{|m|}m\otimes_B s^{-1}(\mathrm{\bar{
d}}^{0,0}_{_B\mathcal{B}(N)}(s\tilde{B}))), \\
&&\mathrm{\bar{d}}^{r,0}_{M \otimes_{B}
{_B\mathcal{B}(N)}}(a_1|\dots|a_r|(m\otimes_B \tilde{B}))=s(
s^{-1}(\mathrm{\bar{d}}^{r,0}_M (a_1|\dots | a_r | m )) \otimes_B
\tilde{B} ),\\
&&\mathrm{\bar{d}}^{0,m}_{M \otimes_{B} {_B\mathcal{B}(N)}}((m\otimes_B
(b\otimes b_1|\dots|b_q\otimes
n))|c_1|\dots|c_m)=(-1)^{|m|+|b|+\sum_{i=1}^q(|b_i|-1)}\\
&&\sum_{q'=0}^q s( m\otimes_B (b\otimes b_1|\dots|b_{q'}|\otimes s^{-1}(
\mathrm{\bar{d}}^{q',n}_{N}( b_{q'+1}|\dots|b_{q}|n|c_1|\dots|c_m ) ) ) ) ,
\end{eqnarray*}
and zero otherwise.

\begin{Remark}
Exchanging the role of $M$ and $N$ in lemma~\ref{grillo!} we can
describe the
strict $A_{\infty}$-$A$-$C$-bimodule isomorphism
\begin{eqnarray*}
M\underline{\otimes}_B N \equiv \mathcal{B}_B(M) \otimes_{B} N. \label{q2}
\end{eqnarray*}
\end{Remark}

\begin{Remark}
In what follows we only consider the triple $(A,K,B)$ of bigraded
$A_{\infty}$-objects with $A_{\infty}$-algebras
$(A,\mathrm d_A)$ and $(B,\mathrm d_ B)$ s.t. $A=S(X^*)$, $B=\wedge(X)$
and $A_\infty$-bimodule $(K,\mathrm d_K)$, $K=\mathbb K$.
\end{Remark}

\subsubsection{On the right derived $A_{\infty}$-module $\underline{K}$}
Let $\mathcal{B}_B(K):= K\underline{\otimes}_B B$ denote the bar
construction of $K$ as right $B$-module.
By definition,
\begin{eqnarray*}
\mathcal{B}_B(K)^i_j:= \bigoplus_{q\geq 0} (K\otimes B[1]^{\otimes
q}\otimes B)^i_j.
\end{eqnarray*}
and
\begin{eqnarray*}
\mathcal{B}_B(K)^i_j=\left\{\begin{array}{cc}
0 & i+j >0,\\
K\otimes (B[1]^{\otimes -(i+j)}\otimes B)^i_j & i+j\leq 0.
\end{array}\right.
\end{eqnarray*}
We have also the isomorphism $\mathcal{B}_B(K)\cong M\otimes B$ in $\C$, where $M^i_j=\bigoplus_{q\geq 0} K\otimes (B[1]^{\otimes
q})^i_j=K\otimes (B[1]^{\otimes -(i+j)})^i_j$.

\begin{Def}
The right derived dual module $\underline{K}$ of $K$ is
the object 
\begin{eqnarray}
\underline{K}=\Homm_{B}(\mathcal{B}_B(K), B). \label{fin_dual}
\end{eqnarray}
in $\C$.
\end{Def}
We recall that, for every pair $M,N$ of  right $B$-modules, then
$\Homm_{B}(M,N)$
is the object in $\C$ with bihomogeneous components $\Homm_{B}^{i,j}(M,N)=\{\varphi\in\Hom2^{i,j}_{\C}(M,N),~\varphi~\mbox{right
$B$-linear}\}$.  \begin{Lem}\label{yuppy}
$\underline{K}$ can be endowed with a strictly unital
$A_{\infty}$-$B$-$A$-bimodule structure
$\mathrm d_{\underline{K}}$ with Taylor components given by

\begin{eqnarray*}
&&\bar{\mathrm d}_{\underline{K}}^{0,0}= -s\circ\partial_{
\underline{K} } \circ s^{-1},\\
&&\bar{\mathrm d}_{\underline{K}}^{1,0}(b|\varphi)=s\circ \mathrm
D_{\underline{K}}^{1,0}(b|\varphi),\\
&&\bar{\mathrm d}_{\underline{K}}^{0,m}(\varphi|a_1|\dots|a_m)=
s\circ \mathrm D_{\underline{K}}^{0,m}(\varphi|a_1|\dots|a_m),\\
\end{eqnarray*}
with
\begin{eqnarray*}
&&\partial_{\underline{K}}(\varphi)=(-1)^{|\varphi|}
\varphi\circ \bar{d}^{0,0}_{\mathcal B_B(K)},\\
&&\mathrm D_{\underline{K}}^{1,0}(b|\varphi)(1,b_1|\dots|b_q,b')=
(-1)^{|b|} b\cdot\varphi(1,b_1|\dots|b_q,b'),\\
&&\mathrm
D_{\underline{K}}^{0,m}(\varphi|a_1|\dots|a_m)(1,b_1|\dots|b_q,b')=\\
&&(-1)^{|\varphi|-1+\sum_{i=1}^m(|a_i|-1)}\sum_{q'=0}^q\varphi(s^{-1}\mathrm{\bar{d}}^{m,q'}_K(
a_1|\dots|a_m|1|b_1|\dots|b_{q'}),b_{q'+1}|\dots|b_q,b),
\end{eqnarray*}
and $\mathrm{\bar{d}}^{0,0}_{\mathcal{B}_B(K)}=-s\circ
\bar{d}^{0,0}_{\mathcal{B}_B(K)}\circ s^{-1}$, $\bar{\mathrm d}_{\underline{K}}^{n,m}=0$, otherwise.
\end{Lem}

\begin{Cor}
$(\underline{K},\mathrm d_{\underline{K}})$ is a strictly unital
differential
bigraded left
$B$-module; we have a strict isomorphism
\begin{eqnarray*}
K\underline{\otimes}_B \underline{K}\equiv\mathcal{B}_B(K)\otimes_B
\underline{K} \label{q3}
\end{eqnarray*}
of strictly unital $A_{\infty}$-$A$-$A$-bimodules.
\end{Cor}

\subsubsection{On the quasi-isomorphism $A\rightarrow
K\underline{\otimes}_B \underline{K}$.}

\begin{Def}
$\Endd_{B}(\mathcal{B}_B(K))$ is the
object in $\C$ with bihomogeneous components
\begin{eqnarray}
\Endd_B^{i,j}(\mathcal{B}_B(K))=
\Homm^{0,0}_{B}( K \underline{\otimes}_B B,
(K \underline{\otimes}_B B)[i]\langle j\rangle ).\label{fin_end}
\end{eqnarray}
\end{Def}

\begin{Lem}
$\Endd_{B}(\mathcal{B}_B(K))$ can be endowed with a strictly unital
$A_{\infty}$-$A$-$A$-bimodule structure $\mathrm
d_{\Endd_{B}(\mathcal{B}_B(K))}$
with Taylor components

\begin{eqnarray*}
&&\bar{\mathrm d}_{\Endd_{B}(\mathcal{B}_B(K))}^{0,0}=
-s\circ\partial_{\Endd_{B}(\mathcal{B}_B(K)) } \circ s^{-1},\\
&&\bar{\mathrm
d}_{\Endd_{B}(\mathcal{B}_B(K))}^{l,0}(a_1|\dots|a_l|\varphi)=s\circ \mathrm
D_{\Endd_{B}(\mathcal{B}_B(K))}^{l,0}(a_1|\dots|a_l|\varphi),\\
&&\bar{\mathrm
d}_{\Endd_{B}(\mathcal{B}_B(K))}^{0,m}(\varphi|a_1|\dots|a_m)=
s\circ \mathrm
D_{\Endd_{B}(\mathcal{B}_B(K))}^{0,m}(\varphi|a_1|\dots|a_m),\\
\end{eqnarray*}
with
\begin{eqnarray*}
&&\partial_{\Endd_{B}(\mathcal{B}_B(K))}(\varphi)=(-1)^{|\varphi|}
\varphi\circ
\bar{d}^{0,0}_{\mathcal{B}_B(K)}-\bar{d}^{0,0}_{\mathcal{B}_B(K)}\circ\varphi,\\
&&\mathrm
D_{\Endd_{B}(\mathcal{B}_B(K))}^{l,0}(a_1|\dots|a_l|\varphi)(1,b_1|\dots|b_q,b)=\\
&&(-1)^{\sum_{i=1}^l(|a_i|-1)-1}s^{-1}(\mathrm{\bar{d}}^{l,0}_{\mathcal{B}_B(K)}(a_1|\dots|a_l|\varphi(1,b_1|\dots|b_q,b))),
(l\geq 1)\\
&&\mathrm
D_{\Endd_{B}(\mathcal{B}_B(K))}^{0,m}(\varphi|a_1|\dots|a_m)(1,b_1|\dots|b_q,b)=\\
&&(-1)^{|\varphi|}\varphi(\sum_{q'=0}^qs^{-1}(\mathrm
d^{m,q'}_{K}(a_1|\dots|a_m|1|b_1|\dots|b_{q'})),b_{q'+1}|\dots|b_q,b),\\
\end{eqnarray*}
and $\mathrm{\bar{d}}^{0,0}_{\mathcal{B}_B(K)}=-s\circ
\bar{d}^{0,0}_{\mathcal{B}_B(K)}\circ s^{-1}$, where
$\mathrm{\bar{d}}^{0,0}_{\mathcal{B}_B(K)}$ is given in proposition~\ref{prop1}, 
$\bar{\mathrm d}_{\Endd_{B}(\mathcal{B}_B(K))}^{n,m}=0$, otherwise.
\end{Lem}

Let $(\underline{\End2}_{B}(K),\mathrm
d_{\underline{\End2}_{B}(K)})$ be the strictly unital $A_\infty$-$A$-$A$-bimodule described in prop~\ref{End-bim}.

We recall that the bar resolution
$\mathcal{B}_B(K))=K\underline{\otimes}_B B$ is homotopy equivalent to $K$
in $\C$ (but not as right bigraded $B$-modules); the maps giving such
homotopy equivalence are the projection
$\mathrm p:K\underline{\otimes}_B B\rightarrow K$
and the inclusion $\mathrm i: K\rightarrow K\underline{\otimes}_B B$,
with $\mathrm p(1,b)=b(0)$ and
$\mathrm p(1,b_1|\dots|b_q,b)=0$ for $q\geq 1$.

\begin{Prop}\label{chiave}
$(\Endd_{B}(\mathcal{B}_B(K)),\mathrm
d_{\Endd_{B}(\mathcal{B}_B(K))})$ and
$(\underline{\End2}_{B}(K),\mathrm d_{\underline{\End2}_{B}(K)})$
are homotopy equivalent as strictly unital $A_{\infty}$-$A$-$A$-bimodules.
\end{Prop}
\begin{proof}
We define the strict (and strictly unital) morphism $\mathcal H:
\Endd_{B}(\mathcal{B}_B(K))\rightarrow\underline{\End2}_{B}(K)$
 of strictly unital
$A_{\infty}$-$A$-$A$-bimodules, {\it via} $\mathcal H=s\circ H\circ
s^{-1}$, where,
for any $(i,j)\in \mathbb{Z}^2$,
$H:\Endd_{B}^{i,j}(\mathcal{B}_B(K))\rightarrow
\underline{\End2}_{B}^{i,j}(K)$
is the composition $H:=(s\circ 1\circ
s^{-1})\circ \mathcal{P}\circ\mathcal{I}$,
denoting by $\mathcal I$ and $\mathcal P$ the morphisms
\begin{eqnarray*}
&&\Homm^{0,0}_{B}(\mathcal{B}_B(K),\mathcal{B}_B(K)[i]\langle j\rangle
)\stackrel{\mathcal{I}}{\rightarrow}\Hom2^{0,0}_{\C}(K\otimes
\T(B[1]),\mathcal{B}_B(K)[i]
\langle j\rangle )\\
&&\stackrel{\mathcal P}{\rightarrow}\Hom2^{0,0}_{\C}(K\otimes
\T(B[1]),K[i]\langle j\rangle)
\stackrel{s\circ 1\circ s^{-1}}{\rightarrow}
\Hom2^{0,0}_{\C}(K[1]\otimes \T(B[1]),(K[1])[i]\langle j\rangle),
\end{eqnarray*}
with $\mathcal{I}(\varphi)(1,b_1|\dots|b_q):=\varphi(1, b_1|\dots|b_q,
1)$, $\mathcal{P}(\psi):=\mathrm p\circ\psi$.
More explicitly, if $\varphi\in \Endd_{B}^{i,j}(\mathcal{B}_B(K))$, then
\begin{eqnarray}
(H\varphi)(1|b_1|\dots|b_{q})=s((\varphi^{0}(1,b_1|\dots|b_{q},1))(0)),\label{qq}
\end{eqnarray}
denoting by $\varphi^{0}(1,b_1|\dots|b_{q},1)$ the projection of
$\varphi(1,b_1|\dots|b_{q},1)$ onto $K\otimes B$.

To prove
\[
\mathcal H\circ \mathrm d_{\Endd_{B}} =\mathrm
d_{\underline{\End2}_{B}(K)}\circ \mathcal H
\]
is straightforward; the only issue is represented by the signs; all
details are contained in \cite{Ferrario2}.
\end{proof}

\begin{Prop}\label{prop12}
$(K\underline{\otimes}_B \underline{K},\mathrm d_{K\underline{\otimes}_B
\underline{K}})$ and
$(\Endd_{B}(\mathcal{B}_B(K)),\mathrm
d_{\Endd_{B}(\mathcal{B}_B(K))})$ are strictly isomorphic as strictly unital
$A_{\infty}$-$A$-$A$-bimodules.
\end{Prop}
\begin{proof}
We recall that $\mathcal{B}_B(K)=M\otimes B$ in $\C$.
The strict isomorphism of $A_{\infty}$-$A$-$A$-bimodules
\[
\mathcal G:\mathcal{B}_B(K)\otimes_B \underline{K}\rightarrow
\Endd_{B}(\mathcal{B}_B(K)),
\]
with $\mathcal G=s\circ G\circ s^{-1}$ is given as follows. The
morphism $G$ is defined by the commutative diagram

\[
\begin{CD}
\mathcal B_B(K)\otimes_B \Homm_B(\mathcal B_B(K),B) @> G >>
\Endd_B(\mathcal B_B(K)) \\
@VV \mathcal T_1 V @| \\
( M\otimes B) \otimes_B \Hom2_{\C}(M,B) @. \Endd_B(\mathcal B_B(K)) \\
@VV \mathcal T_2 V @VV\mathcal I V \\
M\otimes \Hom2_{\C}(M,B) @>\underline{G}>> \Hom2_{\C}(M,M\otimes
B) \\
\end{CD}
\]
in $\C$, where
\begin{eqnarray}
\underline{G}(m,\varphi)(m'):=m\otimes\varphi(m'), \label{unG}
\end{eqnarray}
and $\mathcal T_1$, $\mathcal T_2$, $\mathcal I$ denote the obvious
isomorphisms.
Note the sign in $\mathcal T_2((m\otimes
b)\otimes_B\varphi)=(-1)^{|m|+|b|+|\varphi|}m\otimes
b\varphi$.

More explicitly
\begin{eqnarray}
G((m\otimes b)\otimes_B\varphi)(m'\otimes b'):=(-1)^{|m|+|b|+|\varphi|}
m\otimes
b\varphi(m'\otimes b').\label{Depeche!}
\end{eqnarray}
By definition, $G(Q^i_j)=0$ for every $(i,j)\in\mathbb{Z}^2$,
where $Q^i_j$ is the submodule in $(\mathcal B_B(K)\otimes
\Homm_B(\mathcal B_B(K),B))^i_j$
introduced in the proof of lemma~\ref{grillo!}. So $G$ is well defined,
as morphism
in $\C$. Note that $\mathcal T_2(Q^i_j)=0$, as well.
$\underline{G}$ is an isomorphism in $\C$; so $G$ is an isomorphism
in $\C$ as well; in fact
$M$ is an object in $\C$ with finite dimensional bihomogeneous
components $M^i_j=K\otimes (B[1]^{\otimes -(i+j)})^i_j$,
for every $i,j\in \mathbb{Z}^2$.
We finish the proof of proposition~\ref{prop12} by checking that $G$ is
a chain
map and commutes with
the left and right $A_{\infty}$-$A$-actions on
$\mathcal{B}_B(K)\otimes_B \underline{K}$ and
$\Endd_B(\mathcal B_B(K))$.
The only issue is represented by the signs
appearing in $G$ and in the Taylor components of the codifferentials on
$\mathcal{B}_B(K)\otimes_B \underline{K}$ and
$\Endd_B(\mathcal B_B(K))$. In particular, the non trivial sign in \eqref{Depeche!} is necessary to
prove compatibility between $G$ and the right $A_\infty$-module structures, i.e.
\[
     \mathrm{\bar{d}}^{0,m}_{\underline{\End2}_B(K)}(\mathcal G((m\otimes
b))\otimes_B\varphi),a_1|\dots|a_m)=
\mathcal G ( \mathrm{\bar{d}}^{0,m}_{\mathcal B_B(K)\otimes_B
\underline{K}}( ((m\otimes b)\otimes_B\varphi )|a_1|\dots|a_m),
\]
with l.h.s. equal to ( applying it on $m'\otimes b'$ with
$m'\otimes b'=1,b_1|\dots|b_q\otimes b'$)
\begin{eqnarray*}
\sum_{q'=0}^qm\otimes
b\varphi(s^{-1}(\mathrm{\bar{d}}^{m,q'}_{\underline{K}}(a_1|\dots|a_m|1|b_1|\dots|b_{q'})),b_{q'}|\dots|b_q,b'),
\end{eqnarray*}
and r.h.s. equal to
\begin{eqnarray*}
(-1)^{|m|+|b|}\mathcal G ( (m\otimes b)\otimes_B s^{-1}(\mathrm{\bar{d}}^{0,m}_{\underline{K}}(  \varphi |a_1|\dots|a_m)))=
(-1)^{|\varphi|-1+\sum_{i=1}^m(|a_i|-1) }     m\otimes b\cdot s^{-1}(\mathrm{\bar{d}}^{0,m}_{\underline
K}(\varphi|a_1|\dots|a_m)(1,b_1|\dots|b_q,b')).
\end{eqnarray*}
The Taylor components $\mathrm{\bar{d}}^{0,m}_{\underline K}$ generate
the sign $(-1)^{|\varphi|+\sum_{i=1}^m(|a_i|-1)-1}$; so we are done.

\end{proof}

We summarize the results so far into
\begin{Cor}$(K\underline{\otimes}_B \underline{K},\mathrm
d_{K\underline{\otimes}_B \underline{K}})$ and
$(\underline{\End2}_{B}(K),\mathrm d_{\underline{\End2}_{B}(K))}
)$ are
homotopy equivalent as strictly unital $A_{\infty}$-$A$-$A$-bimodules.
\end{Cor}

\begin{Prop}\label{Policy}
There exists a strictly unital quasi-isomorphism
\[
A\rightarrow K\underline{\otimes}_B\underline{K}
\]
of strictly unital $A_\infty$-$A$-$A$-bimodules.
\end{Prop}
\begin{proof}
Just compose the homotopy equivalence in the above corollary with the
left derived action $\mathrm L_A$.
\end{proof}


\subsubsection{On the quasi-isomorphism $B\rightarrow
\underline{K}\underline{\otimes}_A K$.}

Following the example of $\underline{K}$, we can introduce the left
derived bimodule
\[
\overline{K}=\Homm_A(A\underline{\otimes}_A K, A).
\]
As $A[1]$ is concentrated in cohomological degree $-1$, then
\[
A\underline{\otimes}_A K= A\otimes N
\]
in $\C$, with $N^i_j=0$ if $i>0$, $N^0_0=K$ and
\[
N^i_j=\bigoplus_{j_1+\dots+j_{-i}=j}
\overbrace{(A[1])_{j_1}^{-1}\otimes\dots\otimes
(A[1])_{j_{-i}}^{-1}}^{-i-\mbox{times}}
\]
for any $i<0$. Every bihomogeneous component of $N$ is finite
dimensional. In what follows $_A \mathcal{B}(K):=A\underline{\otimes}_A K$.

By definition, $\overline{K}$ is a strictly unital
$A_\infty$-$B$-$A$-bimodule with codifferential
$\mathrm d_{\overline{K}}$ whose Taylor components are given by

\begin{eqnarray*}
\bar{\mathrm d}_{\overline{K}}^{0,0}= -s\circ\partial_{
\overline{K} } \circ s^{-1}, &
\bar{\mathrm d}_{\overline{K}}^{k,0}(b_1|\dots|b_k|\varphi)=s\circ \mathrm
D_{\overline{K}}^{k,0}(b_1|\dots|b_k|\varphi), &
\bar{\mathrm d}_{\overline{K}}^{0,1}(\varphi|a)=
s\circ \mathrm D_{\overline{K}}^{0,1}(\varphi|a),
\end{eqnarray*}
with
\begin{eqnarray*}
&&\partial_{\overline{K}}(\varphi)=(-1)^{|\varphi|}
\varphi\circ \bar{d}^{0,0}_{_A \mathcal{B}(K)},\\
&&\mathrm D_{\overline{K}}^{k,0}(b_1|\dots|b_k|\varphi)(a,a_1|\dots|a_q,1)=
(-1)^{(|\varphi|+|a|+\sum_{i=1}^q(|a_i|-1)+1)(\sum_{i=1}^k(|b_i|-1)+1 ) }\\
&&\sum_{q'=0}^q\varphi(a,a_1|\dots|a_{q-q'},s^{-1}\mathrm{\bar{d}}^{q',k}_K(
a_{q-q'+1}|\dots|a_q|1|b_1|\dots|b_{k}) ),\\
&&\mathrm
D_{\overline{K}}^{0,1}(\varphi|a')(m)=(-1)^{|\varphi|+|a||m|}\varphi(m)\cdot
a,\\
\end{eqnarray*}
where $\mathrm{\bar{d}}^{0,0}_{_A \mathcal{B}(K)}=-s\circ
\bar{d}^{0,0}_{_A \mathcal{B}(K)}\circ s^{-1}$
and $\bar{\mathrm d}_{\overline{K}}^{n,m}=0$, otherwise.

To check that $(\overline{K},\mathrm d_{\overline{K}})$ is a a strictly
unital $A_\infty$-$B$-$A$-bimodule
is long but straightforward.

\begin{Def}
$\Endd_{A}(_A \mathcal{B}(K))^{op}$ is the
object in $\C$ with bihomogeneous components
\begin{eqnarray*}
\Endd_{A}^{i,j}(_A \mathcal{B}(K))^{op}=
\Homm^{0,0}_{A}( A \underline{\otimes}_A K,
(A \underline{\otimes}_A K)[i]\langle j\rangle ).
\end{eqnarray*}
\end{Def}

\begin{Lem}
$\Endd_{A}(_A \mathcal{B}(K))^{op}$ can be endowed with a strictly unital
$A_{\infty}$-$B$-$B$-bimodule structure $\mathrm
d_{\Endd_{A}(_A \mathcal{B}(K))^{op}}$
with Taylor components

\begin{eqnarray*}
&&\bar{\mathrm d}_{\Endd_{A}(_A \mathcal{B}(K))^{op}}^{0,0}=
-s\circ\partial_{\Endd_{A}(_A \mathcal{B}(K))^{op} } \circ s^{-1},\\
&&\bar{\mathrm
d}_{\Endd_{A}(_A
\mathcal{B}(K))^{op}}^{l,0}(b_1|\dots|b_l|\varphi)=s\circ \mathrm
D_{\Endd_{A}(_A \mathcal{B}(K))^{op}}^{l,0}(b_1|\dots|b_l|\varphi),\\
&&\bar{\mathrm
d}_{\Endd_{A}(_A \mathcal{B}(K))^{op}}^{0,m}(\varphi|b_1|\dots|b_m)=
s\circ \mathrm
D_{\Endd_{A}(_A \mathcal{B}(K))^{op}}^{0,m}(\varphi|b_1|\dots|b_m),\\
\end{eqnarray*}
with
\begin{eqnarray*}
&&\partial_{\Endd_{A}(_A \mathcal{B}(K))^{op}}(\varphi)=(-1)^{|\varphi|}
\varphi\circ
\bar{d}^{0,0}_{_A \mathcal{B}(K)}-\bar{d}^{0,0}_{_A
\mathcal{B}(K)}\circ\varphi,\\
&&\mathrm
D_{\Endd_{A}(_A
\mathcal{B}(K))^{op}}^{l,0}(b_1|\dots|b_l|\varphi)(a,a_1|\dots|a_q,1)=
(-1)^{(|\varphi|+|a|+\sum_{i=1}^q(|a_i|-1)+1)(\sum_{i=1}^l(|b_i|-1)+1)}\\
&& \sum_{q'=0}^q \varphi(a,a_1|\dots|a_{q'},s^{-1}\mathrm
d^{q-q',l}_{K}(a_{q'+1}|\dots|a_q|1|b_1|\dots|b_l)),~~~(l\geq 1)\\
&&\mathrm
D_{\Endd_{A}(_A
\mathcal{B}(K))^{op}}^{0,m}(\varphi|b_1|\dots|b_m)(a,a_1|\dots|a_q,1)=\\
&&(-1)^{(|a|+\sum_{i=1}^q(|a_i|-1) )\sum_{i=1}^m(|b_i|-1) }\mathrm
d^{0,m}_{_A
\mathcal{B}(K)}(\varphi(a,a_1|\dots|a_q,1)|b_1|\dots|b_m),~~~(m\geq 1)\\
\end{eqnarray*}
and $\mathrm{\bar{d}}^{0,0}_{_A \mathcal{B}(K)}=-s\circ
\bar{d}^{0,0}_{_A \mathcal{B}(K)}\circ s^{-1}$, where
$\mathrm{\bar{d}}^{0,0}_{_A \mathcal{B}(K)}$ is given in proposition~\ref{prop1},
and $\bar{\mathrm d}_{\Endd_{A}(_A \mathcal{B}(K))^{op}}^{n,m}=0$,
otherwise.
\end{Lem}

 \begin{Prop}
$(\overline{K}\underline{\otimes}_A K,\mathrm
d_{\overline{K}\underline{\otimes}_A
K})$ and
$(\Endd_{A}(_A \mathcal{B}(K))^{op},\mathrm
d_{\Endd_{A}(_A \mathcal{B}(K))^{op}})$ are strictly isomorphic as
strictly unital
$A_{\infty}$-$B$-$B$-bimodules.
\end{Prop}
\begin{proof}
The proof is similar to the one of prop.~\ref{prop12}, with due
changes.
\end{proof}

\begin{Prop}\label{Truth}
There exists a strictly unital quasi-isomorphism
\[
B\rightarrow \overline{K}\underline{\otimes}_A K
\]
of strictly unital $A_\infty$-$B$-$B$-bimodules.
\end{Prop}
\begin{proof}
Just compose the homotopy equivalence in the above prop. with the right
derived action $\mathrm R_B$.
\end{proof}


\subsubsection{$A_{\infty}$-Morita theory for the triple $(A,K,B)$}\label{stripped}

\subsubsection{On the functors}
Let us consider the functors
\begin{eqnarray*}
F': \Modd_{\infty}(A)\rightarrow
\Modd^{strict}_{\infty}(B), & G': \Modd_{\infty}(B) \rightarrow \Modd^{strict}_{\infty}(A),
\end{eqnarray*}
given by
\begin{eqnarray*}
F'(M):= M\underline{\otimes}_{A}K ,~~~~G'(N):=
N\underline{\otimes}_{B}\underline{K} ,
\end{eqnarray*}
on objects $M\in \Modd_{\infty}(A)$ and $N\in \Modd_{\infty}(B)$, while on
morphisms
$f: M_1\rightarrow M_2$ in $\Modd_{\infty}A$ and $g: N_1\rightarrow
N_2$ in $\Modd_{\infty}B$ we set
\begin{eqnarray*}
F'(f):= M_1\underline{\otimes}_{A}K\rightarrow
M_2\underline{\otimes}_{A}K,~~~~ G'(g):=
N_1\underline{\otimes}_{B}\underline{K}\rightarrow
N_2\underline{\otimes}_{B}\underline{K},
\end{eqnarray*}
with
\begin{eqnarray*}
F'(f):=(s^{-1}\circ F\circ s)\otimes 1, ~~~~ G'(g):=(s^{-1}\circ G\circ
s)\otimes 1.
\end{eqnarray*}
We have denoted by
$F:\mathcal R(M_1)\rightarrow \mathcal R(M_2)$,
respectively $G:\mathcal R(N_1)\rightarrow \mathcal R(N_2)$,
the unique lifting of $f$ (resp. $g$) to a $\T(A[1])$-counital-comodule morphism,
respectively a $\T(B[1])$-counital-comodule morphism.
In this notation, $\mathcal R(M_1):=M_1[1]\otimes \T(A[1])$
and similarly for $\mathcal R(N_1)$, with due changes.
Let $\mathcal F$ and $\mathcal G$ be the functors given by the compositions
\begin{eqnarray*}
\mathcal F:\Modd_{\infty}(A)\stackrel{F'}{\rightarrow}
\Modd^{strict}_{\infty} (B) \stackrel{i}{\hookrightarrow} \Modd_{\infty}(B),
\end{eqnarray*}
and
\begin{eqnarray*}
\mathcal G:\Modd_{\infty}(B)\stackrel{G'}{\rightarrow}
\Modd^{strict}_{\infty} (A) \stackrel{i}{\hookrightarrow} \Modd_{\infty}(A),
\end{eqnarray*}
denoting by $i$ the inclusion of the subcategories
$\Modd^{strict}_{\infty} (A) $ (resp. $\Modd^{strict}_{\infty} (B)$)
in $\Modd_{\infty}A$ (resp. $\Modd_{\infty}(B)$). We remark that
$\Modd^{strict}_{\infty} (A) $ and $\Modd^{strict}_{\infty}(B)$ are not full
subcategories.

If two morphisms $f$ and $g$ in $\Modd_{\infty}(A)$ are ($A_\infty$-)
homotopic,
then we write $f\sim g$. An analogous notation holds true in
$\Modd_{\infty}(B)$.
If the homotopy between $f$ and $g$ is strict, then we write
$f\sim_{strict} g$.

\begin{Lem}
$\bf{a)}$ Let $f\sim g$ in $\Modd_{\infty}A$, resp. in
$\Modd_{\infty}B$. Then
\begin{eqnarray*}
F'(f)\sim_{strict} F'(g),
\end{eqnarray*}
in $\Modd^{strict}_{\infty}B$, resp.
\begin{eqnarray*}
G'(f)\sim_{strict} G'(g),
\end{eqnarray*}
in $\Modd^{strict}_{\infty}A$.

$\bf{b)}$ The functors $F'$ and $G'$ send strictly unital homotopy
equivalences to
strict and strictly unital homotopy equivalences.

\end{Lem}

\begin{proof} Part $\bf{a)}$. Let $f,g: M\rightarrow N$ with $f\sim g$
in $\Modd_{\infty}(A)$,
i.e. $f-g=\mathrm d_Nh+h\mathrm d_M$, where $h: M\rightarrow N$ is a strictly unital
$A_{\infty}$-homotopy.
By definition, $h$ is a degree $-1$ map with components
$h_n: M[1]\otimes A[1]^{\otimes n}\rightarrow N[1]$, $n\geq 0$.
We claim that $H: M\underline{\otimes}_A K\rightarrow
N\underline{\otimes} K$, where
\begin{eqnarray*}
H_0:=(s^{-1}\circ h \circ s) \otimes 1
\end{eqnarray*}
and $H_n=0$ for $n>0$, is a strict $A_{\infty}$-homotopy between $F'(f)$
and $F'(g)$, i.e.
\begin{eqnarray}
F'(f)- F'(g)=\mathrm d_{N\underline{\otimes}_A K}\circ H+H\circ \mathrm
d_{M\underline{\otimes}_A K}. \label{h1}
\end{eqnarray}
Eq. (\ref{h1}) is equivalent to
\begin{eqnarray}
(F'(f)- F'(g))(s(m,a_1|\dots|a_q,1))=(\mathrm d_{N\underline{\otimes}_A
K}\circ H+H\circ \mathrm d_{M\underline{\otimes}_A K})(s(m,a_1|\dots|a_q,1)), \label{h2}
\end{eqnarray}
and
\begin{eqnarray}
0=(\mathrm d_{N\underline{\otimes}_A K}\circ H+H\circ \mathrm
d_{M\underline{\otimes}_A K})((m,a_1|\dots|a_q,1)|b_1|\dots|b_l), \label{h3}
\end{eqnarray}
for every $q,l\geq 0$.
Let us consider at first eq. (\ref{h2});
on the l.h.s. we have terms involving the Taylor components of the
codifferential $\mathrm d_M$ on $M$ and the $A_{\infty}$-homotopy $h$ by
the homotopy hypothesis $f\sim g$;
all we need to prove is that on the r.h.s the terms involving the
Taylor components of the codifferential $\mathrm d_K$ on $K$ cancel.
This is true because these terms appear in
\begin{eqnarray*}
&&\sum_{q_1=0}^q\sum_{q_2=0}^{q_1+1}(-1)^{1+|m|+\sum_{i=1}^{q-q_2}(|a_i|-1)}((h_{q_1}(m,a_1|\dots,a_{q_1})|
a_{q_1+1}|\dots | a_{q-q_2},\mathrm d_{K}^{q_2}(a_{q-q_2+1}| \dots |a_q|1 ) )+\\
&&\sum_{q_1=0}^q\sum_{q_2=0}^{q_1+1}(-1)^{|m|+(|a_i|-1)}
(h_{q_1}(m,a_1|\dots| a_{q_1}),a_{q_1+1}|\dots | a_{q-q_2}, \mathrm d_{K}^{q_2}(a_{q-q_2+1}| \dots |a_q|1 ) )=0,
\end{eqnarray*}
as the $A_{\infty}$-homotopy $h$ has degree $-1$. Eq. $(\ref{h3})$ is
equivalent to
\begin{eqnarray*}
&&0=\sum_{q_1=0}^q \mathrm d_{N\underline{\otimes}_A K}
(s^{-1}(h_{q_1}(m|a_1|\dots | a_{q_1}))|a_{q_1+1}|\dots | a_{q}|1
)+\nonumber \\
&&\sum_{q_2=0}^q
(-1)^{1+|m|+\sum_{i=1}^{q-q_2}(|a_i|-1)}H(m|a_1|\dots|a_{q-q_2},\mathrm
d_{K}^{q_2,l}(a_{q-q_2+1}|\dots|a_q|1|b_1|\dots|b_l ) ),
\end{eqnarray*}
which is verified by the same argument we used for eq.
(\ref{h2}) and (\ref{h3}). The case $f\sim g$ in $\Modd_{\infty}(B)$ is
similar.

Part $\bf{b)}$. The morphism $f:M\rightarrow N$ is a homotopy equivalence in
$\Modd_{\infty}A$ if there exists a morphism $g: N\rightarrow M$ in
$\Modd_{\infty}(A)$ s.t.
$f \circ g \sim 1$ and $ g\circ f\sim 1$. We denote by $h_1:
N\rightarrow N$, resp. $h_2: M\rightarrow M$ the $A_{\infty}$-homotopies
between $f \circ g$ and $1_N$, resp. $g\circ f$ and $1_M$. We want to
prove that
\begin{eqnarray*}
&&F'(f)\circ F'(g)=1+ \mathrm
d_{N\underline{\otimes}_A K}\circ H_1+H_1\circ\mathrm d_{N\underline{\otimes}_A K},\\
&&F'(g)\circ F'(f)=1+ \mathrm
d_{M\underline{\otimes}_A K}\circ H_2+H_2\circ\mathrm d_{M\underline{\otimes}_A K},
\end{eqnarray*}
with strict $A_{\infty}$-homotopies
\begin{eqnarray*}
H_i:=(s^{-1}\circ h_i \circ s) \otimes 1,
\end{eqnarray*}
for $i=1,2$. Using the proof of $\bf{a)}$ we get the statement. The case
for $G'$ is similar.
\end{proof}


\subsubsection{On the derived categories}

In this section we introduce the derived categories
$\DD^{\infty}(A)$, respectively $\DD^{\infty}(B)$, of right unital
$A_{\infty}$-modules over $A$, respectively $B$, with strictly unital $A_\infty$-morphisms. 
Using the theory of closed model categories it is possible to prove
\begin{Thm}[K.~Lefevre-Hasegawa,~\cite{Kenji} ]
Let $A$ be an augmented $A_\infty$-algebra
\footnote{Augmentation w.r.t. a ground field $\mathbb K$ of characteristic $0$.};
quasi-isomorphisms in $\Modd_{\infty}(A)$ are homotopy equivalences of strictly unital $A_\infty$-$A$-modules.
\end{Thm}
This results implies that
\[
\DD^{\infty}(A)=\Modd_{\infty}(A)/\sim,
\]
and similarly for $\DD^{\infty}(B)$. In this setting  quasi-isomorphisms of strictly unital $A_\infty$-modules
are already isomorphisms in the homotopy categories; no localization is needed. The main advantage is represented
by the explicit structure of the morphisms in the derived categories themselves; no ``roofs'' manipulation is needed.

 We discuss now the triangulated structures on the derived categories.
 The direct sum of two objects in $\Modd_{\infty}(A)$
is again a strictly unital $A_{\infty}$-module; the cohomological
grading shift functor $\Sigma(M):=M[1]$ in actually an endofunctor
on $\DD^{\infty}(A)$ and $\DD^{\infty}(B)$. It follows that
$\Sigma(\cdot):=\cdot[1]$ is an autoequivalence of
$\DD^{\infty}(A)$ and $\DD^{\infty}(B)$. More precisely, let $(M,d_M)$ be an object of $\DD^{\infty}(A)$. The bigraded object
$M[1]$ can be endowed with a strictly unital
$A_{\infty}$-$A$ module structure as follows.
The codifferential $\mathrm d_{M[1]}$ has Taylor components
$\mathrm{\bar{d}}_{M[1]}^{l}:
(M[1])[1]\otimes B[1]^{\otimes l}\rightarrow (M[1])[1]$
given by
\begin{eqnarray*}
\mathrm{\bar{d}}_{M[1]}^{l}=-s\circ \mathrm{\bar{d}}_{M}^{l}\circ (
s^{-1}\otimes 1 ).
\end{eqnarray*}
Proving that $\mathrm d_{M[1]}^2=0$ is a straightforward
sign-check. Given any morphism
$F: M[1]\rightarrow N[1]$ in $\DD^{\infty}(A)$ with Taylor
components (of bidegree (0,0))
$\bar{F}^{l}:M[1]\otimes B[1]^{\otimes l}\rightarrow N[1]$, we get the
induced morphism
$\tilde{F}: (M[1])[1] \rightarrow (N[1])[1]$ in $\Modd_{\infty}(A)$
with Taylor components
\begin{eqnarray*}
\tilde{\bar{F}}^{l}=s\circ \bar{F}^{l}\circ ( s^{-1}\otimes 1 ).
\end{eqnarray*}
Once again, the proof of $\tilde{F}\circ
\mathrm d_{M[1]}=\mathrm d_{N[1]}\circ\tilde{F}$ is a straighforward
sign check. Same considerations hold in $\DD^{\infty}(B)$. The inverse functor $\Sigma^{-1}$ is given by $\Sigma^{-1}(\cdot)=\cdot[-1]$.

\begin{Def}[\cite{Kenji}]
The triangulated structure  on the derived category
$\DD^{\infty}(A)$ is given as follows.
The autoequivalence $\Sigma$ is simply the (cohomological) grading shift
functor $\Sigma=[1]$.
The distinguished triangles are isomorphic to those induced by
semi-split sequences of strict $A_{\infty}$-morphisms
\begin{eqnarray*}
M\stackrel{f}{\rightarrow} M'\stackrel{g}{\rightarrow} M''
\end{eqnarray*}
in $\Modd_{\infty}A$, i.e. sequences such that
\begin{eqnarray}
0\rightarrow M\stackrel{f}{\rightarrow} M'\stackrel{g}{\rightarrow}
M''\rightarrow 0 \label{seq}
\end{eqnarray}
is an exact sequence in $\C$, and such that there exists a splitting
$\rho\in\Hom2_{\C}(M', M)$ of $f$ with
\begin{eqnarray*}
\rho\circ \bar{\mathrm d}^{i}_M=\bar{\mathrm d}^i_{M'}\circ (\rho\otimes
1^{\otimes i-1}), ~~~~~i\geq 2.
\end{eqnarray*}
\end{Def}
For the derived category of $B$ the definition is analogous.
The splitting $\rho$ in the exact sequence (\ref{seq})
does not commute with the differentials $\bar{\mathrm d}_M^0$ and
$\bar{\mathrm d}^0_{M'}$, in general. The above exact triangles endow
$\DD^{\infty}(A)$ with a triangulated category structure; the proof is
contained
in thm. 2.4.3.1 in \cite{Kenji}; the idea is induce the triangulated
category structure on $\DD^{\infty}(A)$ by using the one on $\DD(UA)$, denoting by $UA$ the enveloping algebra of $A$; by definition $UA$ is 
a differential (bi)graded algebra we refer to \cite{Kenji}, \cite{Halperin} for all details; its derived category $\DD(UA)$ is a well-known object. The equivalence of categories
$\DD(UA)\rightarrow \DD^{\infty}(A)$ becomes then an
equivalence of triangulated categories.

Let $X\rightarrow Y\rightarrow Z\rightarrow  X[1]$ be a distinguished triangle in $\DD^{\infty}(A)$; it is isomorphic
 to a triangle of the form
$M\stackrel{f}{\rightarrow} M'\stackrel{g}{\rightarrow}
M''\rightarrow M[1]$,
with $M\stackrel{f}{\rightarrow} M'\stackrel{g}{\rightarrow} M''$
satisfying the hypothesis of the above definition. In more detail, let
\begin{eqnarray*}
0\rightarrow M\stackrel{f}{\rightarrow}
M'\stackrel{g}{\rightarrow} M''\rightarrow 0
\end{eqnarray*} be a semi-split exact sequence with $f$, $g$ strict, and splitting
$\rho:M'\rightarrow M$, $\rho\circ f=1$.
This implies that
\[
{M'}^i_j\cong M^i_j\oplus {M''}^i_j
\]
as vector spaces over $\mathbb K$, for any $(i,j)\in\mathbb Z$; in virtue of this
we assume that $M'=(M\oplus M'',\mathrm d_{M\oplus M''})$,
where $\mathrm d_{M\oplus M''}=(\mathrm d_{M}-h,\mathrm d_{M''})$.
It follows that $\mathrm d_{M\oplus M''}\circ \mathrm d_{M\oplus M''}=0$
if and only if
$h:M''\rightarrow M[1]$ defines an $A_\infty$-morphism of strictly
unital $A_\infty$-$A$-modules. Thanks to this, we will consider the semisplit exact sequence
$0\rightarrow M\stackrel{i}{\rightarrow} M'\stackrel{p}{\rightarrow} M''\rightarrow 0$
with $i$ and $p$ the natural inclusion and projection (which are strict morphisms in $\DD^{\infty}(A)$),
and complete it  to the exact triangle
\begin{eqnarray}
M\stackrel{i}{\rightarrow} M'\stackrel{p}{\rightarrow}
M''\stackrel{h}{\rightarrow}  M[1].\label{ciao!}
\end{eqnarray}
A small $memento$; in section~\ref{trianghbar} we will discuss the triangulated structure on some ``deformed'' derived categories
of topologically free modules; some examples will be given:
taking there the ``limit''  $\hbar=0$ we  obtain  further examples of exact triangles in $\DD^{\infty}(A)$ and $\DD^{\infty}(B)$.

\subsubsection{On the functors $\mathcal F$ and $\mathcal G$}
Collecting the results on the derived categories of $A$ and $B$ and the definitions of the functors $\mathcal F$
 and $\mathcal G$ we arrive at the pair of functors
\begin{eqnarray*}
\mathcal F:\DD^{\infty}(A)~ \stackrel{F'}{\rightarrow}
\Modd^{strict}_{\infty} (B) / \sim_{strict}
\stackrel{i}{\hookrightarrow}
\DD^{\infty}(B),
\end{eqnarray*}
and
\begin{eqnarray*}
\mathcal G:\DD^{\infty}(B)~ \stackrel{G'}{\rightarrow}
\Modd^{strict}_{\infty} (A) / \sim_{strict}
\stackrel{i}{\hookrightarrow}
\DD^{\infty}(A),
\end{eqnarray*}
with a little abuse of notation.
\begin{Prop}\label{aah!}
Let $(\mathcal F,\mathcal G)$ be the pair of functors introduced above.
Then $\mathcal F(A) \simeq K$, $\mathcal F(\overline{K})\simeq B$,
in $\DD^{\infty}(B)$, and $\mathcal G(B)\simeq \underline{K}$, $\mathcal G(K)\simeq A$
in $\DD^{\infty}(A)$.

It follows that
\begin{eqnarray*}
\mathcal F(\mathcal G(K)) \simeq K ~in~ \DD^{\infty} (B), ~~\mathcal
G(\mathcal F(A)) \simeq A ~in ~\DD^{\infty}(A).
\end{eqnarray*}

\end{Prop}

\begin{proof}
The quasi-isomorphisms of strictly unital $A_\infty$-bimodules
\begin{eqnarray*}
&&K\rightarrow A\underline{\otimes}_A
K\rightarrow ( K\underline{\otimes}_B
\underline{K})\underline{\otimes}_A K
=\mathcal F(\mathcal G(K))
\end{eqnarray*}
and
\begin{eqnarray*}
A\rightarrow K\underline{\otimes}_B \underline{K} \rightarrow
A\underline{\otimes}_A (K
\underline{\otimes}_B\underline{K}) \equiv (A\underline{\otimes}_A K)
\underline{\otimes}_B\underline{K} = \mathcal G(\mathcal F(A))
\end{eqnarray*}
give both the statements. We used lem.~\ref{Lem2}, prop.~\ref{Policy} and prop.~\ref{Truth}.
\end{proof}

\begin{Lem}\label{hell}
$(\mathcal F,\varphi_1)$ and $(\mathcal G,\varphi_2)$ are exact functors w.r.t the triangulated category
 structures on $\DD^{\infty}(A)$ and $\DD^{\infty}(B)$; for any  $M\in \DD^{\infty}(A)$:
\[
\varphi_1(\mathcal F(M)): (M\underline{\otimes}_A K)[1]\rightarrow
M[1]\underline{\otimes}_A
K,~~~~\varphi_1(\mathcal F(M))(s(m,a_1|\dots|a_l,k)):=(m|a_1|\dots|a_l,k),
\]
and similarly for $\varphi_2$.
\end{Lem}

\begin{proof}
$\mathcal F$ and $\mathcal G$ send quasi-isomorphisms into quasi-isomorphisms as
quasi-isomorphisms in the derived categories
$\DD^{\infty}(A)$ and $\DD^{\infty}(B)$ are homotopy equivalences. To
prove that $\mathcal F$ (and $\mathcal G$)
are exact w.r.t.\ the triangulated structures on the derived categories
 it is sufficient to consider  triangles of the form \eqref{ciao!}, i.e.
 $M\stackrel{i}{\rightarrow}M\oplus M'\stackrel{p}{\rightarrow}M'\stackrel{h}{\rightarrow}M[1]$.

Applying $\mathcal{F}$ to such a triangle, and using the above lemmata we
get the sequence
\begin{eqnarray*}
M\underline{\otimes}_{A}\stackrel{\mathcal F(i)}{\rightarrow}
(M\oplus M')\underline{\otimes}_{A}K\stackrel{\mathcal F(p)}{\rightarrow} M'\underline{\otimes}_{A}K
\stackrel{\varphi_{M}(F(M))}{\rightarrow}  (M\underline{\otimes}_{A}K)[1]
\end{eqnarray*}
in $\DD^{\infty}(B)$; the short exact sequence ($\mathcal F$ is additive)
\begin{eqnarray}
0\rightarrow M\underline{\otimes}_{A}K\stackrel{F(i)}{\rightarrow
}(M\oplus M')\underline{\otimes}_{A}K
\stackrel{F(p)}{\rightarrow} M''\underline{\otimes}_{A}K
\rightarrow 0 \label{tac}
\end{eqnarray}
 is semi-split w.r.t. the splitting
\[
F(\rho):=\rho\otimes 1,
\]
denoting by $\rho: M'\rightarrow M$ the splitting of the short exact
sequence $0\rightarrow M\stackrel{i}{\rightarrow} M\oplus M'\stackrel{p}{\rightarrow} M'\rightarrow 0$.
In fact
\[
F(\rho)\circ F(\alpha)=1,~~ \mbox{and}~~F(\rho) \circ(s^{-1}\circ
\mathrm{\bar{d}}^{0,i}_{M'\underline{\otimes}_A K})=
(s^{-1}\circ \mathrm{\bar{d}}^{0,i}_{M\underline{\otimes}_A K})
\circ ( F(\rho)\otimes 1^{\otimes i}),
\]
for $i\geq 1$. Then (\ref{tac}) can be completed to the distinguished triangle
\begin{eqnarray*}
M\underline{\otimes}_{A}K\stackrel{F(i)}{\rightarrow}
(M\oplus M')\underline{\otimes}_{A}K\stackrel{F(p)}{\rightarrow} M'\underline{\otimes}_{A}K\stackrel{h'}{\rightarrow}
(M\underline{\otimes}_{A}K)[1],
\end{eqnarray*}
with $h':=F(h)$. In summary $\mathcal F$ sends exact triangles into exact triangles. Same considerations
holds true for $\mathcal G$.
\end{proof}

With $ \triang^{\infty}_A(M)$ we denote the full triangulated subcategory
in $\DD^{\infty}(A)$ generated by $\{M[i]\langle j\rangle, i\in\mathbb{Z}\}$.
$ \thick^{\infty}_A(M_A)$, resp. $ \thick^{\infty}_A(N_B)$ are the thick subcategories of direct
summands of objects in $ \triang^{\infty}_A(M_A)$,
resp. $ \triang^{\infty}_B(N_B)$. We refer to Appendix C for all definitions.
 Finally, we can state the main theorem of this section.
\begin{Thm}\label{Thm29}
Let $X$ be a finite dimensional vector space on $\mathbb K=\mathbb R$,
or $\mathbb C$. Let $(A,K,B)$ be the triple of $A_{\infty}$-structures
with $A=S(X^*)$ and $B=\wedge(X)$ Koszul dual augmented differential
bigraded algebras with zero differential and $K=\mathbb K$ endowed with
the bigraded
$A_{\infty}$-$A$-$B$-bimodule structure $\mathrm d_K$ given in \cite{CFFR}.
The triangulated functor
\begin{eqnarray*}
  \mathcal F : \DD^{\infty}(A)\rightarrow
\DD^{\infty}(    B), ~~~~~\mathcal
F(\bullet)=\bullet~\underline{\otimes}_{A} K
\end{eqnarray*}
induces the equivalence of triangulated categories
\begin{eqnarray*}
\triang^{\infty}_{A}(A)\simeq
\triang^{\infty}_{B}(K), &
\thick^{\infty}_{A}(A)\simeq\thick^{\infty}_{B}(K).
\end{eqnarray*}
Let $(\tilde{K},\mathrm d_{\tilde{K}})$ be the
$A_{\infty}$-$B$-$A$-bimodule
with $\tilde{K}=K$ and $\mathrm d_{\tilde{K}}$ obtained from $\mathrm
d_K$ exchanging $A$ and $B$; then the triangulated functor
\begin{eqnarray*}
  \mathcal F^{''} : \DD^{\infty}(B)\rightarrow
\DD^{\infty}(    A), ~~~~~\mathcal
F^{''}(\bullet)=\bullet~\underline{\otimes}_{B}
\tilde{K}
\end{eqnarray*}
induces the equivalence of triangulated categories
\begin{eqnarray*}
\triang^{\infty}_{A}(\tilde{K})\simeq
\triang^{\infty}_{B}(B), &
\thick^{\infty}_{A}(\tilde{K})\simeq\thick^{\infty}_{B}(B).
\end{eqnarray*}
\end{Thm}
\begin{proof}
Appendix B. 
\end{proof}

\section{Deformation Quantization of $A_\infty$-structures}

In this section we study the quantizations
$(A_\hbar,K_\hbar,B_\hbar)$ of the
$A_{\infty}$-structures on the triple $(A,K,B)$. In this contest,
the term ``quantization'', or more properly, ``Deformation
Quantization'' refers to a technique that produces new $A_{\infty}$-structures
from already given $A_{\infty}$-data: the latter are recovered from the
former through a ``limiting'' procedure.  For the original idea we refer to \cite{Q}.
 $A_{\infty}$-structures on bigraded topologically free $\mathbb K\c1$-modules are said to be topological.
The deformations are obtained through certain Feynman diagrams
expansions, a ``two
branes'' Formality theorem and an explicit choice of an $\hbar$-formal
quadratic Poisson
bivector $\pi_{\hbar}=\hbar\pi$ on $X$, the
finite dimensional vector space underlying $A$ and $B$. For the full
construction and the 2-branes formality theorem we refer to
\cite{CFFR}; the diagrammatic techniques there described generalize those introduced in \cite{Kont}. 
The choice of a quadratic Poisson bivector field is motivated by the
necessity of preserving
the internal grading on the Deformation Quantization of triple $(A,K,B)$;
its main consequences are
\begin{itemize}
\item The Deformation Quantizations $(A_\hbar,B_\hbar)$ of $(A,B)$ are
flat bigraded $A_{\infty}$-algebras.

\item The Deformation Quantization $K_\hbar$ of $K$ is a left
$A_\hbar$-module and a right $B_\hbar$-module with zero differential.
\item It is possible to quantize the bimodules $A\underline{\otimes}_A
K$, $K\underline{\otimes}_B B$, $\underline{K}$,
$\underline{\End2}_A(K)$ and $\underline{\End2}_B(K)$
straightforwardly by using the ``classical'' $A_{\infty}$-bimodule
structures
with due changes.
\item The quantized left and right derived actions are
quasi-isomorphisms of topological $A_{\infty}$-algebras $and$
topological $A_{\infty}$-bimodules.

\end{itemize}

\subsubsection{On modules over $\mathbb K\c1$}
We consider the local ring $\mathbb K\c1$ of formal power series with coefficients in $\mathbb K$.
Topological free $\mathbb K\c1$ modules are
modules over $\mathbb K\c1$ isomorphic to $\mathbb K\c1$-modules of
the form $M\c1$,
with $M$ a $\mathbb K$ vector space.
Let $M$ and $N$ be $\mathbb K\c1$-modules. The $\mathbb K\c1$-module
$M\otimes_{\mathbb K\c1}N$ is the quotient of the
tensor product $M\otimes N$ ($\otimes=\otimes_{\mathbb K}$) by the
subspace generated by all elements of the form
$km\otimes n - m\otimes kn$, with $k\in\mathbb K\c1$ and $m\in M$, $n\in N$.
We denote by $\tilde{\otimes}$ the completed tensor product $M\tilde{\otimes}N$ of $M\otimes_{\mathbb
K\c1}N$. If $M$ and $N$ are topologically free, i.e. $M=M_1\c1$ and $N=N_1\c1$, then
$M\c1\tilde{\otimes}N\c1$ is topologically free as well;
in fact $M\tilde{\otimes}N= (M_1\otimes N_1)\c1$.

Let $\Hom2_{\mathbb K\c1}(M\c1,N\c1)$ be the space of $\mathbb
K\c1$-linear morphisms from $M\c1$ to $N\c1$;
there exists an isomorphism $\mathcal I: \Hom2(M,N)\c1\rightarrow \Hom2_{\mathbb K\c1}(M\c1,N\c1)$
of $\mathbb K\c1$-modules. Any $\varphi\in\Hom2_{\mathbb
K\c1}(M\c1,N\c1)$ is uniquely determined
by a formal power series
\[
\sum_{i\geq 0}\varphi_i\hbar^i \in \Hom2(M,N)\c1.
\]
We observe that any $\varphi\in\Hom2_{\mathbb
K\c1}(M\c1,N\c1)$ is continuous w.r.t.\ the
$\hbar$-adic topology on $M\c1$ and $N\c1$. In the sequel we will use the formal power series description of morphisms extensively.

\subsubsection{On the  category $\GD$}
Let $\GD$ be the category of bigraded $\mathbb K\c1$-modules;
an object in $\GD$ is a collection $\{M^i_j\}_{i,j\in\mathbb Z}$ of $\mathbb
K\c1$-modules; the space of morphisms $\Hom2_{\GD}(M,N)$ is the
 object in $\GD$ with bihomogeneous components
\[
\Hom2^{i,j}_{\GD}(M,N)=\prod_{r,s\in\mathbb Z}\Hom2_{\mathbb
K\c1}(M^r_s,N^{i+r}_{j+s} ).
\]


\subsubsection{Topologically free modules in $\GD$}

We say that an object $M_{\hbar}$ in $\GD$ is topologically free if
\begin{eqnarray*}
M_{\hbar}=\{(M_{\hbar})^i_j\}_{(i,j)\in\mathbb{Z}^2}, ~~~\mbox{with}~~~  (M_{\hbar})^i_j=M^i_j\c1.
\end{eqnarray*}
Let $M_\hbar$ and $N_\hbar$ be topologically free objects in $\GD$, with  $M_\hbar=M\c1$
and $N_\hbar=N\c1$, for $M,N$
objects in $\C$; then
$\Hom2_{\GD}(M_\hbar,N_\hbar)$ is the topologically free object in $\GD$
with bihomogeneous components
\[
\Hom2^{i,j}_{\GD}(M_\hbar,N_\hbar)=
\Hom2^{i,j}_{\C}(M,N)\c1.
\]
For any topologically free $M_{\hbar}$ in $\GD$, the objects $M_{\hbar}[k]$ and $M_{\hbar}\langle l\rangle$ in $\GD$
are defined $via$
\begin{eqnarray*}
M_{\hbar}[k]=\{(M_{\hbar}[k])^i_j\}_{(i,j)\in\mathbb{Z}^2},~~~~~(M_{\hbar}[k])^i_j:=M^{i+k}_j\c1;
\end{eqnarray*}
and
\begin{eqnarray*}
M_{\hbar}\langle l\rangle=\{(M_{\hbar}\langle
l\rangle)^i_j\}_{(i,j)\in\mathbb{Z}^2},~~~~~
(M_{\hbar}\langle l\rangle)^i_j:=M^{i}_{j+l}\c1;
\end{eqnarray*}
for any $(k,l)\in\mathbb Z^2$. Topologically free objects in $\GD$ form a full subcategory in $\GD$
which is not abelian;
we endow it with a monoidal structure induced by the completion
$\tilde{\otimes}$, w.r.t the
$\hbar$-adic topology, of the tensor product of topologically free
$\mathbb K\c1$-modules.
More precisely, for any $M_{\hbar}$ and $N_{\hbar}$ topologically free
in $\GD$ and $(i,j)\in\mathbb Z^2$,
we write (with a little abuse of notation)
\[
(M_{\hbar}\tilde{\otimes}N_{\hbar})^i_j=\bigoplus_{i_1+i_2=i,\atop
j_1+j_2=j}
{M^{i_1}_{j_1}}\c1\tilde{\otimes}{N^{i_2}_{j_2}}\c1,
\]
where $\tilde{\otimes}$ on the right hand side is the completed tensor
product of topologically free
$\mathbb K\c1$-modules introduced above.

\section{Topological $A_{\infty}$-structures}\label{quantized-A-infty}

\subsubsection{Topological $A_{\infty}$-algebras}

\begin{Def}
 Let $A_\hbar$ be a topologically free object in $\GD$. The topological tensor coalgebra over $A_\hbar$ is the triple 
$(\T(A_{\hbar}[1]), \Delta_{\hbar}, \epsilon_\hbar)$ where
\[
\T(A_{\hbar}[1]):=\bigoplus_{q\geq 0} A_{\hbar}[1]^{\tilde{\otimes} q}=\T(A[1])\c1
\]
in $\GD$, and
\[
\Delta_{\hbar}\in\Hom2_{\GD}^{0,0}(\T(A_{\hbar}[1]),\T(A_{\hbar}[1])\tilde{\otimes}\T(A_{\hbar}[1]))
\]
given by $\Delta_{\hbar}=\sum_{i\geq
0}\Delta^{(i)}\hbar^i=\Delta^(0)=\Delta$, where $\Delta$ denotes the coproduct on $\T(A[1])$ and
$\epsilon_\hbar=\epsilon$, where $\epsilon$ is the counit in $\T(A[1])$.
\end{Def}
By definition $(1\tilde{\otimes}\Delta_{\hbar})\circ\Delta_\hbar=(\Delta_{\hbar}\tilde{\otimes}1)\circ\Delta_{\hbar}$
 and $(\epsilon_\hbar\tilde{\otimes} 1)\circ\Delta_\hbar= (1\tilde{\otimes}\epsilon_\hbar )\circ\Delta_\hbar= 1$.
\subsubsection{On codifferentials: definitions}

\begin{Def}
A coderivation on  $\T(A_\hbar[1])$ is a morphism $\mathrm d_{A_{\hbar}}
\in\Hom2^{1,0}_{\GD}(\T(A_{\hbar}[1]),\T(A_{\hbar}[1]))$ s. t.
$(1\tilde{\otimes}\mathrm d_{A_{\hbar}}+\mathrm
d_{A_{\hbar}}\tilde{\otimes} 1)\circ\Delta_{\hbar}=
\Delta_{\hbar}\circ\mathrm d_{A_{\hbar}}$
and
\begin{equation}
\mathrm
d_{A_{\hbar}}^2=0\label{braciola}
\end{equation}
\end{Def}
Let $\mathrm d_{A_{\hbar}}$ be the coderivation on $\T(A_{\hbar}[1])$
uniquely determined by the formal power series
\[
\mathrm d_{A_{\hbar}}=\sum_{i\geq 0}\mathrm d^{(i)}_{A_{\hbar}}\hbar^i, ~~
\mathrm d^{(i)}_{A_{\hbar}}\in\Hom2^{1,0}_{\C}(\T(A[1]),\T(A[1])).
\]
Then, by definition of $\mathrm d_{A_\hbar}$, each $\mathrm d^{(i)}_{A_{\hbar}}$ is uniquely determined by the family of Taylor components
 $\mathrm d^{(i),k}_{A_\hbar}=p_{A[1]}\circ \mathrm d^{(i)}_{A_{\hbar}}|_ {A[1]^{\otimes k}}$. The quadratic relations \eqref{braciola} are 
equivalent to a tower of quadratic relations with the Taylor components $\mathrm d^{k,(i)}_{A_\hbar}$, $k\geq 0$, $i\geq 0$.
\begin{Def}
Let $A_{\hbar}$ be a topologically free object in $\GD$. A topological
$A_{\infty}$-algebra structure on $A_{\hbar}$ is the datum of a
coderivation on the topological tensor coalgebra over $A_\hbar$.
\end{Def}
\begin{Lem} Let $( A_{\hbar},\mathrm d_{A_\hbar} )$ be a topological
$A_{\infty}$-algebra. Then
$(A,\mathrm d_A)$, $\mathrm d_A:=\mathrm d^{(0)}_{A_\hbar}$, is
an $A_{\infty}$-algebra (on $\mathbb K$).
\end{Lem}

\subsubsection{Topologically free $A_{\infty}$-modules}

\begin{Def}
Let $M_{\hbar}$ be topologically free module in $\GD$;
$\mathcal R_{\hbar}(M_{\hbar})$ is the object
\[
\mathcal
R_{\hbar}(M_{\hbar}):=M_{\hbar}\tilde{\otimes}\T(A_{\hbar}[1])=
(M[1]\otimes\T(A[1]))\c1
\]
in $\GD$. A right $(\T(A_{\hbar}[1]),\Delta_{\hbar},\epsilon_\hbar)$-counital-comodule structure on
$\mathcal R_{\hbar}(M_{\hbar})$ is the morphism
\[
\Delta^{R}_{\hbar}\in\Hom2_{\GD}^{0,0}(\mathcal
R_{\hbar}(M_{\hbar}),\mathcal R_{\hbar}(M_{\hbar})
\tilde{\otimes}\T(A_{\hbar}[1])),~~~\Delta^{R}_\hbar=\Delta^{R,(0)}_\hbar=\Delta^{R},
\]
satisfying
$(1\tilde{\otimes}\Delta_{\hbar})\circ
\Delta^{R}_{\hbar}=(\Delta^{R}_{\hbar}\tilde{\otimes}1)\circ\Delta^{R}_{\hbar}$ and $(1\tilde{\otimes}\epsilon_\hbar)\circ\Delta^{R}_\hbar=1$,
denoting by
$\Delta^R$ the usual counital-$\T(A[1])$-comodule structure on $M[1]\otimes
\T(A[1])$. 
\end{Def}
\begin{Def}
A codifferential on the right $\T(A_{\hbar}[1])$-comodule
$\mathcal R_{\hbar}(M_{\hbar})$ is
a morphism $\mathrm{d}_{M_{\hbar}}\in\Hom2_{\GD}^{1,0}(\mathcal R_{\hbar}(M_{\hbar}),
\mathcal R_{\hbar}(M_{\hbar}) )$ s.t.
$\Delta^R_{\hbar}\circ\mathrm d_{M_{\hbar}}=(1\tilde{\otimes}\mathrm
d_{M_{\hbar}}+\mathrm d_{A_{\hbar}}
\tilde{\otimes}1)\circ \Delta^R_{\hbar}$
and
\begin{eqnarray}
\mathrm{d}_{M_{\hbar}}^2=0. \label{ee}
\end{eqnarray}
\end{Def}
By definition, if $\mathrm d_{M_{\hbar}}=\sum_{i\geq 0}\mathrm
d^{(i)}_{M_{\hbar}}\hbar^i$, then  each $\mathrm d^{(i)}_{M_{\hbar}}\in \Hom2^{1,0}_{\C}(M[1]\otimes\T(A[1]),M[1]\otimes\T(A[1]))$
is uniquely determined by its Taylor components $\mathrm d^{(i),n}_{M_\hbar}=p_{M[1]}\circ \mathrm
d^{(i)}_{M_{\hbar}}|_ {M[1]\otimes A[1]^{\otimes n}}$, for any $i,n\geq 0$. The quadratic relations \eqref{ee} are equivalent to
a tower of quadratic relations involving the aforementioned maps $\mathrm d^{(i),n}_{M_\hbar}$.
\begin{Def}
Let $M_{\hbar}$ be an object in $\GD$. A topological right
$A_{\infty}$-$A_{\hbar}$-module structure on
$M_{\hbar}$ is the datum of a codifferential $\mathrm {d}_{M_{\hbar}}$
on $\mathcal R_{\hbar}(M_{\hbar})$.
\end{Def}
\begin{Lem}
Let $M_{\hbar}$ be a topological right $A_{\infty}$-$A_{\hbar}$-module. Then
$M$ is a right $A_{\infty}$-$A$-module.
\end{Lem}
In the same spirit, one can define topological left $A_{\infty}$-modules
and topological $A_{\infty}$-bimodules, with due changes.

\subsubsection{On morphisms, quasi-isomorphisms and homotopy equivalences}
\begin{Def}
Let $(M_{\hbar},\mathrm d_{M_{\hbar}})$ and $(N_{\hbar},\mathrm
d_{N_{\hbar}})$ be topological $A_{\infty}$-$A_{\hbar}$-modules,
with $(A_{\hbar},\mathrm d_{A_{\hbar}})$ topological $A_{\infty}$-algebra.

A morphism $f_{\hbar}:M_{\hbar}\rightarrow N_{\hbar}$ of topological
$A_{\infty}$-$A_{\hbar}$-modules
is a map $f_{\hbar}\in\Hom2_{\GD}^{0,0}(\mathcal R_{\hbar}(M_{\hbar}),\mathcal
R_{\hbar}(N_{\hbar}) )$
which is a morphism of $\T(A_{\hbar}[1])$-counital-comodules s.t.
\begin{eqnarray*}
\mathrm d_{N_{\hbar}}\circ f_{\hbar}=f_{\hbar}\circ\mathrm
d_{M_{\hbar}}.
\end{eqnarray*}
\end{Def}
Such a morphism is uniquely determined by a formal power series
$f_{\hbar}=\sum_{i\geq 0}f^{(i)}_{\hbar}\hbar^i$,
with $f^{(i)}_{\hbar}\in\Hom2^{0,0}_{\C}(M[1]\otimes \T(A[1]),N[1]\otimes
\T(A[1]))$, for any $i\geq 0$.
Each component $f^{(i)}_\hbar$ is a morphism of
counital-$\T(A[1])$-comodules, and so it admits an explicit description by Taylor components $f^{(i),n}: M[1]\otimes A[1]^{\otimes n}\rightarrow N[1]$,
for any $i,n\geq 0$.

\begin{Lem}
Let $f_{\hbar}:M_{\hbar}\rightarrow N_{\hbar}$, $f_{\hbar}=\sum_{i\geq
0}f^{(i)}_{\hbar}\hbar^i$
be a morphism of topological $A_{\infty}$-$A_{\hbar}$-modules.
Then $f^{(0)}_{\hbar}: M\rightarrow N$ is a morphism of
$A_{\infty}$-$A$-modules.
\end{Lem}

\begin{Def}
Let $f_\hbar,g_\hbar:M_\hbar\rightarrow N_\hbar$ be morphisms of
topological $A_\infty$-$A_\hbar$-modules; we say that
they are topological $A_\infty$-homotopy equivalent (alternatively: top.
$A_\infty$-homotopic) if there exists a topological
$A_\infty$-homotopy between them, i.e. a map $H_\hbar:
M_\hbar\rightarrow N_\hbar$ of $\T(A_\hbar[1])$-counital-comodules  with
\begin{eqnarray*}
&&H_\hbar=\sum_{i\geq 0} H^{(i)}_\hbar \hbar^{i}, \\
&&H^{(i),n}_\hbar\in \Hom2^{-1,0}_{\C}(M[1]\otimes A[1]^{\otimes
n},M[1]), ~~n\geq 0,
\end{eqnarray*}
such that
\[
f_\hbar-g_\hbar=\mathrm d_{N_\hbar}\circ H_\hbar+H_\hbar\circ \mathrm
d_{M_\hbar}
\]
holds true, order by order in $\hbar$.
\end{Def}

\subsubsection{On units}

Let $A_{\hbar}$ in $\GD$ be a topological $A_{\infty}$-algebra with
codifferential $\mathrm d_{A_{\hbar}}$.
We say that the right $A_{\infty}$-$A_{\hbar}$-module structure
$\mathrm d_{M_{\hbar}}$ on $M_{\hbar}$ is strictly unital if
\[
\mathrm{d}^{(i),n}_{M_{\hbar}}(m|a_1|\dots|\eta|\dots|a_n)=0,
\]
for any $n\geq 2$ and $i\geq 0$.. We have denoted by $\eta$ the unit in
$A$ and by
$ \mathrm{d}^{n}_{M_{\hbar}}$ the $n$-th Taylor component of $\mathrm
d^{(i)}_{M_{\hbar}}$.

A morphism $f_{\hbar}:M_{\hbar}\rightarrow N_{\hbar}$ of topological
$A_{\infty}$-$A_{\hbar}$-modules is strictly unital if
\[
f_{\hbar}^{(i),n}(m|a_1|\dots|\eta|\dots|a_n)=0,
\]
for any $n\geq 1$, $i\geq 0$, where $f_{\hbar}^{(i),n}$ is the $n$-th Taylor
component
of $f_{\hbar}^{(i)}$.

Strictly unital homotopies are defined similarly.


\subsubsection{Quantizing $(A,K,B)$ {\it via} quadratic Poisson structures}

By $X$ we denote a finite dimensional vector space of dimension $n$ on
$\mathbb
K=\mathbb R$ or $\mathbb C$.
Let $(T_{poly}(X)\c1,[\cdot,\cdot]_{\hbar})$ be the trivial deformation
of $(T_{poly}(X),[\cdot,\cdot])$, with Schouten-Nijenhuis bracket
$[\cdot,\cdot]_{\hbar}$ obtained by extending $[\cdot,\cdot]$
$\mathbb K\c1$-linearly. Let $\{x_i\}_{i\in I}$ be a set of
global coordinates on $X$, with $\sharp I=n$.
We say that the Poisson bivector $\pi\in T_{poly}(X)$ is quadratic if it
can be written as
\[
\pi=\sum_{i,j=1}^n \pi^{ij} \frac{\partial}{\partial
x_i}\wedge\frac{\partial}{\partial x_j},~~~\pi^{ij}= \sum_{k,l=1}^n
c^{ij}_{kl} x_kx_l,
\]
for some constant coefficients $c^{ij}_{kl}\in\mathbb K$ such that
$c^{ij}_{kl}=-c^{ji}_{kl}$, for any
$k,l\in I$.
In \cite{CFFR} a 2-branes Formality theorem is proved; we refer to
\cite{CFFR} for all details; here we
sketch the construction in the special case in which the triple
$(A,K,B)$ appears.

To the $A_\infty$-triple $(A,K,B)$ it is possible to associate a unital
$A_\infty$-category $\Cat_\infty
(A,K,B)$; its objects are the branes $U=X$ and $V=\{0\}$ in the vector
space $X$ and the spaces of morphisms
are given by $\Hom2(U,U)=A$, $\Hom2(V,V)=B$, $\Hom2(U,V)=K$ and
$\Hom2(V,U)=0$. In this ``local'' setting the branes are linear
subspaces on the ambient space $X$.
The unital $A_\infty$-category structure on $\Cat_\infty
(A,K,B)$ is induced by the associative algebra structures on $A$, $B$ and
the $A_\infty$-$A$-$B$-bimodule structure on $K$.
The 2-branes Formality theorem states the existence of a
quasi-isomorphism of $L_\infty$-algebras
\[
\mathcal U:(\T^{\bullet+1}_{poly}(X),[\cdot,\cdot],0)\rightarrow
(C^{\bullet+1}(\Cat_\infty
(A,K,B)),[\cdot,\cdot]_G,\partial )
\]
between the differential graded Lie algebra (shortly, DGLA)
of polynomial polyvector fields on $X$ and the DGLA of
Hochschild cochains on the $A_{\infty}$-category $\Cat_\infty(A,K,B)$
endowed with the Gerstenhaber bracket $[\cdot,\cdot]_G$ and Hochschild
differential $\partial$. As a graded object $C^{\bullet+1}(\Cat_\infty
(A,K,B))$ decomposes in the direct sum of three components:
$C^{\bullet+1}(A,A)$, $C^{\bullet+1}(B,B)$ and $C^{\bullet+1}(A,K,B)$.
$C^{\bullet+1}(A,A)$ and $C^{\bullet+1}(B,B)$ are the DGLAs of
Hochschild cochains of $A$ and $B$; they are sub complexes of
$C^{\bullet+1}(\Cat_\infty(A,K,B))$. $C^{\bullet+1}(A,K,B)$ is given by
\[
C^{n}(A,K,B)=\bigoplus_{p+q+r=n-1}\Hom2^q(A^{\otimes p}\otimes K\otimes
B^{\otimes r},K).
\]

The proof of the 2-branes Formality theorem is based on Stokes' theorem
on manifolds with corners and the properties of the 4-color propagators (\cite{CF0},\cite{CFFR},\cite{AF})
at the boundary components. In the general case, we have to consider
short loops in the Feynman diagrams describing the
$L_\infty$-quasi-isomorphism $\mathcal U$.

The $L_\infty$-quasi-isomorphism $\mathcal U$ induces an isomorphism
between the sets of Maurer-Cartan elements (MCEs) on the DGLAs
$(\T^{\bullet+1}_{poly}(X),[\cdot,\cdot],0)$ and $(C^{\bullet+1}(\Cat_\infty
(A,K,B)),[,\cdot,\cdot]_G,\partial)$. MCEs in $\T_{poly}(X)$ are
Poisson structures on $X$; they are mapped to MCEs on
$C^{\bullet+1}(A,A)$ and
$C^{\bullet+1}(B,B)$ which are $A_{\infty}$-deformations of the graded
associative algebra structures on $A$ and $B$ and to an
$A_{\infty}$-deformation of the $A_{\infty}$-$A$-$B$-bimodule structure
on $K$.

Let $\hbar\pi$ be a MCE in $\T_{poly}(X)\c1$; it satisfies
\[
[\hbar\pi,\hbar\pi]_\hbar=0
\]
order by order in $\hbar$, denoting by $[\cdot,\cdot]_\hbar$ the Lie
bracket on $\T_{poly}(X)\c1$
obtained by extending $\mathbb K\c1$-linearly the Lie bracket
$[\cdot,\cdot]$ on $\T_{poly}(X)$.

Let  $\mathcal U(\hbar\pi)=\mathcal U_A(\hbar\pi)+\mathcal
U_B(\hbar\pi)+\mathcal U_K(\hbar\pi)$ be the MCE in
$C^{\bullet+1}(\Cat_\infty
(A,K,B))\c1$ where
$\mathcal U_A(\hbar\pi)$ is the component of $\mathcal U$ on
$C^{\bullet+1}(A,A)\c1$,
$\mathcal U_B(\hbar\pi)$ is the one on $C^{\bullet+1}(B,B)\c1$ and
$\mathcal U_K(\hbar\pi)$ on
$C^{\bullet+1}(A,K,B)\c1$. The three components are
defined through an expansion in Feynman graphs in which ``areal
vertices'' appear.

\begin{Prop}[\cite{CFFR}, section 8.1]

Let $\hbar\pi$ be an $\hbar$-formal quadratic MCE in $\T_{poly}(X)\c1$.

\begin{itemize}
\item The $A_\infty$-deformations of $A$, resp. $B$ are given by
$(A\c1,\cdot+\mathcal U_A(\hbar\pi)),~~~(B\c1,\wedge+\mathcal U_B(\hbar\pi))$.

In other words, $A$ and $B$ are deformed into bigraded associative
algebras with zero differential.
The deformed products preserves the internal grading.
\item The $A_{\infty}$-$A\c1$-$B\c1$-bimodule deforming
$K$ is given by $(K\c1,\mathrm d_{K_{\hbar}}),~~~ \mathrm d_{K_{\hbar}}=\mathrm
d_K+\mathcal U_K(\hbar\pi)$.

The codifferential $\mathrm d_{K_{\hbar}}$ is such that
\begin{eqnarray*}
\mathrm{d}^{(i),n,0}_{K_{\hbar}} =\mathrm{d}^{(i),0,m}_{K_{\hbar}}=0
\end{eqnarray*}
if either $m=n=0$ or $m,n,\geq 2$, for any $i\geq 0$.
\end{itemize}
\end{Prop}
Choosing a general Poisson structure $\pi$ on $X$ we obtain different
quantizations $A_\hbar$ resp. $B_\hbar$ of $A$, resp. $B$,
in general curved as $A_{\infty}$-algebras. For a curved example we
refer to \cite{CFR}.

\subsubsection{Quantizing bimodules}

In section 2 we have defined left and right bar resolutions of
$A_{\infty}$-bimodules.

The Taylor components of the codifferential on such resolutions are
given by
the formul\ae~\eqref{eq-tayl-tens}. We quantize the resolutions
considering the triple $(A_\hbar,K_\hbar,B_\hbar)$.
\begin{Def2}
Let $(M_\hbar,\mathrm d_{M_\hbar})$ be a topological
$A_{\infty}$-$A_\hbar$-$B_\hbar$ bimodule.
The left topological bar resolution of $M_\hbar$ is the object
$(A_\hbar\underline{\tilde{\otimes}}_{A_\hbar}M_\hbar =(A_\hbar\tilde{\otimes}
\T(A_\hbar[1])\tilde{\otimes}
M_\hbar)$ in $\GD$.
 It is a topological $A_\infty
$-$A_\hbar$-$B_\hbar$-bimodule with codifferential $\mathrm
d_{A_\hbar\underline{\tilde{\otimes}}_{A_\hbar}M_\hbar}=\sum_{i\geq 0}
\mathrm d^{(i)}_{A_\hbar\underline{\tilde{\otimes}}_{A_\hbar}M_\hbar}
\hbar^i$. For any $i\geq 0$,
the $(k,l)$-th Taylor component
\[
\mathrm
d^{(i),k,l}_{A_\hbar\underline{\tilde{\otimes}}_{A_\hbar}M_\hbar}\in\Hom2^{1,0}_{\C}((A[1]^{\otimes
k}\otimes
(A\underline{\otimes}_{A}M)[1])\otimes B[1]^{\otimes
l},(A\underline{\otimes}_{A}M)[1])
\]
of $\mathrm
d^{(i)}_{A_\hbar\underline{\tilde{\otimes}}_{A_\hbar}M_\hbar}$
is
given by the formul\ae~\eqref{eq-tayl-tens} with the insertion of
the operators $\mathrm d^{(i),2}_{A_\hbar}$, $\mathrm
d^{(i),\cdot,\cdot}_{M_\hbar}$ and
$\mathrm d^{(i),2}_{B_\hbar}$.
\end{Def2}
\begin{proof}
It is easy but quite long to check that 
$\mathrm d_{A_\hbar\underline{\tilde{\otimes}}_{A_\hbar}M_\hbar}\circ \mathrm d_{A_\hbar\underline{\tilde{\otimes}}_{A_\hbar}M_\hbar}=0$
follows from associativity of   the products on $A_\hbar$, $B_\hbar$ and the quadratic relations
$\mathrm d^2_{M_\hbar}=0$.
\end{proof}

We can define right topological bar resolutions, or bar resolutions of
topological $A_{\infty}$-bimodules, with due changes.
In the sequel we will consider the bimodule $K_\hbar$ and the
topological bar resolutions $A_\hbar\underline{\tilde{\otimes}}_{A_\hbar}K_\hbar~~\mbox{and}~~
K_\hbar\underline{\tilde{\otimes}}_{B_\hbar}B_\hbar$.

Let $\underline{K}$ and $\overline{K}$ be the
$A_{\infty}$-$B$-$A$-bimodules introduced in section~\ref{homee}.

\begin{Def2}
The quantization $\underline{K}_\hbar$ of the
$A_{\infty}$-$B$-$A$-bimodule $\underline{K}$
is the object
\[
\underline{K}_\hbar=\Homm_{B_\hbar}(K_\hbar\underline{\tilde{\otimes}}_{B_\hbar}B_\hbar,B_\hbar)
\]
 in $\GD$. It is a strictly unital topological $A_{\infty}$-$B_\hbar$-$A_\hbar$-bimodule
with codifferential $\mathrm d_{\underline{K}_\hbar}=\sum_{i\geq 0}
\mathrm d^{(i)}_{\underline{K}_\hbar} \hbar^i$. For any $i\geq 0$,
the $(k,l)$-th Taylor component
\[
\mathrm
d^{(i),k,l}_{\underline{K}_\hbar}\in\Hom2^{1,0}_{\C}((A[1]^{\otimes
k}\otimes
\underline{K}[1]\otimes B[1]^{\otimes
l},\underline{K}[1])
\]
of $\mathrm
d^{(i)}_{\underline{K}_\hbar}$
is given by the formul\ae~ in lemma~\ref{yuppy}, with the insertion of
the operators $\mathrm d^{(i),2}_{A_\hbar}$, $\mathrm
d^{(i),\cdot,\cdot}_{K_\hbar}$ and
$\mathrm d^{(i),2}_{B_\hbar}$.
\end{Def2}

\begin{proof}
The proof of the topological $A_\infty$-bimodule structure is similar to the one for
$A_\hbar\underline{\tilde{\otimes}}_{A_\hbar}M_\hbar$; we use the
associativity of the
products on $A_\hbar$, $B_\hbar$ and the topological
$A_{\infty}$-bimodule structure on $K_\hbar$.
\end{proof}
The definition of $\overline{K}_\hbar$ is analogous, with due changes.
Similarly, we can introduce the quantizations
$\underline{\End2}_{B_\hbar}(K_\hbar)$, resp.
$\underline{\End2}_{A_\hbar}(K_\hbar)$ of $\underline{\End2}_{B}(K)$,
resp. $\underline{\End2}_{A}(K)$; their
topological $A_{\infty}$-algebra structures are induced by the
topological $A_{\infty}$-structures on $(A_\hbar,K_\hbar,B_\hbar)$.
We use the same classical formul\ae~ introduced in section 6, with due
changes. Such quantizations are topologically free objects in $\GD$.
Let
$\mathcal{B}_{A_\hbar}(K_\hbar)=A_\hbar\underline{\tilde{\otimes}}_{A_\hbar}K_\hbar$,
$\mathcal{B}_{B_\hbar}(K_\hbar)=K_\hbar\underline{\tilde{\otimes}}_{B_\hbar}B_\hbar$
and $\Endd_{A_\hbar}(\mathcal{B}_{A_\hbar}(K_\hbar))$,
$\Endd_{B_\hbar}(\mathcal{B}_{B_\hbar}(K_\hbar))$
be the topologically free objects in $\GD$
\begin{eqnarray*}
&&\Endd_{B_\hbar}(\mathcal{B}_{B_\hbar}(K_\hbar))=\{\varphi\in
\End2_{\GD}(\mathcal{B}_{B_\hbar}(K_\hbar)), B_{\hbar}-linear\},\\
&&\Endd_{A_\hbar}(\mathcal{B}_{A_\hbar}(K_\hbar))=\{\varphi\in
\End2_{\GD}(\mathcal{B}_{A_\hbar}(K_\hbar)), A_{\hbar}-linear\}
\end{eqnarray*}
They are canonically endowed with topological $A_\infty$-bimodule
structures; the formul\ae~are induced by the
classical constructions presented in the previous sections. The constructions used
to quantize $\underline{\End2}_A(K)$ and
$\underline{\End2}_B(K)$ are replied here, with due changes.

\begin{Prop}
The quantized derived actions
\begin{eqnarray*}
&&\mathrm L_{A_\hbar}: A_\hbar\rightarrow
\underline{\End2}_{B_\hbar}(K_\hbar), \\
&&\mathrm R_{B_\hbar}: B_\hbar\rightarrow
\underline{\End2}_{A_\hbar}(K_\hbar)^{op},\\
\end{eqnarray*}
are quasi-isomorphisms of topological $A_{\infty}$-algebras.
\end{Prop}
\begin{proof}
The proposition is proved in \cite{CFFR}, section 8.1. 
\end{proof}

\begin{Cor}
$\mathrm L_{A_\hbar}$ and $\mathrm R_{B_\hbar}$ descend to
quasi-isomorphisms of topological $A_{\infty}$-bimodules.
\end{Cor}
\begin{proof}
The bimodule structures on $A_\hbar$, $B_\hbar$,
$\underline{\End2}_{B_\hbar}(K_\hbar)$ and
$\underline{\End2}_{A_\hbar}(K_\hbar)$ are those described in section 6
for $A$, $B$,
$\underline{\End2}_{B}(K)$ and $\underline{\End2}_{A}(K)$ with due changes.
\end{proof}

\subsubsection{Some  quasi-isomorphisms of quantized bimodules}

\begin{Prop}\label{Shuffle}
\begin{itemize}
\item There exist quasi-isomorphisms
\begin{eqnarray*}
\mu_{K_\hbar}:K_\hbar\underline{\tilde{\otimes}}_{B_\hbar}B_\hbar\rightarrow
K_\hbar, &\mu'_{K_\hbar}:A_\hbar\underline{\tilde{\otimes}}_{A_\hbar}K_\hbar\rightarrow
K_\hbar, \\
\Phi_{K_\hbar}:K_\hbar\rightarrow
K_\hbar\underline{\tilde{\otimes}}_{B_\hbar}B_\hbar
,& \Phi'_{K_\hbar}: K_\hbar\rightarrow
A_\hbar\underline{\tilde{\otimes}}_{A_\hbar}K_\hbar.
\end{eqnarray*}
of strictly unital topological $A_{\infty}$-$A_\hbar$-$B_\hbar$-bimodules.
\item There exists an isomorphism
\[
\Theta^1_\hbar:
K_\hbar\underline{\tilde{\otimes}}_{B_\hbar}\underline{K}_\hbar\rightarrow
\Endd_{B_\hbar}(\mathcal{B}_{B_\hbar}(K_\hbar)),
\]
of strictly unital topological $A_{\infty}$-$A_\hbar$-$A_\hbar$-bimodules.

\item There exists an isomorphism
\begin{eqnarray*}
\Theta^2_\hbar:\overline{K}_\hbar\underline{\tilde{\otimes}}_{A_\hbar}K_\hbar\rightarrow
\Endd_{A_\hbar}(\mathcal{B}_{A_\hbar}(K_\hbar))^{op}
\end{eqnarray*}
of strictly unital topological $A_{\infty}$-$B_\hbar$-$B_\hbar$-bimodules.
\end{itemize}
\end{Prop}
\begin{proof}
On $\mu_{K_\hbar}$. The morphism is defined using the
formul\ae~ for the morphism $\mu$ in proposition 2, section 2, with due
changes. So the compatibility with the topological
$A_{\infty}$-bimodule structures follows.
In other words,
\[
\mu_{K_\hbar}(\sum_{i\geq 0} a_i \hbar^i )=\sum_{n\geq
0}\sum_{i+j=n}\mu^{(i)}_1(a_j)\hbar^{n},
\]
with $\mu^{(i)}_1\in\Hom2^{0,0}_{\C}(\T(A[1])\otimes
(A\underline{\otimes}_{A}K)[1]\otimes \T(B[1]),
K[1])$ uniquely determined by the Taylor components
$\mu^{(i),k,l}_1$, with
\[
\mu^{(i),k,l}_1(a_1|\dots|a_k,s(1|b_1|\dots|b_q|b),b'_1|\dots|b'_l)=\pm\mathrm
d^{(i),k,q+1+l}_{K_\hbar}
(a_1|\dots|a_k|1|b_1|\dots|b_q|b|b'_1|\dots|b'_l),
\]
 The relations $\mathrm d_{K_\hbar}\circ \mu_{K_\hbar}=\mu_{K_\hbar}\circ \mathrm
d_{K_\hbar\underline{\tilde{\otimes}}_{B_\hbar}B_\hbar}$
are equivalent to $\sum_{i+j=n}\mathrm d^{(j)}_{K_\hbar}\circ\mathrm
d^{(i)}_{K_\hbar}=0$,
for any $n\geq 0$. We recall that in the classical case the relations
\[
\mathrm d_{K}\circ \mu=\mu\circ \mathrm d_{K\underline{\otimes}_{B}B}
\]
are equivalent to $\mathrm d_K^2=0$.

It is also clear that
${\mu_{K_\hbar}}|_{\hbar=0}={\mu_{K_\hbar}}|_{\hbar=0}=\mu$; as $\mu$ is a
quasi-isomorphism, then $\mu_{K_\hbar}$ is a quasi-isomorphism as well.
Similar considerations hold for $\mu'_{K_\hbar}$.
For $\Phi_{K_\hbar}$ and $\Phi'_{K_\hbar}$ the conclusions are similar,
with due changes: all we need is to consider
the trivially ``quantized'' version of the formul\ae~introduced in
lemma~\ref{Lem2}.
$\Theta^1_\hbar$ is given by the formal power series
$\Theta^1_\hbar=\sum_{i\geq 0}\Theta^{1,(i)}_\hbar \hbar^i$,
with $\Theta^{1,(i)}_\hbar=0$ for $i\geq 1$ and
$\Theta^{1,(0)}_\hbar=\mathcal I\circ G$,
where $\mathcal I: K\underline{\otimes}_B \underline{K}\rightarrow
\mathcal B_B(K)\otimes \underline{K}$
is an isomorphism of $A_\infty$-$A$-$A$-bimodules and $G:\mathcal B_B(K)\otimes \underline{K}\rightarrow
\Endd_{A}(\mathcal{B}_{B}(K))$
is the isomorphism $A_{\infty}$-$A$-$A$-bimodules given in prop.~\ref{prop12}.

On any element $\sum_{i\geq 0 }v_i\hbar^i$ in
$(\T(A[1])\otimes(K\underline{\otimes}_{B}\underline{K})[1]
\otimes \T(A[1]) )\c1$ the relations
\[
\mathrm d_{\Endd_{B_\hbar}(\mathcal{B}_{B_\hbar}(K_\hbar))}\circ
\Theta^1_\hbar=\Theta^1_\hbar\circ
\mathrm d_{K_\hbar\underline{\tilde{\otimes}}_{B_\hbar}\underline{K}_\hbar}
\]
are equivalent to
\[
\sum_{i+j=n}\mathrm
d^{(j)}_{\Endd_{B_\hbar}(\mathcal{B}_{B_\hbar}(K_\hbar))}(\Theta^{1,(0)}_\hbar(v_i))=\\
\sum_{i+j=n}\Theta_\hbar^{1,(0)}(\mathrm
d^{(j)}_{K_\hbar\underline{\tilde{\otimes}}_{B_\hbar}\underline{K}_\hbar}(v_i))
\]
for any $n\geq 0$; the above relations are verified (projecting both
sides onto
$\Endd_{A_\hbar}(\mathcal{B}_{B_\hbar}(K_\hbar))$, as usual )
as $\Theta^{1,(0)}_\hbar$ commutes with $\mathrm
d^{(j)}_{\Endd_{B_\hbar}(\mathcal{B}_{B_\hbar}(K_\hbar))}$,
$\mathrm
d^{(j)}_{K_\hbar\underline{\tilde{\otimes}}_{B_\hbar}\underline{K}_\hbar}$,
for any $j\geq 0$.

We recall that the Taylor components $\mathrm
d^{(j),\cdot,\cdot}_{\Endd_{B_\hbar}(\mathcal{B}_{B_\hbar}(K_\hbar))}$
and $\mathrm
d^{(j),\cdot,\cdot}_{K_\hbar\underline{\tilde{\otimes}}_{B_\hbar}\underline{K}_\hbar}$
are
given by the Taylor components
$\mathrm d^{\cdot,\cdot}_{\Endd_{B}(\mathcal{B}_{B}(K))}$ and $\mathrm
d^{\cdot,\cdot}_{ K\underline{\otimes}_{B}\underline{K} }$
{\it via} the substitutions
$\mathrm d^2_{A}\mapsto \mathrm d^{(j),2}_{A_\hbar}$,
$\mathrm d^2_{B}\mapsto \mathrm d^{(j),2}_{B_\hbar}$ and
$\mathrm d^{\cdot,\cdot}_{K_\hbar}\mapsto\mathrm
d^{(j),\cdot,\cdot}_{K_\hbar}$.

$\Theta^{1,(0)}_\hbar$ is an isomorphism in $\C$; it follows that
$\Theta^1_\hbar$ is an isomorphism in $\GD$. For $\Theta^2_\hbar$
similar considerations hold, with due changes.

\end{proof}

\begin{Cor}\label{miaumiau}
\begin{itemize}
\item There exists a homotopy equivalence
\[
K_\hbar\underline{\tilde{\otimes}}_{A_\hbar}\underline{K}_\hbar\rightarrow
\underline{\End2}_{B_\hbar}(K_\hbar),
\]
of strictly unital topological $A_{\infty}$-$A_\hbar$-$A_\hbar$-modules.
\item There exists a homotopy equivalence
\[
\overline{K}_\hbar
\underline{\tilde{\otimes}}_{A_\hbar}K_\hbar\rightarrow
\underline{\End2}_{A_\hbar}(K_\hbar)^{op}
\]
of strictly unital topological $A_{\infty}$-$B_\hbar$-$B_\hbar$-modules.
\end{itemize}
\end{Cor}
\begin{proof}
The classical homotopy equivalence $H:
\Endd_{B}(\mathcal{B}_{B}(K))\rightarrow \underline{\End2}_{B}(K)$
defined in prop.~\ref{chiave} induces
the homotopy equivalence $H_\hbar:
\Endd_{B_\hbar}(\mathcal{B}_{B_\hbar}(K_\hbar))\rightarrow
\underline{\End2}_{B_\hbar}(K_\hbar)$,
$H_\hbar=H$. The check is immediate; in fact the codifferentials on
$\Endd_{B_\hbar}(\mathcal{B}_{B_\hbar}(K_\hbar)$ and
$\underline{\End2}_{B_\hbar}(K_\hbar)$ are constructed by using $\mathrm
d^{(i),k,l}_{K_\hbar}$,
$\mathrm d^{(i),j}_{A_\hbar}$ and $\mathrm d^{(i),j}_{B_\hbar}$. It is
necessary to prove the
compatibility of the quantized homotopy equivalence with these
operators. But this goes on like in the classical case.
Similar considerations hold for the second statement, with due changes.
\end{proof}
Composing with the quantized derived action we arrive at
\begin{Cor}\label{party2}
\begin{itemize}
\item There exists a quasi-isomorphism
\[
A\rightarrow
K_\hbar\underline{\tilde{\otimes}}_{A_\hbar}\underline{K}_\hbar,
\]
of strictly unital topological $A_{\infty}$-$A_\hbar$-$A_\hbar$-modules.
\item There exists a quasi-isomorphism
\[
B_\hbar\rightarrow \overline{K}_\hbar
\underline{\tilde{\otimes}}_{A_\hbar}K_\hbar.
\]
of strictly unital topological $A_{\infty}$-$B_\hbar$-$B_\hbar$-modules.
\end{itemize}
\end{Cor}

\subsubsection{On the categories $\Modd^{\infty}_{tf}(A_{\hbar})$ and
$\Modd^{\infty}_{tf}(B_{\hbar})$}

\begin{Def}
Let $(A_{\hbar},K_{\hbar},B_{\hbar})$ be the triple quantizing $(A,K,B)$
w.r.t.\ an $\hbar$-formal
quadratic Poisson bivector $\pi_{\hbar}$.
\begin{itemize}
\item $\Modd^{\infty}_{tf}(A_{\hbar})$ is the category of strictly unital
topological $A_{\infty}$-right-$A_{\hbar}$-modules.
\item $\Modd^{\infty}_{tf}(B_{\hbar})$ is the category of strictly unital
topological $A_{\infty}$-right-$B_{\hbar}$-modules.
\end{itemize}
\end{Def}
$\Modd^{\infty}_{tf}(A_{\hbar})$ and $\Modd^{\infty}_{tf}(B_{\hbar})$
are additive categories. The direct sum $M_\hbar\tilde{\oplus}N_\hbar$
of objects $(M_\hbar\mathrm d_{M_\hbar})$ and $(N_\hbar,\mathrm
d_{N_\hbar})$ in $\Modd^{\infty}_{tf}(A_{\hbar})$ (or
$\Modd^{\infty}_{tf}(B_{\hbar}) )$ is the topologically free module
\[
M_\hbar\tilde{\oplus}N_\hbar:= (M\oplus N)\c1
\]
if $M_\hbar=M\c1$ and $N_\hbar=N\c1$ in $\GD$, endowed with the strictly
unital topological $A_\infty$-module structure given
by the codifferential $\mathrm d_{M_\hbar}\tilde{\oplus}\mathrm
d_{N_\hbar}$.
The natural inclusion and projection
\[
\mathrm i_\hbar :M_\hbar \rightarrow M_\hbar\tilde{\oplus}N_\hbar,
~~~\mathrm p_\hbar: M_\hbar\tilde{\oplus}N_\hbar\rightarrow N_\hbar,
\]
are the strict topological $A_\hbar$-module morphisms defined {\it via}
\[
\mathrm i_\hbar=\mathrm i^{(0)}_\hbar,~~~\mathrm i^{(0)}_\hbar=\mathrm
i^{(0),0}_\hbar=\mathrm i: M\rightarrow M\oplus N,~~m\mapsto m\oplus 0,
\]
and
\[
\mathrm p_\hbar=\mathrm p^{(0)}_\hbar,~~~\mathrm p^{(0)}_\hbar=\mathrm
p^{(0),0}_\hbar=\mathrm p: M\oplus N\rightarrow N,~~m\oplus n\mapsto n.
\]

\subsubsection{On quasi-isomorphisms in $\Modd^{\infty}_{tf}(A_{\hbar})$
and $\Modd^{\infty}_{tf}(B_{\hbar})$ }
The category $\GD$ of {\it all} bigraded $\mathbb K\c1$-modules is
abelian; clearly
$\Modd^{\infty}_{tf}(A_{\hbar})$ and $\Modd^{\infty}_{tf}(B_{\hbar})$
are (not full) subcategories of $\GD$.

In general, the cohomology of a topologically free differential bigraded
$\mathbb K\c1$-module is not topologically free; so we introduce the following definition.

\begin{Def}
A quasi-isomorphism $f_\hbar: N_\hbar\rightarrow M_\hbar$
of objects in $\Modd^{\infty}_{tf}(A_{\hbar})$ (or
$\Modd^{\infty}_{tf}(B_{\hbar})$)
is a morphism of topological $A_\infty$-modules s. t.
$H(f_\hbar): H(N_\hbar)\rightarrow H(M_\hbar)$
is an isomorphism in the abelian category $\GD$.
\end{Def}
Quasi-isomorphisms of strictly unital top. $A_\infty$-modules are not, in general, homotopy equivalences.
A counterexample is given by
\begin{Example}[B.~Keller, \cite{Keller-mail}]
Let $A_\hbar=\mathbb K\c1$ and  $(M_\hbar,\mathrm d_{M_\hbar}) $ be the object in $\Modd^{\infty}_{tf}(A_{\hbar})$ given by
\[
M_\hbar=\{M^0_0\c1, M^1_0\c1, M^2_0\c1 \}, ~~~~M^0_0=M^1_0=\mathbb K,~~ M^2_0=0,
\]
with  codifferential $\mathrm d_{M_\hbar}$ s.t. $\mathrm d^{(0),0}_{M_\hbar}:M^0_0\c1\rightarrow M^1_0\c1$ is the multiplication by $\hbar$, and
$\mathrm d^{(0),1}_{M_\hbar}:M_\hbar\tilde{\otimes}A_\hbar\rightarrow M_\hbar$ is the multiplication in $\mathbb K\c1$. All the other
components are set to be zero.
The strict quasi-isomorphism ($\mathbb K$ is concentrated in bidegree $(0,0)$) of  strictly unital top. $A_\infty$-$A_\hbar$-modules
\[
f_\hbar: M_\hbar\rightarrow \mathbb K
\]
admits no inverse $g_\hbar$ (up to homotopy); in fact  such an inverse would have a $\mathbb K\c1$-linear  $(0,0)$-th component 
$g^{0,0}_\hbar: \mathbb K\rightarrow \mathbb K\c1$. But this implies that $g_\hbar$ is the zero map. 
\end{Example}

\subsubsection{On $\mathcal H^{tf}_\infty (A_\hbar)$, $\mathcal
H^{tf}_\infty (B_\hbar)$ and their triangulated structures. }\label{trtrtr}

As $\Modd^{\infty}_{tf}(A_{\hbar})$ and $\Modd^{\infty}_{tf}(B_{\hbar})$
are additive categories, we can introduce the homotopy categories
\[
\mathcal H^{tf}_\infty (A_\hbar):= \mathcal
H(\Modd^{tf}_{\infty}(A_{\hbar})), ~~~~~~~ \mathcal H^{tf}_\infty
(B_\hbar):=\mathcal
H(\Modd^{tf}_{\infty}(B_{\hbar})).
\]

The objects in $\mathcal H^{tf}_\infty (A_\hbar)$, resp. $\mathcal
H^{tf}_\infty (B_\hbar)$ are the same objects of
$\Modd^{\infty}_{tf}(A_{\hbar})$, resp.
$\Modd^{\infty}_{tf}(B_{\hbar})$. The morphisms are equivalence classes
w.r.t. the equivalence relation
defined as follows; two morphisms
$f_\hbar, g_\hbar :X_\hbar\rightarrow Y_\hbar$ in
$\Modd^{\infty}_{tf}(A_{\hbar})$, resp. $\Modd^{\infty}_{tf}(B_{\hbar})$
are equivalent, i.e. $f_\hbar \sim g_\hbar$, if there exists a strictly unital topological $A_\infty$-homotopy
$H_\hbar$ (see 8.3.1., subsubsection
``On morphisms, quasi-isomorphisms and homotopy equivalences '') s.t.
$f_\hbar-g_\hbar= \mathrm d_{Y_\hbar}\circ H_\hbar+H_\hbar\circ \mathrm
d_{X_\hbar}$,
denoting by $\mathrm d_{X_\hbar}$, resp. $\mathrm d_{Y_\hbar}$ the
codifferentials on $X_\hbar$, resp. $Y_\hbar$.
$\sim$ is an equivalence relation on morphisms in
$\Modd^{\infty}_{tf}(A_{\hbar})$, resp. $\Modd^{\infty}_{tf}(B_{\hbar})$.
We want to prove that $\mathcal H^{tf}_\infty (A_\hbar)$ and $\mathcal
H^{tf}_\infty (B_\hbar)$ are triangulated categories.

\subsubsection{Triangulated structure on $ \mathcal H^{tf}_\infty (A_\hbar)$,
$ \mathcal H^{tf}_\infty (B_\hbar)$ }

We endow the categories $ \mathcal H^{tf}_\infty (A_\hbar)$ and $
\mathcal H^{tf}_\infty (B_\hbar)$ with a triangulated structure such
that, for
$\hbar=0$ it reduces to the triangulated structure on $\mathcal H_\infty
(A)$ and $\mathcal H_\infty (B)$. We refer to Appendix A for
the notation on triangulated categories. As usual we give the definition
for $ \mathcal H^{tf}_\infty (A_\hbar)$; it applies to
$ \mathcal H^{tf}_\infty (B_\hbar)$ as well, with due changes.

Let
$0\rightarrow M_{\hbar}\stackrel{f_{\hbar}}{\rightarrow}
M'_{\hbar}\stackrel{g_{\hbar}}{\rightarrow} M''_{\hbar}\rightarrow 0$
be a short exact sequence of objects in $ \mathcal H^{tf}_\infty (A_\hbar)$
with $f_\hbar$ and $g_\hbar$ strict.
This means that, for any $(i,j)\in\mathbb Z^2$, then
\[
0\rightarrow (M_{\hbar})^i_j\stackrel{f_{\hbar}}{\rightarrow}
(M'_{\hbar})^i_j\stackrel{g_{\hbar}}{\rightarrow}
(M''_{\hbar})^i_j\rightarrow 0
\]
is short exact as sequence of $\mathbb K\c1$-modules.

\begin{Def}
The triangulated structure on the additive category
$ \mathcal H^{tf}_\infty (A_\hbar)$ is given as follows.
The endofunctor $\Sigma$ is simply the (cohomological) grading shift
functor $\Sigma=[1]$.
The distinguished triangles are isomorphic to those induced by
semi-split sequences of strict morphisms

\begin{eqnarray*}
M_{\hbar}\stackrel{f_{\hbar}}{\rightarrow}
M'_{\hbar}\stackrel{g_{\hbar}}{\rightarrow} M''_{\hbar}
\end{eqnarray*}
in $\Modd^{\infty}_{tf}(A_{\hbar})$, i.e. sequences such that
$0\rightarrow M_{\hbar}\stackrel{f_{\hbar}}{\rightarrow}
M'_{\hbar}\stackrel{g_{\hbar}}{\rightarrow} M''_{\hbar}\rightarrow 0$
is an exact sequence in $\GD$, and such that there exists a strict splitting
\[
\rho_\hbar=\sum_{k\geq
0}\rho^{(k)}\hbar^k,~~~\rho^{(k)}=\rho^{(k),0}\in\Hom2^{0,0}_{\C}(M'[1],M[1])
\]
of $f_{\hbar}$, i.e.
\begin{eqnarray}
\rho_\hbar\circ f_{\hbar}=1_\hbar,\label{splitting}
\end{eqnarray}
with
\begin{eqnarray}
\rho_\hbar\circ \mathrm d_{M_\hbar}=\mathrm d_{M'_{\hbar}} \circ
(\rho_\hbar\tilde{\otimes} 1_\hbar^{\tilde{\otimes} i-1}),~~~~~~~~i\geq 2.
\label{Rainbow}
\end{eqnarray}
\end{Def}

By the very definition if the triangulated structure on $\mathcal H^{tf}_\infty (A_\hbar)$ we have
\begin{Cor}
The ``evaluation at $\hbar=0$'' functor $(E_\hbar,1)$, $E_\hbar: \mathcal H^{tf}_\infty (A_\hbar)\rightarrow \Modd_{\infty}(A)/\sim$ with
$E_\hbar(M_\hbar):=M_\hbar /\hbar M_\hbar$, is exact w.r.t. the triangulated category structures on $\mathcal H^{tf}_\infty (A_\hbar)$
and $\Modd_{\infty}(A)/\sim$.
\end{Cor}

\subsubsection{Characterization of exact triangles in $\mathcal
H^{tf}_\infty(A_\hbar)$}\label{trianghbar}

Before proving that the endofunctor $[1]$ and the class of exact
triangles in the above definition endow $\mathcal H^{tf}_\infty(A_\hbar)$
with a triangulated category structure, let us better characterize the
exact triangles. What follows is a suitable topological
$A_\infty$-version of the analysis contained in \cite{GM} on the
triangulated structure of the homotopy category $\mathcal K(\mathcal A)$ of
any additive category $\mathcal A$. The goal is to show that it is
possible to lift those computations to the aforementioned
topological $A_\infty$-case.

\subsubsection{Cones and cylinders. Exact sequences of topologically
free modules}

We recall that, given a topological $A_\infty$-module $M_\hbar$, the
bigraded object
$M_\hbar[\pm 1]$ can be endowed with a topological
$A_{\infty}$-module structure $via$
\begin{eqnarray*}
\mathrm{\bar{d}}_{M_\hbar[\pm 1]}^{(i),l}=-s\circ
\mathrm{\bar{d}}_{M_\hbar}^{(i),l}\circ (
s^{-1}\otimes 1 ).
\end{eqnarray*}

Let $f_{\hbar}:(M_{\hbar},\mathrm d_{M_{\hbar}}) \rightarrow
(N_{\hbar},\mathrm d_{M_{\hbar}})$
be a morphism in $\Modd^{tf}_\infty(A_\hbar)$;
$f_{\hbar}$, $\mathrm d_{M_{\hbar}}$ and $\mathrm d_{N_{\hbar}}$ are
uniquely determined by formal power
series whose $i$-th components are
$f^{(i)}_{\hbar}$, $\mathrm d^{(i)}_{M_{\hbar}}$
and $\mathrm d^{(i)}_{N_{\hbar}}$.
\begin{Def}
A cone $C(f_{\hbar})$ of $f_\hbar$ is the object
\[
C(f_{\hbar}):=M_\hbar[1]\tilde{\oplus}N_\hbar
\]
with topological $A_\infty$-structure given by the differential $\mathrm d_{C(f_\hbar)}$, such that 
\begin{eqnarray}
D_{C(f_\hbar)}=\left(\begin{array}{cc}
\mathrm d_{M_\hbar[1]} & 0 \\
s^{-1}\circ f_{\hbar} & s^{-1}\circ \mathrm d_{N_\hbar}\circ s \\
\end{array}\right)\label{tictac}
\end{eqnarray}

\end{Def}

\begin{Def}
The $A_\infty$-cylinder $\Cyl(f_\hbar)$ is the object
\[
\Cyl(f_\hbar)=M_\hbar\tilde{\oplus} M_\hbar[1]\tilde{\oplus} N_\hbar,
\]
with codifferential $\mathrm d_{\Cyl(f_\hbar)}$ given by
\begin{eqnarray}
D_{\Cyl(f_\hbar)}=\left(\begin{array}{ccc}
s^{-1}\circ \mathrm d_{M_\hbar}\circ s & -\mathrm i_\hbar\circ s^{-1} & 0 \\
0 & \mathrm d_{M_\hbar[1]} & 0 \\
0 & s^{-1}\circ f_{\hbar} &s^{-1}\circ \mathrm d_{N_\hbar}\circ s \\
\end{array}\right).\label{tictac2}
\end{eqnarray}
The natural inclusion
\[
\mathrm i_\hbar: M_\hbar\rightarrow \Cyl(f_\hbar), ~~\mathrm
i_\hbar=\sum_{i\geq 0}\mathrm i^{(i)}_\hbar\hbar^i,
\]
with
\[
\mathrm i^{(i),n}_\hbar=0, ~~for~~i,n\geq 1,
\]
and $\mathrm i^{(0),0}_\hbar=\mathrm i: M\rightarrow M\oplus M[1]\oplus
N$, $m\mapsto (m,0,0)$,
is a strict morphism of topological $A_\infty$-$A_\hbar$-modules.

\end{Def}
\begin{Prop}
For any morphism $f_\hbar$ in $\Modd^{tf}_\infty(A_\hbar)$
\[
\mathrm d_{C(f_\hbar)}\circ \mathrm d_{C(f_\hbar)}=\mathrm
d_{\Cyl(f_\hbar)}\circ \mathrm d_{\Cyl(f_\hbar)}=0.
\]
\end{Prop}
\begin{proof}
Using \eqref{tictac} and \eqref{tictac2}, the proof is immediate.
\end{proof}

For any morphism $f_{\hbar}:(M_{\hbar},\mathrm d_{N_{\hbar}})
\rightarrow (N_{\hbar},\mathrm d_{M_{\hbar}})$
we consider the sequence
\begin{eqnarray}
0\rightarrow M_{\hbar}\stackrel{\mathrm i_{\hbar}}{\rightarrow}
\Cyl(f_{\hbar})\stackrel{\pi_{\hbar}}{\rightarrow}
C(f_{\hbar})\rightarrow 0 \label{canseq}
\end{eqnarray}
The natural projection $\pi_{\hbar}$ is the strict morphism of
topological $A_\infty$-modules
$ \pi_\hbar=\sum_{i\geq 0} \pi^{(i)}_\hbar\hbar^i$, with
$\pi^{(i)}_\hbar=0$, for $i\geq 1$ and
$\pi^{(0),0}_\hbar=\pi: M\oplus M[1]\oplus N\rightarrow M[1]\oplus N$.
The sequence \eqref{canseq} is exact in $\GD$; actually more can be
said: as $\ker2 \pi_\hbar=M_\hbar=\im2 \mathrm i_\hbar$
then \eqref{canseq} is exact in $\Modd^{tf}_{\infty}(A_{\hbar})$.

\begin{Prop}
Let $(M_{\hbar},\mathrm d_{M_{\hbar}})$, $(N_{\hbar}\mathrm
d_{N_{\hbar}})$ and $(L_{\hbar},\mathrm d_{L_{\hbar}})$
be objects in $\Modd^{tf}_\infty(A_\hbar)$ and $f_\hbar:
M_\hbar\rightarrow N_\hbar$, $g_\hbar: N_\hbar\rightarrow L_\hbar$
be strict morphisms in $\Modd^{tf}_\infty(A_\hbar)$.
Any short exact sequence
\[
0\rightarrow M_{\hbar}\stackrel{f_{\hbar}}{\rightarrow}
N_{\hbar}\stackrel{g_\hbar}{\rightarrow} L_{\hbar}\rightarrow 0
\]
in $\GD$ is quasi-isomorphic in $\Modd^{tf}_\infty(A_\hbar)$ to the
short exact sequence $0\rightarrow M_{\hbar}\stackrel{\mathrm i_{\hbar}}{\rightarrow}
\Cyl(f_{\hbar})\stackrel{\pi_{\hbar}}{\rightarrow}
C(f_{\hbar})\rightarrow 0$.

\end{Prop}

\begin{proof}
Like in \cite{GM}, prop. 5, section III, with due changes.
\end{proof}

For any morphism $f_{\hbar}:M_{\hbar} \rightarrow N_{\hbar}$ in
$\Modd^{tf}_\infty(A_\hbar)$ the sequence

\begin{eqnarray*}
0\rightarrow M_{\hbar}\stackrel{\mathrm i_{\hbar}}{\rightarrow}
\Cyl(f_{\hbar})\stackrel{\pi_{\hbar}}{\rightarrow}
C(f_{\hbar})\rightarrow 0
\end{eqnarray*}
is exact. Something more can be said; in
fact the sequence is semi-split with strict splitting
$\rho_\hbar: \Cyl(f_\hbar) \rightarrow M_\hbar$ given by
\[
\rho^{(0)}_\hbar(m,sm',l)=m,~~\rho^{(i)}_\hbar=0~~\mbox{for}~~i\geq 1.
\]
It is important to note that $\rho_\hbar$ does not commute with the
components $\mathrm d^{(i),0}_{\Cyl(f_\hbar)}$ and
$\mathrm d^{(i),0}_{M_\hbar}$ of the codifferentials on $\Cyl(f_\hbar)$
and $M_\hbar$ for any $i\geq 0$, but
\[
\rho^{(0)}_\hbar(s^{-1}(\mathrm
d^{(j),n}_{\Cyl(f_\hbar)}(m,sm',l|a^{\otimes n} ) ))=
\mathrm d^{(j),n}_{M_\hbar}(\rho^{(0)}_\hbar(m,sm',l )|a^{\otimes
n}),~~\mbox{for any}~~n\geq 1.
\]
Like in the classical case, the presence of the inclusion $\mathrm
i_\hbar$ in the definition of the codifferential
$\mathrm d_{\Cyl(f_\hbar)}$ plays a major role.
In summary,
\[
M_{\hbar}\stackrel{\mathrm i_{\hbar}}{\rightarrow}
\Cyl(f_{\hbar})\stackrel{\pi_{\hbar}}{\rightarrow}
C(f_{\hbar})\rightarrow M_\hbar[1]
\]
is an exact triangle in $\mathcal H^{tf}_\infty(A_\hbar)$
for any morphism $f_{\hbar}:M_{\hbar} \rightarrow L_{\hbar}$ in
$\Modd^{tf}_\infty(A_\hbar)$.

Let \begin{eqnarray}
M_\hbar\stackrel{f_\hbar}{\rightarrow} L_\hbar
\stackrel{p'_\hbar}{\rightarrow}
C(f_\hbar)\stackrel{r_\hbar}{\rightarrow} M_\hbar[1];\label{q}
\end{eqnarray}
be a sequence in $\mathcal H^{tf}_\infty(A_\hbar)$; it is isomorphic in
$\mathcal H^{tf}_\infty(A_\hbar)$ to the exact triangle
$M_{\hbar}\stackrel{\mathrm i_{\hbar}}{\rightarrow}
\Cyl(f_{\hbar})\stackrel{\pi_{\hbar}}{\rightarrow}
C(f_{\hbar})\stackrel{r_\hbar}{\rightarrow} M_\hbar[1]$ $via$
\begin{eqnarray}
\begin{CD}
M_\hbar @> f_\hbar >> L_\hbar @> p'_\hbar >> C(f_\hbar) @> r_\hbar >>
M_\hbar[1] \\
@VV 1_\hbar V @VV \alpha_\hbar V @VV 1_\hbar V @VV 1_\hbar V \\
M_\hbar @>\mathrm i_{\hbar} >> \Cyl(f_\hbar) @> \pi_\hbar >> C(f_\hbar)
@> r_\hbar >> M_\hbar[1]
\end{CD}\label{d}
\end{eqnarray}
with strict $A_\infty$-morphism
\[
\alpha_\hbar: L_\hbar\rightarrow \Cyl(f_\hbar),~~~
\alpha^{(0),0}_\hbar(l)=(0,0,l)
\]
and $\alpha^{(i)}_\hbar=0$, for $i\geq 1$.
In summary, \eqref{q} is an exact triangle in $\mathcal
H^{tf}_\infty(A_\hbar)$, as well.

\subsubsection{Other exact triangles: using the splitting}

Let
\begin{eqnarray*}
0\rightarrow M_{\hbar}\stackrel{f_{\hbar}}{\rightarrow}
N_{\hbar}\stackrel{g'_{\hbar}}{\rightarrow} Q_{\hbar}\rightarrow 0
\end{eqnarray*} be a semi-split exact sequence like in def. 25, with
splitting $\rho_\hbar$ and
$N_\hbar=N\c1$, $M_\hbar=M\c1$, $Q_\hbar=Q\c1$ in $\GD$.

At the order $\hbar^0$
eq. \eqref{splitting} is equivalent to $\rho^{(0)}_\hbar\circ
f^{(0)}_\hbar=1$. This implies that $N^i_j\cong M^i_j\oplus Q^i_j$
as $\mathbb K$-modules, for any $(i,j)\in\mathbb Z^2$; in virtue of this
we assume that $N_\hbar=((M\oplus Q)\c1,\mathrm d_{M_\hbar\tilde{\oplus} Q_\hbar})$,
where
\[
\mathrm d_{M_\hbar\tilde{\oplus} Q_\hbar}=\left(\begin{array}{cc}
\mathrm d_{M_\hbar} & -f_\hbar \\
0 & \mathrm d_{Q\hbar} \\
\end{array}\right).
\]
and $M_\hbar\tilde{\oplus} Q_\hbar\equiv (M\oplus Q)\c1$ in $\GD$. It follows that $\mathrm d_{M_\hbar\tilde{\oplus} Q_\hbar}\circ \mathrm
d_{M_\hbar\tilde{\oplus} Q_\hbar}=0$ if and only if
\[
f_\hbar:Q_\hbar\rightarrow M_\hbar[1]
\]
defines an $A_\infty$-morphism of topological $A_\infty$-$A_\hbar$-modules.
By definition of the triangulated structure on $\mathcal
H^{tf}_\infty(A_\hbar)$, the sequence $M_\hbar\stackrel{i_\hbar}{\rightarrow} M_\hbar\tilde{\oplus} Q_\hbar
\stackrel{p_\hbar}{\rightarrow} Q_\hbar\stackrel{f_\hbar}{\rightarrow}M_\hbar[1]$ is an exact triangle.
\begin{Thm}
The homotopy categories $\mathcal H^{tf}_\infty(A_\hbar)$ and $\mathcal
H^{tf}_\infty(B_\hbar)$ are triangulated; the triangulated
structure is the one given in def. 25.
\end{Thm}
\begin{proof}
\cite{GM}, pag. 246, with due changes; we sketch the proof for sake of
clarity. On the axiom $(T1)$ (see the Appendix);
the sequence
\[
X_\hbar\stackrel{1_\hbar}{\rightarrow} X_\hbar\rightarrow 0\rightarrow
X_\hbar[1]
\]
is isomorphic to $X_\hbar\stackrel{1_\hbar}{\rightarrow}
X_\hbar\rightarrow C(1_\hbar)\rightarrow X_\hbar[1]$ as the zero morphism
$0\rightarrow C(1_\hbar)$ is homotopic to the identity morphism
$1'_\hbar: C(1_\hbar)\rightarrow C(1_\hbar)$ on $C(1_\hbar)$;
in fact
\[
1'_\hbar=H_\hbar\circ \mathrm d_{C(1_\hbar)}+ \mathrm
d_{C(1_\hbar)}\circ H_\hbar,
\]
with strict topological $A_\infty$-homotopy $H_\hbar=H^{(0),0}$,
$H^{(0),0}(sx\oplus x')=(x',0)$.
Compatibility with the codifferentials follows easily.
Axiom $(T2)$ is proved similarly. Let
\[
X_\hbar\stackrel{u_\hbar}{\rightarrow}
Y_\hbar\stackrel{v_\hbar}{\rightarrow} C(u_\hbar)
\stackrel{p_\hbar}{\rightarrow} X_\hbar[1]
\]
be an exact triangle. We want to prove that the sequence
\[
Y_\hbar\stackrel{v_\hbar}{\rightarrow} C(u_\hbar)
\stackrel{p_\hbar}{\rightarrow} X_\hbar[1]
\stackrel{-u_\hbar[1]}{\rightarrow} Y_\hbar[1]
\]
is isomorphic to the exact triangle
\[
Y_\hbar\stackrel{v_\hbar}{\rightarrow} C(u_\hbar)
\stackrel{s_\hbar}{\rightarrow} C(v_\hbar)
\stackrel{-u_\hbar[1]}{\rightarrow} Y_\hbar[1].
\]
All we need is to introduce the topological $A_\infty$-homotopy equivalence
\[
\theta_\hbar: X_\hbar[1]\rightarrow C(v_\hbar),
\]
with
\[
\theta^{(0),0}_\hbar(sx)=(-su^{(0),0}_\hbar(sx),sx,0),~~~~~~\theta^{(i),n}_\hbar(x|a_1|\dots|a_n)=(-su^{(i),n}((x|a_1|\dots|a_n)),0,0),
~~n\geq1
\]
and to check that $s_\hbar\circ 1_\hbar-\theta_\hbar\circ p_\hbar= \mathrm
d_{C(v_\hbar)}H_\hbar+H_\hbar\mathrm d_{C(u_\hbar)}$ with strict $A_\infty$-homotopy $H_\hbar: C(u_\hbar)\rightarrow
C(v_\hbar)$, $H^{(0),0}(sx,y)=(y,0,0)$, $H^{(i),n}=0$ otherwise.

The computations are long but straightforward; $\theta_\hbar$ is a
homotopy equivalence because it admits the strict inverse
\[
\psi_\hbar :C(v_\hbar)\rightarrow X_\hbar[1],~~~~
\psi_\hbar^{(0),0}(sy,sx,y')=sx,
\]
and $\psi_\hbar^{(i),n}=0$ otherwise. Clearly
$\psi_\hbar\circ\theta_\hbar=1_\hbar$, but $\theta_\hbar\circ\psi_\hbar=1_\hbar+\mathrm
d_{C(v_\hbar)}{H'}_\hbar+{H'}_\hbar\mathrm d_{C(v_\hbar)}$,
with ${H'}^{(0),0}_\hbar(sy,sx,y')=(y',0,0)$ and zero otherwise.

Axiom $(T3)$ is proved by using cones and $(T4)$ follows by using semi
split exact sequences.

\end{proof}

\subsubsection{Localizing w.r.t.\ topological
$A_\infty$-quasi-isomorphisms: on the derived categories
$\DD^{\infty}_{tf}(A_{\hbar})$ and
$\DD^{\infty}_{tf}(B_{\hbar})$}
 In \cite{GM}, def.6, section III, localizing classes of morphisms are defined. In our setting we have
\begin{Prop}
The class $Qis$ of quasi-isomorphisms in the homotopic categories $
\mathcal H^{tf}_\infty (A_\hbar)$ and $\mathcal H^{tf}_\infty (B_\hbar)$
is localizing.
\end{Prop}
\begin{proof} We prove the statement for $ \mathcal H^{tf}_\infty
(A_\hbar)$.
We refer to the proof of thm. 4, pag 160 in \cite{GM}. We
``translate'' it in our topological $A_\infty$-case, with due changes.
\end{proof}

Thanks to the above proposition the following definition makes sense.

\begin{Def}
The localizations
\[
\DD^{\infty}_{tf}(A_{\hbar}):= \mathcal H^{tf}_\infty
(A_\hbar)[Qis^{-1}],~~resp.
~~\DD^{\infty}_{tf}(B_{\hbar}):= \mathcal H^{tf}_\infty (B_\hbar)[Qis^{-1}]
\]
are said to be the derived categories of
$\Modd^{tf}_\infty(A_\infty)$, resp. $\Modd^{tf}_\infty(B_\infty)$.
\end{Def}

The objects in $\DD^{\infty}_{tf}(A_{\hbar})$, resp.
$\DD^{\infty}_{tf}(B_{\hbar})$ are the same objects of
$ \mathcal H^{tf}_\infty (A_\hbar)$, resp. $ \mathcal H^{tf}_\infty
(B_\hbar)$ while the morphisms are described through the equivalence
classes of ``roofs'', as in \cite{GM}. We use the notation $\mathcal
D=\DD^{\infty}_{tf}(A_{\hbar}), \DD^{\infty}_{tf}(B_{\hbar})$.
Any morphism $\phi_\hbar: X_\hbar\rightarrow Y_\hbar$ in $\mathcal D$
is represented by an equivalence class of roofs; 
if two roofs $(s_\hbar, \bar{\phi}_\hbar)$ and $(t_\hbar, \bar{\psi}_\hbar)$ representing the same
 morphism in $\mathcal D$  are
equivalent, we will simply write
$(s_\hbar, \bar{\phi}_\hbar)=(t_\hbar, \bar{\psi}_\hbar)$. 

In what follows we will state that the morphism $\phi_\hbar:
X_\hbar\rightarrow Y_\hbar$ in $\mathcal D$
is represented by {\it the} roof $(s_\hbar, \bar{\phi_\hbar})$, for
simplicity. The identity morphism $1_\hbar: X_\hbar\rightarrow X_\hbar$ in $\mathcal
D$ is represented by
\begin{diagram}
& & X_\hbar & & \\
& \ldTo^{1_\hbar} & & \rdTo^{1_\hbar} & \\
X_\hbar & & & & X_\hbar \\
\end{diagram}
for any $X_\hbar$ in $\mathcal D$.
The composition
\begin{eqnarray}
\psi_\hbar\circ \phi_\hbar \label{Fantozzo}
\end{eqnarray}
of morphisms $\phi_\hbar: X_\hbar\rightarrow Y_\hbar$, $\psi_\hbar:
Y_\hbar\rightarrow Z_\hbar$ in $\mathcal D$
represented by the roofs $(s_\hbar, \bar{\phi}_\hbar)$ and $(t_\hbar,
\bar{\psi}_\hbar)$  will be
denoted also by
\[
(t_\hbar, \bar{\psi}_\hbar) \circ (s_\hbar, \bar{\phi}_\hbar).
\]

\begin{Cor}
The class of quasi-isomorphisms in $\mathcal H^{tf}_\infty(A_\hbar)$ and
$\mathcal H^{tf}_\infty(B_\hbar)$ is compatible with triangulation;
the derived categories
\[
\DD^{\infty}_{tf}(A_{\hbar}),~~\DD^{\infty}_{tf}(B_{\hbar}).
\]
are triangulated.
\end{Cor}
\begin{proof}
See \cite{GM}; the proofs there apply here with straightforward changes.
\end{proof}

\subsubsection{On the quantized Functors }

Let us define the functors
\begin{eqnarray*}
F_{\hbar}: \Modd^{tf}_{\infty}(A_{\hbar})\rightarrow
\Modd^{tf,strict}_{\infty}(B_{\hbar}), &G_{\hbar}:
\Modd^{tf}_{\infty}(B_{\hbar}) \rightarrow
\Modd^{tf,strict}_{\infty}(A_{\hbar}),
\end{eqnarray*}
{\it via}
\begin{eqnarray*}
F_{\hbar}(M_{\hbar}):=
M_\hbar\underline{\tilde{\otimes}}_{A_\hbar}K_{\hbar}, & G_{\hbar}(N_{\hbar}):=
N_\hbar\underline{\tilde{\otimes}}_{A_\hbar}\underline{K}_{\hbar},
\end{eqnarray*}
on objects $M_{\hbar}\in \Modd^{tf}_{\infty}(A_{\hbar})$ and $N_{\hbar}\in
\Modd^{tf}_{\infty}(B_{\hbar})$.
Let $f_\hbar:M_{\hbar}\rightarrow N_{\hbar}$ be a morphism in
$\Modd^{tf}_{\infty}(A)$. Then $F_{\hbar}(f_\hbar)$ is the strict morphism in
$\Modd^{tf}_{\infty}(B)$ given by

\begin{eqnarray*}
F_{\hbar}(f_\hbar)=\sum_{i\geq 0}F^{(i)}_{\hbar}(f_\hbar)\hbar^i, &
F^{(i),0}_{\hbar}(f_\hbar)=\sum_{k\geq 0}f^{(i),k}_\hbar\otimes 1.
\end{eqnarray*}

and zero otherwise.
Similar definition holds true for $G'_{\hbar}$. Here
$\Modd^{tf,strict}_{\infty}(A_{\hbar})$ denotes the subcategory of
$\Modd^{tf}_{\infty}(A_{\hbar})$ with same objects and strict
topological $A_\infty$-morphisms. Same convention holds true for
$\Modd^{tf,strict}_{\infty}(B_{\hbar})$.

\begin{Lem}
Let $f_{\hbar}:M_\hbar\rightarrow N_\hbar$ be a quasi-isomorphism in
$\Modd^{tf}_{\infty}(A_{\hbar})$; then
$F_{\hbar}(f_\hbar):M_{\hbar}\underline{\tilde{\otimes}}_{A_{\hbar}}K_{\hbar}\rightarrow
N_{\hbar}\underline{\tilde{\otimes}}_{A_{\hbar}}K_{\hbar}$ is a 
quasi-isomorphism in
$\Modd^{tf,strict}_{\infty}(A_{\hbar})$.
\end{Lem}
Similar considerations hold for the functor $G_{\hbar}$.

\subsubsection{The quantized functors on the derived categories;
compatibility with the triangulated structures}

We discuss now the above quantized functors lifting them on the derived
categories $\DD^{\infty}_{tf}(A_{\hbar})$ and
$\DD^{\infty}_{tf}( B_{\hbar})$.

\begin{Def}
$\mathcal F_\hbar$ is the unique functor
$\mathcal F_\hbar: \DD^{\infty}_{tf}(A_{\hbar})\rightarrow
\DD^{\infty}_{tf}( B_{\hbar})$ s.t.
\[
\mathcal F_\hbar\circ Q_{A_\hbar}=\mathcal T_\hbar,
\]
denoting by $Q_{A_\hbar}: \mathcal H^{tf}_\infty(A_\hbar)\rightarrow
\DD^{\infty}_{tf}(A_{\hbar})$ the canonical functor
\[
Q_{A_\hbar}(X)=X,~~~~~~Q_{A_\hbar}(f_\hbar)=(1,f_\hbar),
\]
and by $\mathcal T_\hbar: \mathcal H^{tf}_\infty(A_\hbar)\rightarrow
\DD^{\infty}_{tf}( B_{\hbar})$
the composition
\[
\mathcal T_\hbar= Q_{B_\hbar}\circ \bar{F}_{\hbar},
\]
where $\bar{F}_{\hbar}:\mathcal H^{tf}_\infty(A_\hbar)\rightarrow
\mathcal H^{tf}_\infty(B_\hbar)$ is the functor induced by $F_\hbar$
on the homotopy categories of $\Modd^{tf}_{\infty}(A_{\hbar})$ and
$\Modd^{tf}_{\infty}(B_{\hbar})$.
\end{Def}
The functor $\mathcal G_\hbar: \DD^{\infty}_{tf}(B_{\hbar})\rightarrow
\DD^{\infty}_{tf}(A_{\hbar})$ is defined similarly.
By definition
\[
\mathcal F_\hbar(X_\hbar)=\bar{F}_\hbar(X_\hbar)=F_\hbar(X_\hbar),
\]
on every object $X_\hbar\in \DD^{\infty}_{tf}(A_{\hbar})$ and on any morphism
$(s_\hbar,f_\hbar)$ in $\DD^{\infty}_{tf}(A_{\hbar})$:
\begin{eqnarray}
\mathcal
F_\hbar(s_\hbar,f_\hbar)=(\bar{F}_{\hbar}(s_\hbar),\bar{F}_\hbar(f_\hbar)).
\end{eqnarray}
Both $(\mathcal F_\hbar,{\varphi}^1_\hbar)$ and $(\mathcal G_\hbar,{\varphi}^2_\hbar)$ are exact functors
w.r.t. the triangulated
structure on $\DD^{\infty}_{tf}(A_{\hbar})$ and
$\DD^{\infty}_{tf}(B_{\hbar})$. Here
$\varphi^1_\hbar :\mathcal F_\hbar\circ [1]\rightarrow [1]\circ\mathcal F_\hbar$ denotes the obvious morphism of functors, and similarly for
$\varphi^2_\hbar$.

\begin{Prop}\label{marroni}
Let $(\mathcal F_\hbar,\mathcal G_\hbar)$ be the pair of functors
defined above.
\begin{itemize}
\item $A_\hbar$ is isomorphic to $\mathcal G_\hbar(\mathcal
F_\hbar(A_\hbar))$ in $\DD^{\infty}_{tf}(A_{\hbar})$.
\item $K_\hbar$ is isomorphic to $\mathcal F_\hbar(\mathcal
G_\hbar(K_\hbar))$ in $\DD^{\infty}_{tf}(B_{\hbar})$.
\end{itemize}
\end{Prop}
\begin{proof}
The first isomorphism is represented by 
\begin{diagram}
& & A_\hbar & & \\
& \ldTo^{1_\hbar} & & \rdTo^{\bar{\psi}_{A_\hbar}} & \\
A_\hbar & & & & \mathcal G_\hbar(\mathcal
F_\hbar(A_\hbar)) \\
\end{diagram}
with 
\[
\bar{\psi}_{A_\hbar}:A_\hbar\rightarrow K_\hbar
\underline{\tilde{\otimes}}_{B_\hbar} \underline{K}_\hbar
\rightarrow
(A_\hbar\underline{\tilde{\otimes}}_{A_\hbar} K_\hbar )
\underline{\tilde{\otimes}}_{B_\hbar} \underline{K}_\hbar=
\mathcal G_\hbar(\mathcal
F_\hbar(A_\hbar));
\]
the second isomorphism is represented by 
\begin{diagram}
& & K_\hbar & & \\
& \ldTo^{1_\hbar} & & \rdTo^{\bar{\psi}_{K_\hbar}} & \\
K_\hbar & & & & \mathcal F_\hbar(\mathcal
G_\hbar(K_\hbar)) \\
\end{diagram}
with
\[
\bar{\psi}_{K_\hbar}: K_\hbar  \rightarrow A_\hbar\underline{\tilde{\otimes}}_{A_\hbar} K_\hbar 
\rightarrow (K_\hbar\underline{\tilde{\otimes}}_{B_\hbar}
\underline{K}_\hbar)
\underline{\tilde{\otimes}}_{A_\hbar} K_\hbar=
\mathcal F_\hbar(\mathcal
G_\hbar(K_\hbar))
\]
$\bar{\psi}_{A_\hbar}$ and $\bar{\psi}_{K_\hbar}$ are defined in cor.~\ref{party2}.
\end{proof}

\section{Main result}
Denoting by $\triang^{\infty}_{A_{\hbar}}(A_{\hbar})$ the triangulated
subcategory of $\DD^{\infty}_{tf}(A_{\hbar})$
generated by   $A_\hbar[i]\langle j\rangle$  and by
$\triang^{\infty}_{B_{\hbar}}(K_{\hbar})$ the triangulated subcategory
of $\DD^{\infty}_{tf}(B_{\hbar})$
generated by  $K_\hbar[i]\langle j\rangle$,
$i,j\in\mathbb Z$,
we arrive at the main result of these notes.
\begin{Thm}\label{Thm30}
Let $X$ be a finite dimensional vector space over $\mathbb K=\mathbb R$,
or $\mathbb C$ and $(A,K,B)$ be the triple of bigraded
$A_{\infty}$-structures
introduced in section 6. By $\hbar\pi\in
(T_{poly}(X)\c1,0,[\cdot,\cdot]_{\hbar})$ we denote an
$\hbar$-formal quadratic Poisson bivector on $X$ such that
$(A_{\hbar},K_{\hbar},B_{\hbar})$ is the Deformation Quantization on
$(A,K,B)$ w.r.t.\ $\hbar\pi$.
The functor
\begin{eqnarray*}
\mathcal F_\hbar : \DD^{\infty}_{tf}(A_{\hbar})\rightarrow
\DD^{\infty}_{tf}( B_{\hbar}), ~~~~~\mathcal
F_\hbar(\bullet)=\bullet~\underline{\tilde{\otimes}}_{A_\hbar} K_\hbar
\end{eqnarray*}
induces equivalences of triangulated categories
\begin{eqnarray*}
\triang^{\infty}_{A_{\hbar}}(A_{\hbar})\simeq
\triang^{\infty}_{B_{\hbar}}(K_{\hbar}), &
\thick^{\infty}_{A_{\hbar}}(A_{\hbar})\simeq\thick^{\infty}_{B_{\hbar}}(K_{\hbar}).
\end{eqnarray*}
Let $(\tilde{K},\mathrm d_{\tilde{K}})$ be the
$A_{\infty}$-$B$-$A$-bimodule
with $\tilde{K}=K$ and $\mathrm d_{\tilde{K}}$ obtained from $\mathrm
d_K$ exchanging $A$ and $B$ and
$(\tilde{K}_{\hbar},\mathrm d_{\tilde{K}_{\hbar}})$ be its quantization
w.r.t.\ $\pi_{\hbar}$; the functor
\begin{eqnarray*}
\mathcal F^{''}_\hbar : \DD^{\infty}_{tf}(B_{\hbar})\rightarrow
\DD^{\infty}_{tf}( A_{\hbar}), ~~~~~\mathcal
F^{''}_\hbar(\bullet)=\bullet~\underline{\tilde{\otimes}}_{B_\hbar}
\tilde{K}_\hbar
\end{eqnarray*}
induces the equivalence of triangulated categories
\begin{eqnarray*}
\triang^{\infty}_{A_{\hbar}}(\tilde{K}_{\hbar})\simeq
\triang^{\infty}_{B_{\hbar}}(B_{\hbar}),&
\thick^{\infty}_{A_{\hbar}}(\tilde{K}_{\hbar})\simeq\thick^{\infty}_{B_{\hbar}}(B_{\hbar}).
\end{eqnarray*}
\end{Thm}

\appendix

\section{Proof of prop.~\ref{End-bim} }

\begin{itemize}
 \item  {\it On the quadratic relations $\mathrm d^2_{\underline{\End2}_{B}(K)}=0$.}
\end{itemize}

First of all we note that the maps $\mathrm
D_{\underline{\End2}_{B}(K)}^{n,0}(a_1|\dots|a_n|\varphi)$ and
$\mathrm D_{\underline{\End2}_{B}(K)}^{0,m}(\varphi|a_1|\dots|a_m)$ have
cohomological degree $2$; we have already remarked that
$\mathcal L_A(a_1|\dots|a_n)$ has cohomological degree $1$, instead. The relations 
$\mathrm{\bar{d}}^{0,0}_{\underline{\End2}_{B}(K)}(\mathrm{\bar{d}}^{0,0}_{\underline{\End2}_{B}(K)}(s\varphi))=0$ are immediate to prove.
 We  prove the case $n\geq 2$, $m=0$, i.e.
\begin{eqnarray*}
&&\sum_{j=1}^n(-1)^{\sum_{i=1}^{j-1}(|a_i|-1)}\mathrm{\bar{d}}^{n-1,0}_{\underline{\End2}_{B}(K)}
(a_1|\dots|a_{j-1},\mathrm{\bar{d}}_A^2(a_j|a_{j+1})|a_{j+2}|\dots|a_n|\varphi)+\\
&&\sum_{n'=1}^{n-1}(-1)^{\sum_{i=1}^{n-n'}(|a_i|-1)}\mathrm{\bar{d}}^{n-n',0}_{\underline{\End2}_{B}(K)}
(a_1|\dots|a_{n-n'},\mathrm{\bar{d}}^{n',0}_{\underline{\End2}_{B}(K)}(a_{n-n'+1}|\dots|a_n|\varphi)+\\
&&\mathrm{\bar{d}}^{0,0}_{\underline{\End2}_{B}(K)}(\mathrm{\bar{d}}^{n,0}_{\underline{\End2}_{B}(K)}(a_1|\dots|a_n|\varphi))+
(-1)^{\sum_{i=1}^{n}(|a_i|-1)}
\mathrm{\bar{d}}^{n,0}_{\underline{\End2}_{B}(K)}(a_1|\dots|a_n,\mathrm{\bar{d}}^{0,0}_{\underline{\End2}_{B}(K)}(\varphi))=0;
\end{eqnarray*}
such quadratic relations are equivalent to
\begin{eqnarray}
&&\sum_{j=1}^n(-1)^{\sum_{i=1}^{j-1}(|a_i|-1)+\sum_{i=1}^n(|a_i|-1)}
\mathcal
L_A(a_1|\dots|a_{j-1},\mathrm{\bar{d}}_A^2(a_j|a_{j+1})|a_{j+2}|\dots|a_n)\circ\varphi+\nonumber\\
&&\sum_{n'=1}^{n-1}(-1)^{\sum_{i=n-n'+1}^{n}(|a_i|-1)}\mathcal
L_A(a_1|\dots|a_{n-n'})\circ\left(
\mathcal L_A(a_{n-n'+1}|\dots|a_n|)\circ\varphi\right)+\nonumber\\
&&(-1)^{\sum_{i=1}^{n}(|a_i|-1)}\partial_{\underline{\End2}_{B}(K)}
(\mathcal L_A(a_1|\dots|a_{n})\circ \varphi)+
\mathcal L_A(a_1|\dots|a_{n})\circ
\partial_{\underline{\End2}_{B}(K)}(\varphi)
=0\nonumber\\
\label{questa}
\end{eqnarray}
The last contributions on the l.h.s. of \eqref{questa} can be written as
\begin{eqnarray*}
&&(-1)^{\sum_{i=1}^{n}(|a_i|-1)}\partial_{\underline{\End2}_{B}(K)}
(\mathcal L_A(a_1|\dots|a_{n})\circ \varphi)+
\mathcal L_A(a_1|\dots|a_{n})\circ
\partial_{\underline{\End2}_{B}(K)}(\varphi)=\\
&&(-1)^{\sum_{i=1}^{n}(|a_i|-1)}\partial_{\underline{\End2}_{B}(K)}
(\mathcal L_A(a_1|\dots|a_{n}))\circ \varphi=\\
&&(\mathcal L_A(a_1|\dots|a_{n})\circ\mathrm
d_{K}+(-1)^{\sum_{i=1}^{n}(|a_i|-1)}\mathrm d_K\circ \mathcal
L_A(a_1|\dots|a_{n}) )\circ \varphi,
\end{eqnarray*}
because $|\mathcal L_A(a_1|\dots|a_{n})|=\sum_{i=1}^{n}(|a_i|-1)+1$.

At the end, multiplying both sides of \eqref{questa}
by $(-1)^{\sum_{i=1}^{n}(|a_i|-1)}$ and using the associativity of the
product $\circ$
we get that \eqref{questa}
is equivalent to a finite sum of equations of the type $\mathrm
d_K^2(a_1|\dots|a_n|\varphi(\dots))=0$.

We continue with the case $n=0$, $m\geq 2$, i.e.

\begin{eqnarray*}
&&\sum_{j=1}^m(-1)^{|\varphi|-1+\sum_{i=1}^{j-1}(|a_i|-1)}\mathrm{\bar{d}}^{0,m-1}_{\underline{\End2}_{B}(K)}
(\varphi|a_1|\dots|a_{j-1},\mathrm{\bar{d}}_A^2(a_j|a_{j+1})|a_{j+2}|\dots|a_m)+\\
&&\sum_{m'=1}^{m-1}\mathrm{\bar{d}}^{0,m-m'}_{\underline{\End2}_{B}(K)}(
\mathrm{\bar{d}}^{0,m'}_{\underline{\End2}_{B}(K)}(\varphi|a_1|\dots|a_{m'}),a_{m'+1}|\dots|a_m))+\\
&&\mathrm{\bar{d}}^{0,0}_{\underline{\End2}_{B}(K)}(\mathrm{\bar{d}}^{0,m}_{\underline{\End2}_{B}(K)}(\varphi|a_1|\dots|a_m))+
\mathrm{\bar{d}}^{0,m}_{\underline{\End2}_{B}(K)}(\mathrm{\bar{d}}^{0,0}_{\underline{\End2}_{B}(K)}(\varphi),a_1|\dots|a_m)=0.
\end{eqnarray*}
The above relations are equivalent to
\begin{eqnarray*}
&&\sum_{j=1}^m(-1)^{\sum_{i=1}^{j-1}(|a_i|-1)}
\varphi\circ \mathcal
L_A(a_1|\dots|a_{j-1},\mathrm{\bar{d}}_A^2(a_j|a_{j+1})|a_{j+2}|\dots|a_m)+\nonumber\\
&&\sum_{m'=1}^{m-1}(-1)^{\sum_{i=1}^{m-m'}(|a_i|-1)}\varphi\circ\left(\mathcal
L_A(a_1|\dots|a_{m-m'})\circ
\mathcal L_A(a_{m-m'+1}|\dots|a_m|)\right)+\\
&&(-1)^{|\varphi|+1}\partial_{\underline{\End2}_{B}(K)} (\varphi\circ
\mathcal L_A(a_1|\dots|a_{m}))+
(-1)^{|\varphi|}\partial_{\underline{\End2}_{B}(K)}
(\varphi)\circ\mathcal L_A(a_1|\dots|a_{m})=0,
\end{eqnarray*}
which are easily verified, as
\begin{eqnarray*}
&& (-1)^{|\varphi|+1}\partial_{\underline{\End2}_{B}(K)} (\varphi\circ
\mathcal L_A(a_1|\dots|a_{m}))+
(-1)^{|\varphi|}\partial_{\underline{\End2}_{B}(K)}
(\varphi)\circ\mathcal L_A(a_1|\dots|a_{m})=\\
&&-\varphi\circ \partial_{\underline{\End2}_{B}(K)} (\mathcal
L_A(a_1|\dots|a_{m})=\\
&&\varphi\circ \left( (-1)^{\sum_{i=1}^{n}(|a_i|-1)+1}\mathcal
L_A(a_1|\dots|a_{m})\circ \mathrm d_K-
\mathrm d_K\circ\mathcal L_A(a_1|\dots|a_{m}) \right).
\end{eqnarray*}
Note the overall $-1$ sign, which plays no role.

The equations expressing compatibility between the left and right
actions on $\underline{\End2}_{B}(K)$ (for $n,m\geq 1$ ) are
\begin{eqnarray*}
&&(-1)^{\sum_{i=1}^n(|a_i|-1)}\mathrm{\bar{d}}^{n,0}_{\underline{\End2}_{B}(K)}
(a_1|\dots|a_n,\mathrm{\bar{d}}^{0,m}_{\underline{\End2}_{B}(K)}(\varphi|\bar{a}_1|\dots|\bar{a}_m))+\\
&&
\mathrm{\bar{d}}^{0,m}_{\underline{\End2}_{B}(K)}(\mathrm{\bar{d}}^{n,0}_{\underline{\End2}_{B}(K)}
(a_1|\dots|a_n,|\varphi)|\bar{a}_1|\dots|\bar{a}_m)=0;
\end{eqnarray*}
they are equivalent to
\begin{eqnarray*}
&& \mathcal L_A(a_1|\dots|a_n)\circ \mathrm
D^{0,m}_{\underline{\End2}_{B}(K)}(\varphi|\bar{a}_1|\dots|\bar{a}_m)+\\
&&(-1)^{|\varphi|+\sum_{i=1}^n(|a_i|-1)}\mathrm
D^{n,0}_{\underline{\End2}_{B}(K)}(a_1|\dots|a_n|\varphi)\circ
\mathcal L_A(\bar{a}_1|\dots|\bar{a}_m)=0,
\end{eqnarray*}
or
\begin{eqnarray*}
&& (-1)^{|\varphi|+1}\mathcal L_A(a_1|\dots|a_n)\circ\left(\varphi\circ
\mathcal L_A(\bar{a}_1|\dots|\bar{a}_m) \right) +\\
&&(-1)^{|\varphi|}\left(\mathcal
L_A(a_1|\dots|a_n)\circ\varphi\right)\circ \mathcal
L_A(\bar{a}_1|\dots|\bar{a}_m)=0.
\end{eqnarray*}

We finish by checking the compatibility of the actions with the
differential, i.e.
\begin{eqnarray}
(-1)^{|a|-1} \mathrm{\bar{d}}^{1,0}_{\underline{\End2}_{B}(K)}
(sa,\mathrm{\bar{d}}^{0,0}_{\underline{\End2}_{B}(K)}(\varphi))+\mathrm{\bar{d}}^{0,0}_{\underline{\End2}_{B}(K)}
(\mathrm{\bar{d}}^{1,0}_{\underline{\End2}_{B}(K)}(a|\varphi))=0,\label{compsx}
\end{eqnarray}
and
\begin{eqnarray}
\mathrm{\bar{d}}^{0,1}_{\underline{\End2}_{B}(K)}
(\mathrm{\bar{d}}^{0,0}_{\underline{\End2}_{B}(K)}(\varphi)|a
)+\mathrm{\bar{d}}^{0,0}_{\underline{\End2}_{B}(K)}
(\mathrm{\bar{d}}^{0,1}_{\underline{\End2}_{B}(K)}(\varphi|a)
)=0;\label{compdx}
\end{eqnarray}
\eqref{compsx} is equivalent to
\[
\mathcal
L_A(sa)\circ\partial_{\underline{\End2}_{B}(K)}(\varphi)+(-1)^{|a|-1}\partial_{\underline{\End2}_{B}(K)}(\mathcal
L_A(sa)\circ\varphi)=0;
\]
\eqref{compdx} gives
\[
\partial_{\underline{\End2}_{B}(K)}(\varphi\circ\mathcal
L_A(sa))-\partial_{\underline{\End2}_{B}(K)}(\varphi)\circ\mathcal
L_A(sa)=0.
\]
Both relations are satisfied by checking that
\[
\partial _{\underline{\End2}_{B}(K)}(\mathcal L_A(sa))=0.
\]

\begin{itemize}
 \item {\it $\mathrm L_A$ descends to a morphism of
$A_{\infty}$-$A$-$A$-bimodules.}
\end{itemize}
We  prove that 
\begin{eqnarray}
\mathrm L_A\circ \mathrm {\tilde{d}}_A=\mathrm
d_{\underline{\End2}_{B}(K)} \circ \mathrm L_A. \label{rell}
\end{eqnarray}
We check \eqref{rell} on strings
$(a_1|\dots|a_k|\bar{a}|\tilde{a}_1|\dots|\tilde{a}_l)\in A[1]^{\otimes
k+l+1}$ and
$(1|b_1|\dots|b_q)\in K[1]\otimes B[1]^{\otimes
q}$, for any $k,l,q \geq 0$. If $k,l\geq 1$, the l.h.s of
\eqref{rell} is
\begin{eqnarray*}
&&\sum_{j=1}^k(-1)^{\sum_{i=1}^{j-1}(|a_i|-1)}
\mathrm d_{K}^{k+l,q}(a_1|\dots|a_{j-1},\mathrm
d_A(a_j|a_{j+1})|a_{j+2}|\dots|a_k|\bar{a}|\tilde{a}_1|\dots|\tilde{a}_l|1|b_1|\dots|b_q)+\\
&&\mathrm
d_{K}^{k+l,q}(-1)^{\sum_{i=1}^{k-1}(|a_i|-1)}(a_1|\dots|a_{k-1},\mathrm
d_A(a_k|\bar{a})|\tilde{a}_1|\dots|\tilde{a}_l|1|b_1|\dots|b_q)+\\
&&\mathrm d_{K}^{k+l,q}(-1)^{\sum_{i=1}^{k}(|a_i|-1)}
(a_1|\dots|a_{k},\mathrm
d_A(\bar{a}|\tilde{a}_1)|\tilde{a}_2|\dots|\tilde{a}_l|1|b_1|\dots|b_q)+\\
&&\sum_{j=1}^k(-1)^{\sum_{i=1}^k(|a_1|-1)+|\bar{a}|-1+\sum_{i=1}^{j-1}(|\tilde{a}_i|-1)}\\
&&\mathrm
d_{K}^{k+l,q}(a_1|\dots|a_k|\bar{a}|\tilde{a}_1|\dots|\tilde{a}_{j-1},\mathrm
d_A(\tilde{a}_j|\tilde{a}_{j+1})|\tilde{a}_{j+2}|\dots
|\tilde{a}_l|1|b_1|\dots|b_q).
\end{eqnarray*}

As $\mathrm L_A(a_1|\dots|a_n)=s\circ\mathcal L_A(a_1|\dots|a_n)$ then
the r.h.s. of \eqref{rell} is
given by the terms
\begin{eqnarray*}
&&\bar{\mathrm d}_{\underline{\End2}_{B}(K)}^{0,0}(\mathrm
L_A(a_1|\dots|a_k|\bar{a} |
\tilde{a}_{1}|\dots|\tilde{a}_l))(1|b_1|\dots|b_q)=(-1)^{\sum_{i=1}^k(|a_i|-1)+|\bar{a}|-1+\sum_{i=1}^l(|a_i|-1)+1}\\
&&\left(\mathrm d_K^{k+l+1,q-1}(a_1|\dots|a_k|\bar{a} |
\tilde{a}_{1}|\dots|\tilde{a}_l,\mathrm
d_K^{0,1}(1|b_1)|b_2|\dots|b_q)+\right.\\
&&\left. \sum_{j=1}^q(-1)^{1+\sum_{i=1}^{j-1}(|b_i|-1)}
\mathrm d_K^{k+l+1,q-1}(a_1|\dots|a_k|\bar{a} |
\tilde{a}_{1}|\dots|\tilde{a}_l|1|b_1|\dots|b_{j-1},\mathrm
d_B^{2}(b_j|b_{j+1})|b_2|\dots|b_q)\right)+\\
&&-\mathrm d_K^{0,1}(\mathrm d_K^{k+l+1,q-1}(a_1|\dots|a_k|\bar{a} |
\tilde{a}_{1}|\dots|\tilde{a}_l|1|b_1|b_2|\dots|b_{q-1})|b_q),
\end{eqnarray*}
the sum over $k'\in\{1,\dots,k\}$ of terms
\begin{eqnarray*}
&&\bar{\mathrm
d}_{\underline{\End2}_{B}(K)}^{k',0}(a_1|\dots|a_{k'},\mathrm
L_A(a_{k'+1}|\dots|a_k|\bar{a}|\tilde{a}_1|\dots|\tilde{a}_l)
)(1|b_1|\dots|b_q)=\\
&&\sum_{q'=1}^{q} (-1)^{\sum_{i=1}^{k'}(|a_i|-1)+1}
\mathrm d_K^{k',q-q'}(a_1|\dots|a_{k'},\\
&&\mathrm d_K^{k-k'+l+1,q'}(a_{k'+1}|\dots|a_k|\bar{a}|\tilde{a}_1|\dots|
\tilde{a}_l|1|b_1|\dots|b_{q'})|b_{q'+1}|\dots|b_q),
\end{eqnarray*}


and the sum over $l'\in\{0,\dots,l-1\}$ of terms

\begin{eqnarray*}
&&\bar{\mathrm d}_{\underline{\End2}_{B}(K)}^{0,l-l'}(\mathrm
L_A(a_1|\dots|a_k|\bar{a}|\tilde{a}_1|\dots|\tilde{a}_{l'})|
\tilde{a}_{l'+1}|\dots|\tilde{a}_l)(1|b_1|\dots|b_q)=\\
&&\sum_{q'=0}^{q}
(-1)^{\sum_{i=1}^{k}(|a_i|-1)+|\bar{a}|-1+\sum_{i=1}^{l'}(|\tilde{a}_i|-1)+1}
\mathrm
d_K^{k+l'+1,q-q'}(a_1|\dots|a_{k}|\bar{a}|\tilde{a}_1|\dots|\tilde{a}_{l'},\\
&&\mathrm d_K^{l-l',q'}(\tilde{a}_{l'+1}|\dots|
\tilde{a}_l|1|b_1|\dots|b_{q'})|b_{q'+1}|\dots|b_q),
\end{eqnarray*}
i.e. those contributions
in the r.h.s. of \eqref{rell} corresponding to the right
actions on elements in $\underline{\End2}_{B}(K)$.

Moving the terms on the r.h.s. of \eqref{rell} (note the overall $-1$
sign) to the l.h.s , we get that \eqref{rell}
are equivalent to
\[
\mathrm d_K^2(a_1|\dots|a_k|\bar{a}|
\tilde{a}_1|\dots|\tilde{a}_{l}|1|b_1|\dots|b_q)=0.
\]
The other cases, i.e. $k=0, l\geq 1$, $k\geq 1, l=0$ and $k=l=0$ are a
trivial sign check. We are done.

\section{Proof of thm.~\ref{Thm29}}

The proof of thm.~\ref{Thm29} is shown in detail. We note that all the proof is based on checking the commutativity of diagrams 
in which objects belonging the classes $\mathcal S_1$, $\mathcal S'_1$ appear (see below). Commutativity of the other diagrams follows from
these two special case. Moreover, we do not need to perform any explicit computation; we just need to apply the definition of the 
$A_\infty$-morphisms we introduced in section~\ref{homee}.
\begin{proof}
{\it On objects in $\mathcal
S_1$, $\mathcal{S}^{'}_1$.}
We introduce
$\mathcal S_1=\{ A\langle i\rangle[n], i,n\in\mathbb Z\}$, and
$\mathcal S'_1=\{K\langle i\rangle[n], i,n\in\mathbb Z\}$ .
By definition of the functors $(\mathcal F,\mathcal G)$ and by
proposition~\ref{aah!}  and \ref{hell} we get
$\mathcal G(\mathcal F(X))\simeq X$, $\mathcal F(\mathcal G(Y))\simeq Y$, for every $X\in\mathcal S_1$, $Y\in\mathcal S'_{1}$, with
$\mathcal F(X)\in\mathcal S'_1$ for every $X\in\mathcal S_1$, and
$\mathcal G(Y)\in\mathcal S_{1}$ for every $Y\in\mathcal S'_{1}$.

{\it On morphisms of objects in $\mathcal S_1$.}

We want to prove that
\begin{eqnarray*}
\begin{CD}
X @> f >> Y \\
@VV \varphi_X V @VV \varphi_Y V \\
\mathcal G(\mathcal F(X)) @> \mathcal G(\mathcal F(f)) >>\mathcal
G(\mathcal F(Y))\\
\end{CD}\label{diag1}
\end{eqnarray*}
commutes, for every $X$ and $Y$ in $\mathcal S_1$ and natural
quasi-isomorphisms $\varphi_X$, $\varphi_Y$ in $\DD^{\infty}(A)$.

Let $f: A\langle k_1\rangle [i]\rightarrow A\langle k_2\rangle[j]$ be a
morphism in $\DD^{\infty}(A)$, for
any $k_1,k_2,i,j\in\mathbb Z$. As usual
$\bar{f}_n: (A\langle k_1\rangle[i])[1]\otimes A[1]^{\otimes
n}\rightarrow (A\langle k_2\rangle[j])[1]$ denotes
its $n$-th Taylor component, for $n\geq 0$.
As $A[1]$ is concentrated in cohomological degree $-1$ and the morphism
$f$ is of bidegree $(0,0)$, then
$\bar{f}_n\neq 0$ if and only if $n=j-i$, i.e.
there exists one and only one non trivial Taylor component, if $j-i\geq 0$.
We denote by $\bar{f}_n=s\circ f_n\circ s^{-1}$ its desuspension.
If $n=0$, then $f_0$ is a right $A$-linear map.

If $X=A\langle k_1\rangle [i]$ and $Y=A\langle k_2\rangle [j]$ we need to
check the commutativity
of the diagram

\begin{eqnarray}
\begin{CD}
A\langle k_1\rangle[i] @>f_{j-i}>>
A\langle k_2\rangle[j]\\
@VV \mathcal{V}_1 V @VV \mathcal{V}_2 V \\
(A\underline{\otimes}_A A)\langle k_1\rangle[i] @.
(A\underline{\otimes}_A A)\langle k_2\rangle[j]\\
@VV \mathcal{V}_3 V @VV
\ \mathcal{V}_4 V \\
(A\underline{\otimes}_A \underline{\End2}_B{K})\langle k_1\rangle[i] @.
(A\underline{\otimes}_A \underline{\End2}_B{K})\langle k_2\rangle[j] @.\\
@VV \mathcal{V}_5 V @VV \mathcal{V}_6 V \\
(A\underline{\otimes}_A (K\underline{\otimes}_B\underline{K}))\langle
k_1\rangle[i] @.
(A\underline{\otimes}_A (K\underline{\otimes}_B\underline{K}))\langle
k_2\rangle[j]\\
@VV || V @VV || V \\
A\langle k_1\rangle[i]\underline{\otimes}_A
(K\underline{\otimes}_B\underline{K})
@>\mathcal G(\mathcal F(f))>>A\langle
k_2\rangle[j]\underline{\otimes}_A(K\underline{\otimes}_B\underline{K})\\
\end{CD}\label{diag2}
\end{eqnarray}
Let us describe it in some detail. We give the definitions of the maps
$\mathcal V_i$, up to suspensions and desuspensions w.r.t. the
cohomological and internal degree.
The strict quasi-isomorphism $\mathcal{V}_1$ is induced by the
$A_\infty$ quasi-isomorphism $\Phi : A\rightarrow A\underline{\otimes}_A A$, described in lem.~\ref{Lem2}.
A similar formula holds true for $\mathcal{V}_2$. The quasi-isomorphism $\mathcal{V}_3$ is given by
$\mathcal{V}_3= 1\otimes \mathrm{L}_A$ and the morphism $\mathcal{V}_4$ is defined similarly.
$\mathcal V_5$ (and similarly $\mathcal V_6$) is
\[
\mathcal V_5=1\otimes \mathcal T
\]
where $\mathcal T$ is the homotopy equivalence of
$A_\infty$-$A$-$A$-bimodules $\mathcal T: \underline{\End2}_B(K)\rightarrow
\Endd_B(\mathcal{B}_B(K))\rightarrow K\underline{\otimes}_B\underline{K}$
given in prop.~\ref{chiave} and prop.~\ref{prop12}.

Let us prove commutativity of \eqref{diag2}. The morphism $f$ has a
unique non trivial Taylor component $f_n$,
for $n=j-i\geq 0$. For any
\[
(a,a_1,\dots,a_{n'})\in A\otimes A^{\otimes n'},
\]
we distinguish the following cases.
\begin{itemize}
\item $n'<n$. Going east and then south in \eqref{diag2} we get $0$;
going south-east, instead, we arrive at
\[
\mathcal G(\mathcal F(f))( a,a_1,\dots,a_{n'},\mathcal T(\mathrm L_A(1)))=0,
\]
because $\mathrm L_A$ is strictly unital (here we write $\mathrm
L_A(1)=\mathrm L^{1}_A(1)$ );
all Taylor components $\mathrm L^{m'+m''+1}_A(\dots|1|\dots)$ such that
with $1\leq m'+m''$ identically vanish.

\item $n'=n$. Going east and then south in \eqref{diag2} we arrive at
\begin{eqnarray}
f_{n}(a,a_1|\dots|a_{n})\otimes \mathcal T(\mathrm L_A(1)),\label{aswell}
\end{eqnarray}
as
\[
\mathcal V_2(f_{n}(a,a_1,\dots,a_{n}))=f_{n}(a,a_1,\dots,a_{n})\otimes
1\in (A\otimes A)^{0}_r\subset
(A\underline{\otimes}_A A)^0_r,
\]
denoting by $r\geq 0$ the internal degree of the string
$(a,a_1,\dots,a_{n})$.
Going south-east in \eqref{diag2} we arrive at \eqref{aswell} as well.
In fact
\begin{eqnarray*}
\mathcal V_1(a,a_1,\dots,a_n)&=&(a,a_1|\dots
|a_n,1)+\sum_{n'=0}^{n-1}(\mathcal V^{n'}_1(a,a_1,\dots,a_{n'})\otimes
a_{n'+1},\dots,a_n)\in A\underline{\otimes}_A A,
\end{eqnarray*}
and
\[
\mathcal G(\mathcal F(f))(\mathcal V_5(\mathcal V_3(\mathcal
V_1(a,a_1|\dots |a_n,1)))=
\mathcal G(\mathcal F(f))((a,a_1|\dots|a_n)\otimes\mathcal T(\mathrm
L_A(1)));
\]
this is because $\mathrm L_A$ is strictly unital and $\mathcal
G(\mathcal F(f))$ strict.
The diagram \eqref{diag2} commutes.

\item $n'>n$. Going east and then south in \eqref{diag2} we arrive at
\[
\mathcal G(\mathcal F(f))(
a,a_1,\dots,a_{n},a_{n+1},\dots,a_{n'},\mathcal T(\mathrm L_A(1)))=
f_n( a,a_1,\dots,a_{n}),a_{n+1},\dots,a_{n'},\mathcal T(\mathrm L_A(1)),
\]
as $\mathcal G(\mathcal F(f))$ is strictly unital.
\end{itemize}

In summary, \eqref{diag2} commutes.

{\it Induction: $\mathcal S_r$}

We denote by
\[
\mathcal S_r= \underbrace{\mathcal S_1\star\dots\star\mathcal
S_1}_{r-\mbox{times}},
\]
the $r^{th}$ extensions of $\mathcal S_1$, for every $r\geq 1$.
We want to prove the commutativity of any diagram
\begin{eqnarray}
\begin{CD}
X @> f_1 >> Y \\
@VV \varphi_X V @VV \varphi_Y V \\
\mathcal G(\mathcal F(X)) @> \mathcal G(\mathcal F(f_1)) >>\mathcal
G(\mathcal F(Y))\\
\end{CD}\label{paglia}
\end{eqnarray}
with $X,Y\in\triang^{\infty}_A(A)$ 
and $\varphi_X$, $\varphi_Y$ isomorphisms in $\DD^{\infty}(A)$. 
If $X,Y\in \triang^{\infty}_A(A)$, then by definition there exist
$r,r'\geq 1$ s.t. $X\in\mathcal S_r$ , $Y\in\mathcal S_{r'}$ and
exact triangles
\[
X_1\rightarrow X\rightarrow X'_{r-1}\stackrel{f}{\rightarrow} X_1[1]
\]
and
\[
Y_1\rightarrow Y\rightarrow Y'_{r'-1}\stackrel{g}{\rightarrow} Y_1[1]
\]
in $\DD^{\infty}(A)$ for some morphisms $f$ and $g$ and
$X_1,Y_1\in\mathcal S_1$, $X'_{r-1}\in\mathcal S_{r-1}$ and
$Y'_{r'-1}\in\mathcal S_{r'-1}$. They are isomorphic to the exact triangles
$X_1\stackrel{i_{X_1}}{\rightarrow}X_1\oplus X'_{r-1}\stackrel{p_{X'_{r-1}}}{\rightarrow}X'_{r-1} \stackrel{f}{\rightarrow} X_1[1]$ 
and $Y_1\stackrel{i_{Y_1}}{\rightarrow}Y_1\oplus Y'_{r'-1}\stackrel{p_{Y'_{r'-1}}}{\rightarrow}Y'_{r'-1} \stackrel{f}{\rightarrow} Y_1[1]$ 
for some   isomorphisms
$\rho_{X}: X\rightarrow X_1\oplus X'_{r-1}$ in $\triang^{\infty}_A(A)$ and
$\rho_{Y}:Y\rightarrow Y_1\oplus Y'_{r'-1}$ in $\triang^{\infty}_A(A)$.
Commutativity of \eqref{paglia}
is equivalent to the commutativity of
\begin{eqnarray}
\begin{CD}
X_1\oplus X'_{r-1} @> \tilde{f}_1 >> Y_1\oplus Y'_{r'-1}\\
@VV \varphi_{X_1\oplus X'_{r-1}} V @VV \varphi_{Y_1\oplus Y'_{r'-1}} V \\
\mathcal G(\mathcal F(X_1\oplus X'_{r-1})) @> \mathcal G(\mathcal
F(\tilde{f}_1)) >>\mathcal
G(\mathcal F(Y_1\oplus Y'_{r'-1}))\\
\end{CD} \label{paglia2}
\end{eqnarray}
where
\begin{eqnarray*}
\tilde{f}_1= \rho_Y\circ f_1\circ \rho^{-1}_X, & \varphi_X=\mathcal G(\mathcal F(\rho_X))^{-1}\circ\varphi_{X\oplus X'_{r-1}}\circ\rho_X
\end{eqnarray*}
and similarly for $\varphi_Y$.
Let us discuss  \eqref{paglia2} and the isomorphisms
\begin{eqnarray*}
&&\varphi_{X_1\oplus X'_{r-1}}: X_1\oplus X'_{r-1}\rightarrow \mathcal
G(\mathcal F(X_1\oplus X'_{r-1}))\equiv
\mathcal G(\mathcal F(X_1))\oplus \mathcal G(\mathcal F(X'_{r-1})),\\
&&\varphi_{Y_1\oplus Y'_{r'-1}}: X_1\oplus Y'_{r'-1}\rightarrow \mathcal
G(\mathcal F(Y_1\oplus Y'_{r'-1}))\equiv
\mathcal G(\mathcal F(Y_1))\oplus \mathcal G(\mathcal F(Y'_{r'-1}));
\end{eqnarray*}
we distinguish two cases.
\begin{itemize}
\item If $r=2$, and $r'=2$, then such isomorphisms are simply
\end{itemize}
\begin{eqnarray*}
\varphi_{X_1\oplus X'_{1}}&:=&\varphi_{X_1}\oplus\varphi_{X'_{1}},\\
\varphi_{Y_1\oplus Y'_{1}}&:=&\varphi_{Y_1}\oplus\varphi_{Y'_{1}};
\end{eqnarray*}
we have
$\varphi_{X_1\oplus X'_{1}}\circ\mathrm d_{X_1\oplus X'_{1}}
=\mathrm d_{\mathcal G(\mathcal F(X_1\oplus X'_{1}))}\circ
\varphi_{X_1\oplus X'_{1}}$, i.e.
\begin{eqnarray}
\varphi_{X_1\oplus X'_{1}}\circ\left(\begin{array}{cc}
\mathrm d_{X_1} & -f \\
0 & \mathrm d_{X'_{1}}
\end{array}\right)=
\left(\begin{array}{cc}
\mathrm d_{ \mathcal G(\mathcal F(X_1))} & - \mathcal G(\mathcal F(f)) \\
0 & \mathrm d_{ \mathcal G(\mathcal F(X'_{1})) }
\end{array}\right)\circ \varphi_{X_1\oplus X'_{1}}\label{R'lyeh}
\end{eqnarray}
as $\varphi_{X_1}$ and $\varphi_{X'_{1}}$ are morphisms in
$\triang^{\infty}_A(A)$ such that
$\mathcal G(\mathcal F(f))\circ \varphi_{X'_1}=\varphi_{X_1}\circ f$
holds true: we proved this last equality
in the previous subsection, for any $X_1,X'_1\in\mathcal S_1$ and
morphism $f: X'_1\rightarrow X_1[1]$.
Similar considerations hold for $\varphi_{Y_1\oplus Y'_{1}}$.
Moreover $\varphi_{X_1\oplus X'_{1}}$ and $\varphi_{Y_1\oplus Y'_{1}}$
are homotopy equivalences
as $\varphi_{X_1}$ $\varphi_{X'_{1}}$, $\varphi_{Y_1}$ and
$\varphi_{Y'_1}$ are homotopy equivalences.
We recall that
\[
\varphi_{X}: X\rightarrow \mathcal G(\mathcal F(X))
\]
is explicitly given (up to suspensions and desuspensions) by
\begin{eqnarray}
\varphi_{X}(x|a_1|\dots|a_n)=\sum_{n'=0}^n ((x,a_1|\dots|a_{n'}|\mathcal
T(\mathrm L_A(1))) |a_{n'+1}|\dots|a_n),\label{dorma}
\end{eqnarray}
for any $X\in\mathcal S_1$. We are left to prove
\begin{eqnarray}
\varphi_{Y_1\oplus Y'_{1}} \circ \tilde{f}_1=\mathcal G(\mathcal
F(\tilde{f}_1))\circ \varphi_{X_1\oplus X'_{1}}.\label{tognazzi}
\end{eqnarray}
The strategy is clear: we use the same techniques introduced in the
previous subsection, for the $\mathcal S_1$ case.
We just need to consider any string
\[
((x_1\oplus x'_1)| a_1|\dots|a_n)\in (X_1\oplus X'_{1})[1]\otimes
A[1]^{\otimes n};
\]
\eqref{tognazzi} is equivalent to
\begin{eqnarray}
&&\sum_{n'\geq 0} \varphi_{Y_1\oplus Y'_{1}}
(\tilde{f}_{1,n'}((x_1\oplus x'_1)|
a_1|\dots|a_{n'})|a_{n'+1}|\dots|a_n))=\nonumber\\
&&\sum_{n'\geq 0}\left[ \varphi_{Y_1}(\tilde{f}^{Y_1}_{1,n'}((x_1\oplus
x'_1)| a_1|\dots|a_{n'})|a_{n'+1}|\dots|a_n))\oplus\nonumber\right.\\
&&\left.\varphi_{Y'_1}( \tilde{f}^{Y'_1}_{1,n'}((x_1\oplus x'_1)|
a_1|\dots|a_{n'})|a_{n'+1}|\dots|a_n))\right]\stackrel{!}{=}\nonumber\\
&&\mathcal G(\mathcal F(\tilde{f}_1))(\varphi_{X_1}(x_1|
a_1|\dots|a_{n}))\oplus
\mathcal G(\mathcal F(\tilde{f}_1))(\varphi_{X'_1}(x'_1|
a_1|\dots|a_{n}))=\nonumber\\
&&\mathcal G(\mathcal F(\tilde{f}_1))(x_1| a_1|\dots|a_{n},\mathcal
T(\mathrm L_A(1)))\oplus
\mathcal G(\mathcal F(\tilde{f}_1))( x'_1| a_1|\dots|a_{n},\mathcal
T(\mathrm L_A(1))),
\label{corelli}
\end{eqnarray}
where $\tilde{f}^{Y_1}_{1,n'}$. resp. $\tilde{f}^{Y'_1}_{1,n'}$ denotes
the projection of $\tilde{f}_{1,n'}$ onto $Y_1$, resp. $Y'_1$
. In the last equality in \eqref{corelli} we have used the definition
of $\varphi_{X_1}$, $\varphi_{X'_1}$, following \eqref{dorma}; note that
$\mathcal G(\mathcal F(\tilde{f}_1))$ is a strict morphism; in fact
$\mathcal G(\mathcal F(\tilde{f}_1))=\tilde{f}_1\otimes 1$.

By definition \eqref{dorma}, $\varphi_{Y_1}$
and $\varphi_{Y'_{1}}$ leave the contributions
$\tilde{f}^{Y_1}_{1,n'}((x_1\oplus x'_1)| a_1|\dots|a_{n'})\in Y_1$ and
$\tilde{f}^{Y'_1}_{1,n'}((x_1\oplus x'_1)| a_1|\dots|a_{n'})\in Y'_1$ in
\eqref{corelli}
unchanged: then \eqref{tognazzi} follows and commutativity of
\eqref{paglia} is proven.

\begin{itemize}
\item If $r\geq 3$ or $r'\geq 3$ in \eqref{paglia2}, one needs to
further decompose

the objects $X'_{r-1}$ and $Y'_{r'-1}$ using the above techniques, i.e.
introducing suitable exact triangles, arriving
at the isomorphisms
\end{itemize}
\begin{eqnarray}
\rho'_X: X\rightarrow \rightarrow X_1\oplus
X''_1\oplus\dots X''_{r-1}, & \rho'_Y: Y\rightarrow Y_1\oplus
Y''_1\oplus\dots Y''_{r'-1};\label{dec}
\end{eqnarray}
with $X_1,X''_1,\dots,X''_{r-1}, Y_1,Y''_1,\dots,Y'_{r'-1}\in \mathcal
S_1$. We have reduced our problem to a finite direct sum of the $\mathcal S_1$
case. There is no substantial difference with the $r=2$, $r'=2$ case, both
conceptually and computationally. We conclude that $\mathcal G\circ\mathcal F\simeq 1$ on
$\triang^{\infty}_A(A)$ .

{\it On morphisms of objects in $\mathcal {S}'_1$.}

Let us consider $\mathcal{S}'_1$ and a strictly unital $A_{\infty}$-morphism
$g: K\langle i \rangle[l]\rightarrow K\langle j\rangle [r] $ with Taylor
components
$\bar{g}_n: K\langle i \rangle[l+1]\otimes B[1]^{\otimes n}\rightarrow
K\langle j\rangle[r+1]$, for $n\geq 0$.
Once again, as $K\langle j\rangle[r+1]$ is concentrated in bidegree
$(-r-1,-j)$ and $g$ is of bidegree $(0,0)$, then $\bar{g}_n=0$ for
$n \neq r-l+j-i$: there exists one and only one non trivial
component $\bar{g}_n$. We check the commutativity of any diagram

\begin{eqnarray}
\begin{CD}
X @> g >> Y \\
@VV \psi_X V @VV \psi_Y V \\
\mathcal F(\mathcal G(X)) @> \mathcal F(\mathcal G(g)) >>\mathcal
F(\mathcal G(Y)) \\
\end{CD} \label{diag3}
\end{eqnarray}

with $X$ and $Y$ in $\mathcal{S}'_1$ and $\psi_X$, $\psi_Y$ natural
quasi-isomorphisms in $\DD^{\infty}(B)$. We need some preliminary results to prove this statement.

We recall that the $A_\infty$-$B$-$A$-bimodules $\underline{K}$,
$\overline{K}$ are such that $A_\infty$ quasi-isomorphisms of strictly
unital $A_\infty$-bimodules $\nu_A:A\rightarrow K\underline{\otimes}_B\underline{K}$, and
$\nu_B: B\rightarrow \overline{K}\underline{\otimes}_A K$ exist . Introducing the functor
\begin{eqnarray}
\bar{\mathcal{G}}:\DD^{\infty}(B)\rightarrow
\DD^{\infty}(A), M\mapsto
\bar{\mathcal{G}}(M):=M\underline{\otimes}_B\overline{K} , \label{Gbar}
\end{eqnarray}
we prove that

\begin{Lem}\label{AHA}

\begin{itemize}
\item $\bar{\mathcal{G}}(M)\simeq \mathcal G(M)$ in $\DD^{\infty}(A)$
\footnote{More precisely, the quasi-isomorphisms are all of strictly
unital $A_\infty$-$A$-$A$-bimodules.}, for any $M\in\DD^{\infty}(B)$.
\item The functor $\bar{\mathcal{G}}$ is exact w.r.t. the triangulated
structures on
$\DD^{\infty}(B)$ and $\DD^{\infty}(A)$;
\item $\bar{\mathcal{G}}(M)\in \triang^{\infty}_A(A)$, for any object
$M$ in $\triang^{\infty}_B(K)$.

\end{itemize}

\end{Lem}
\begin{proof}
Let us consider the following diagram:
\begin{eqnarray}
\begin{diagram}
& & & & (M\underline{\otimes}_B \overline{K})\underline{\otimes}_A
(K\underline{\otimes}_B\underline{K}) & & & & \\
& & & \ruTo^{1\otimes \nu_A} & & \luTo^{1\otimes\nu_B\otimes 1} & & & \\
& & (M\underline{\otimes}_B \overline{K})\underline{\otimes}_A A & & & &
M\underline{\otimes}_B (B\underline{\otimes}_B\underline{K})& & \\
& \ruTo^{\Phi_{M\underline{\otimes}_B \overline{K}}} & & & & & &
\luTo^{1\otimes\Phi_{\underline{K}}} & \\
M\underline{\otimes}_B \overline{K}=\bar{\mathcal{G}}(M) & & & & & & & &
\mathcal G(M)=M\underline{\otimes}_B\underline{K}
\end{diagram} \nonumber\\ \label{WTF}
\end{eqnarray}
where
\[
\Phi_{M\underline{\otimes}_B \overline{K}}: M\underline{\otimes}_B
\overline{K}\rightarrow (M\underline{\otimes}_B
\overline{K})\underline{\otimes}_A A
\]
and similarly for $\Phi_{\underline{K}}$
are the maps described in lemma~\ref{Lem2}.
All arrows are quasi-isomorphisms of strictly unital
$A_\infty$-$A$-$A$-bimodules, i.e. homotopy equivalences.
We denote by
\[
\varphi_M: \bar{\mathcal{G}}(M)\rightarrow \mathcal{G}(M)
\]
the above quasi-isomorphism in $\DD^{\infty}(A)$, obtained by inverting
the quasi-isomorphisms
\[
\eta_{M,B}:=1\otimes\nu_B\otimes 1,
\]
and $1\otimes\Phi_{\underline{K}}$; their inverses
exist as quasi-isomorphisms in $\DD^{\infty}(A)$ and $\DD^{\infty}(B)$
are homotopy equivalences. In other words,
\begin{eqnarray}
\varphi_M=(1\otimes\Phi_{\underline{K}})^{-1}\circ \eta^{-1}_{M,B}
\circ (1\otimes \nu_A) \circ (\Phi_{M\underline{\otimes}_B
\overline{K}}):={(\varphi^2_M)}^{-1}\circ\varphi^1_M \label{bigwhite}.
\end{eqnarray}
i.e. $\varphi^1_M$ is the composition on the l.h.s. of
\eqref{WTF}, while $\varphi^2_M$ is the one on the r.h.s.
We also introduce the notation
\begin{eqnarray}
T(M)=(M\underline{\otimes}_B \overline{K})\underline{\otimes}_A
(K\underline{\otimes}_B\underline{K})\cong
M\underline{\otimes}_B (\overline{K}\underline{\otimes}_A
K)\underline{\otimes}_B\underline{K}.\label{dopeshow}
\end{eqnarray}

By definition $\varphi_M$ is not strict.

The first statement of the lemma follows by the very definition of
$\varphi_M$.

The second statement is proved easily using the same techniques showing
that $\mathcal G$
is an exact functor w.r.t. the triangulated structures on
$\DD^{\infty}(B)$ and $\DD^{\infty}(A)$. One the third statement; choosing $M=K$, then \eqref{WTF} implies
$\bar{\mathcal{G}}(K)\simeq A$, or
$\bar{\mathcal{G}}(X)\in\mathcal S_1=\{ A[i]\langle j\rangle,
i,j\in\mathbb Z\}$, for any $X\in \mathcal S'_1=\{K[i]\langle j\rangle, i,j\in\mathbb Z\}$.
By definition of $\mathcal S_r$ and $\mathcal S'_r$ we have $\bar{\mathcal{G}}(X)\in\mathcal S_r$,
for any $X\in \mathcal S'_r$.

\end{proof}

\begin{Cor}
 $\mathcal F\circ\bar{\mathcal{G}}\simeq 1$ on $\triang^{\infty}_B(K)$.
\end{Cor}
\begin{proof}

Thanks to lemma~\ref{AHA}, if $M\in \triang^{\infty}_B(K)$, then $\mathcal
G(M)\simeq \bar{\mathcal{G}}(M)$ and
$\mathcal F(\mathcal G(M))\simeq \mathcal F(\bar{\mathcal{G}}(M))$,
as $\mathcal F$ sends quasi-isomorphisms to quasi-isomorphisms.
It is easy to prove that
$\mathcal F\circ\bar{\mathcal{G}}\simeq 1$
on $\triang^{\infty}_B(K)$ by following the subsection ``{\it Induction: $\mathcal S_r$}'' above.
Checking  the commutativity of
\begin{eqnarray*}
\begin{CD}
X @> g >> Y \\
@VV \psi_X V @VV \psi_Y V \\
\mathcal F(\bar{\mathcal{G}}(X)) @> \mathcal F(\bar{\mathcal{G}}(g))
 >>\mathcal
F(\bar{\mathcal{G}}(Y)) \\
\end{CD}
\end{eqnarray*}
for any $X,Y\in\mathcal S'_1$ is immediate; note that $\psi_X$ and $\psi_Y$ are
explicit; in fact
\[
\psi_X: X\rightarrow \mathcal F(\mathcal{\bar{G}}(X)),~~~~~~X\in\mathcal
S'_1
\]
is given by (up to suspensions and desuspensions)
\begin{eqnarray}
\psi_X(x|b_1|\dots|b_n)=\sum_{n'\geq 0} ((x,b_1|\dots|b_n'|\mathcal
N(\mathrm R_B(1)))|b_{n'+1}|\dots|b_n),\label{trovatore}
\end{eqnarray}
where $\mathcal N: \underline{\End2}_A(K)\rightarrow
\overline{K}\underline{\otimes}_A K$ is s.t.
$\mathcal N(\mathrm R_B(1))=\varphi\otimes 1$, with  $\varphi: (A\underline{\otimes}_A K)^0_0\rightarrow K$, and $\varphi(1\otimes 1)=1$.
This follows from the very definition of $\mathrm R_B(1)$; in fact
$\mathrm R_B(1)(a_1|\dots|a_l|1)=0 $
for $l\geq 1$  and  $\mathrm R_B(1)(1)=1\cdot 1=1$.
The relations
\[
\mathcal F(\bar{\mathcal{G}}(g))\circ\psi_X=\psi_Y\circ g,
\]
i.e.
\begin{eqnarray}
(g\otimes 1)\circ\psi_X=\psi_Y\circ g \label{rusticana}
\end{eqnarray}
follow, using \eqref{trovatore}. Then one moves to diagrams in which $X,Y\in\triang^{\infty}_B(K)$ and
any morphism $g:X\rightarrow Y$ in $\DD^{\infty}(B)$ appear;
the proof of commutativity is done as in subsection ``{\it Induction: $\mathcal S_r$}''.
We decompose objects in $\mathcal S'_r$, $r\geq 2$ into direct sums of
objects in $\mathcal S'_1$; as the case for $r=1$ is explicit,
thanks to \eqref{trovatore} and \eqref{rusticana},
then we can repeat {\it verbatim} the considerations in ``{\it Induction: $\mathcal S_r$}'', ending the proof of the equivalence
$\mathcal F\circ\bar{\mathcal{G}}\simeq 1$ on $\triang^{\infty}_B(K)$.

\end{proof}

We are left to prove
\[
\mathcal F \circ \mathcal G\simeq 1
\]
on $\triang^{\infty}_B(K)$ by induction on $\mathcal S'_r$, $r\geq1$; we
begin with the case $r=1$.
Let $X$ and $Y$ be in $\mathcal S'_1$; any diagram
\begin{eqnarray*}
\begin{CD}
X @> g >> Y \\
@VV \rho_X V @VV \rho_Y V \\
\mathcal F(\mathcal{G}(X)) @> \mathcal F(\mathcal{G}(g)) >>\mathcal
F(\mathcal{G}(Y)) \\
\end{CD}
\end{eqnarray*}
can be decomposed into the subdiagrams
\begin{eqnarray}
\begin{CD}
X @> g >> Y \\
@VV \psi_X V @VV \psi_Y V \\
\mathcal F(\bar{\mathcal{G}}(X)) @> \mathcal F(\bar{\mathcal{G}}(g))
 >>\mathcal
F(\bar{\mathcal{G}}(Y)) \\
@VV \mathcal F(\varphi_X) V @VV \mathcal F(\varphi_Y) V \\
\mathcal F(\mathcal{G}(X)) @> \mathcal F(\mathcal{G}(g)) >>\mathcal
F(\mathcal{G}(Y)) \\
\end{CD}\label{lastday}
\end{eqnarray}
where $\psi_X$, $\psi_Y$ are given by \eqref{trovatore}, $\varphi_X$ and
$\varphi_Y$ by \eqref{bigwhite}
and $\rho_X= \mathcal F(\varphi_X) \circ \psi_X$, $\rho_Y=\mathcal
F(\varphi_Y) \circ \psi_Y$ are quasi-isomorphisms.

We have already proved that the upper subdiagram in \eqref{lastday}
commutes;
the lower one commutes if we prove
the commutativity of the
diagram
\begin{eqnarray}
\begin{CD}
\bar{\mathcal{G}}(X) @> \bar{\mathcal{G}}(g)>>\bar{\mathcal{G}}(Y) \\
@VV \varphi_X V @VV \varphi_Y V \\
\mathcal{G}(X) @> \mathcal{G}(g) >>\mathcal{G}(Y) \\
\end{CD}\label{posthuman}
\end{eqnarray}
as $\mathcal F$ is a functor.
Once again, using the definition \eqref{bigwhite} of $\varphi_X$ and
$\varphi_Y$ we decompose \eqref{posthuman} into
\begin{eqnarray*}
\begin{CD}
\bar{\mathcal{G}}(X) @> \bar{\mathcal{G}}(g)>>\bar{\mathcal{G}}(Y) \\
@V (1\otimes \nu_A) \circ (\Phi_{X\underline{\otimes}_B \overline{K}})VV
@VV(1\otimes \nu_A) \circ
(\Phi_{Y\underline{\otimes}_B \overline{K}}) V \\
T(X) @>T(g) >> T(Y)\\
@V(1\otimes\Phi_{\underline{K}})^{-1}\circ \eta^{-1}_{X,B} VV @VV
(1\otimes\Phi_{\underline{K}})^{-1}\circ \eta^{-1}_{Y,B} V \\
\mathcal{G}(X) @> \mathcal{G}(g) >>\mathcal{G}(Y) \\
\end{CD}
\end{eqnarray*}
where the morphisms appear in the definition \eqref{bigwhite} and
$T(X)$ (similarly for $T(Y)$) is defined in \eqref{dopeshow} .

The map $T(g):T(X)\rightarrow T(Y)$ is simply $T(g)=g\otimes 1$.

But
\begin{eqnarray}
(1\otimes \nu_A) \circ (\Phi_{Y\underline{\otimes}_B \overline{K}})
\circ \bar{\mathcal{G}}(g)=T(g)\circ
(1\otimes \nu_A) \circ (\Phi_{X\underline{\otimes}_B
\overline{K}}),\label{somebody}
\end{eqnarray}
as one can easily check just using the definitions of the morphisms; in
fact the identity has to be verified on any string, say
\[
(x,b_1|\dots|b_q,\overline{\phi})|a_1|\dots|a_n)\in (X
\underline{\otimes}_B \overline{K})[1]\otimes \T(A[1]);
\]
The l.h.s. of \eqref{somebody} is equal to (up to signs)
\[
g((x,b_1|\dots|b_q'),b_{q'+1}|\dots|b_q,\overline{\phi})|a_1|\dots|a_n,\nu_A(1)
,
\]
as the morphism $g: X\rightarrow Y$ ($X$ and $Y$ are in $\mathcal S'_1$ )
has only one non trivial Taylor component $g_{q'}$, $q'\geq 0$, due to
the bigrading on $X$ and $Y$;
we already used this fact in the proof of thm.~\ref{Thm29}. The
r.h.s. of \eqref{somebody} gives the same result,
by definition of $T(g)$, which is clearly a strict $A_\infty$-morphism.

On the other hand

\[
\mathcal{G}(g) \circ(1\otimes\Phi_{\underline{K}})^{-1}\circ
\eta^{-1}_{X,B}=
(1\otimes\Phi_{\underline{K}})^{-1}\circ \eta^{-1}_{Y,B} \circ T(g)
\]
holds true if and only if
\begin{eqnarray}
\eta_{Y,B}\circ(1\otimes\Phi_{\underline{K}})\circ \mathcal{G}(g) =
T(g)\circ \eta_{X,B}\circ (1\otimes\Phi_{\underline{K}});\label{tiritera}
\end{eqnarray}
\eqref{tiritera} is easily verified, as we did for \eqref{somebody}; so
\eqref{lastday} commutes.
Let us verify \eqref{tiritera} explicitly, on any string
\[
(x,b_1|\dots|b_q,\underline{\phi})|a_1|\dots|a_n)\in (X
\underline{\otimes}_B \underline{K})[1]\otimes \T(A[1]),~~~n\geq 0.
\]
If $n\geq 1$, as all morphisms in \eqref{tiritera} are strict, then
\eqref{tiritera} is trivially verified.
Note that, by definition, $\Phi_{\underline{K}}$ ``sees'' the left
$B$-module structure on $\underline{K}$, i.e.
$\Phi_{\underline{K}}:\underline{K}\rightarrow B\underline{\otimes}_B
\underline{K}$.
If $n=0$, recalling that $g$ has only one non trivial Taylor component,
say $g_{\bar{q}}$, $\bar{q}\geq 0$, due to the bigradings on $X$ and
$Y$, we arrive at

\[
\sum_{q'=\bar{q}+1}^q
g((x,b_1|\dots|b_{\bar{q}}),b_{\bar{q}+1}|\dots|b_{q'},\nu_B(1),b_{q'+1}|\dots|b_q,\overline{\phi}),
\]
for the l.h.s. of \eqref{tiritera} (up to suspensions and
desuspensions). We recall that
$\nu_B:B\rightarrow \overline{K}\underline{\otimes}_A K$ is strictly
unital, so $\nu_B(b_1|\dots|b_j|1|b_{j+1}|\dots|b_{j'})=0$ if
$j'\geq 1$.
The r.h.s. of \eqref{tiritera} gives the same result, as
$T(g)=g\otimes 1$.

The first step of the induction is proven. If $X\in\mathcal S_r$ and $Y\in\mathcal S_{r'}$, then we prove the commutativity of
\begin{eqnarray}
\begin{CD}
X @> g >> Y \\
@VV \psi_X V @VV \psi_Y V \\
\mathcal F(\bar{\mathcal{G}}(X)) @> \mathcal F(\bar{\mathcal{G}}(g))
 >>\mathcal
F(\bar{\mathcal{G}}(Y)) \\
@VV \mathcal F(\varphi_X) V @VV \mathcal F(\varphi_Y) V \\
\mathcal F(\mathcal{G}(X)) @> \mathcal F(\mathcal{G}(g)) >>\mathcal
F(\mathcal{G}(Y)) \\
\end{CD}\label{lastday2}
\end{eqnarray}
introducing the isomorphisms
\begin{eqnarray*}
\rho_X: X\rightarrow \bigoplus_{i=1}^{r}X_i,~~~\rho_Y: Y\rightarrow
\bigoplus_{j=1}^{r'}Y_j\\
\end{eqnarray*}
in $\DD^{\infty}(B)$, for some $X_1,\dots,X_r,Y_1,\dots,Y_{r'}$ in $\mathcal S_1$. The
considerations that lead us to prove \eqref{R'lyeh} hold here, with due
changes; we are just considering finite direct sums of $A_\infty$-modules and  homotopy equivalences.

Exchanging $A$ and $B$; i.e. using the $A_{\infty}$-$B$-$A$-bimodule
$(\tilde{K},\mathrm d_K)$ and
the new functors $\mathcal F^{''}= \cdot\underline{\otimes}_B \tilde{K}$
and $\mathcal G^{''}=\cdot\underline{\otimes}_A\underline{\tilde{K}}$ we
can prove the equivalence of the triangulated categories
$\triang^{\infty}_B(B) $ and $\triang^{\infty}_A(\tilde{K})$ with the
same techniques introduced above.

\subsubsection{On thick subcategories}
The statement on the thick subcategories follows by additivity of
$\mathcal F$
and $\mathcal G$ ($\mathcal F^{''}$ and $\mathcal G^{''}$ as well),
w.r.t. the coproduct in $\DD^{\infty}(A)$ and
$\DD^{\infty}(B)$, i.e. the direct sum of strictly unital $A_{\infty}$-modules.

More precisely, let $X\in \thick^{\infty}_A(A)$; there exists a $Z\in\triang^{\infty}_A(A)$  s.t.
\[
Z\simeq X\oplus Y,
\]
for some $Y\in \DD^{\infty}(A)$. Let us call such isomorphism $\varphi_X$, i.e. $\varphi_X: Z\rightarrow X\oplus Y$.
It follows that $\mathcal F(X)\in\thick^{\infty}_B(K)$, as $\mathcal F$ is additive and preserves quasi-isomorphisms.
For any morphism $f:X_1\rightarrow X_2$, with $X_1,X_2\in\thick^{\infty}_A(A)$, we want to prove that the diagram
\[
\begin{CD}
X_1 @>f >> X_2  \\
@V\psi_1 VV @VV \psi_2 V \\
\mathcal G(\mathcal F(X_1)) @>\gamma_1 >> \mathcal G(\mathcal F(X_2)) \\
\end{CD}
\]
commutes, for some isomorphisms $\psi_i:X_i\rightarrow\mathcal G(\mathcal F(X_i))$.
All we need is to check the commutativity of the diagram

\begin{eqnarray}
\begin{CD}
@.  X_1 @>f >> X_2 @.  \\
@.  @VV i_1 V @VV i_2 V @.\\
Z_1 @>\varphi_1 >>X_1\oplus Y_1 @>i_2\circ f\circ \pi_1 >> X_2\oplus Y_2 @>\varphi_2^{-1} >> Z_2 \\  
@VV\psi_{Z_1} V  @VV \rho_1 V @VV \rho_2 V @VV\psi_{Z_2} V\\
\mathcal G(\mathcal F(Z_1)) @>\mathcal G(\mathcal F(\varphi_1)) >>
\mathcal G(\mathcal F(X_1\oplus Y_1)) @> 
\mathcal G(\mathcal F(i_2\circ f\circ\pi_1)) >> \mathcal G(\mathcal F(X_2\oplus Y_2)) @>
\mathcal G(\mathcal F(\varphi_2^{-1})) >> \mathcal G(\mathcal F(Z_2)) \\
@.  @VV \mathcal G(\mathcal F(\pi_1)) V @VV \mathcal G(\mathcal F(\pi_2)) V @.\\
@. \mathcal G(\mathcal F(X_1)) @>\mathcal G(\mathcal F(f))>>  \mathcal G(\mathcal F(Y_1)) @.
\end{CD}\label{thick-diag}
\end{eqnarray}
with
\[
\rho_1:= \mathcal G(\mathcal F(\varphi_1))\circ \psi_{Z_1}\circ\varphi_1^{-1}
\]
and similarly for $\rho_2$. The isomorphisms $\psi_{Z_i}$ do exist as $Z_i\in\triang^{\infty}_A(A)$.
The maps $\pi_j :X_j\oplus Y_j\rightarrow X_j$ and $i_j:X_j\rightarrow X_j\oplus Y_j$ are morphisms in $\DD^{\infty}(A)$.
We want to prove that the central square in (\ref{thick-diag}) commutes, i.e.
\begin{eqnarray}
\rho_2\circ (i_2\circ f\circ\pi_1)=\mathcal G(\mathcal F(i_2\circ f\circ\pi_1))\circ \rho_1; \label{gunz}
\end{eqnarray}
clearly 
\[
\begin{CD}
 Z_1 @>\varphi_2^{-1}\circ(i_2\circ f\circ \pi_1)\circ\varphi_1 >>Z_2  \\  
@VV\psi_{Z_1} V  @VV\psi_{Z_2} V\\
\mathcal G(\mathcal F(Z_1)) @>\mathcal G(\mathcal F(\varphi_2^{-1}\circ(i_2\circ f\circ \pi_1)\circ\varphi_1))  >> \mathcal G(\mathcal F(Z_2)) \\
\end{CD}
\]
commutes, as $Z_i\in\triang^{\infty}_A(A)$; in other words
\begin{eqnarray*}
\psi_{Z_2}\circ \varphi_2^{-1}\circ(i_2\circ f\circ \pi_1)\circ\varphi_1&=&
\mathcal G(\mathcal F(\varphi_2^{-1}\circ(i_2\circ f\circ \pi_1)\circ\varphi_1))\circ\psi_{Z_1}\Leftrightarrow\\
\psi_{Z_2}\circ \varphi_2^{-1}\circ(i_2\circ f\circ \pi_1)&=&
\mathcal G(\mathcal F(\varphi_2^{-1}\circ(i_2\circ f\circ \pi_1)\circ\varphi_1))\circ\psi_{Z_1}\circ\varphi_1^{-1}\Leftrightarrow\\
\mathcal G(\mathcal F(\varphi_2)\circ\psi_{Z_2}\circ \varphi_2^{-1}\circ(i_2\circ f\circ \pi_1)&=&
\mathcal G(\mathcal F((i_2\circ f\circ \pi_1)\circ\varphi_1))\circ\psi_{Z_1}\circ\varphi_1^{-1},
\end{eqnarray*}
i.e. \eqref{gunz}. The upper central and the lower central squares in \eqref{thick-diag} clearly commute; the morphisms
\[
\psi_j: X_j\rightarrow  \mathcal G(\mathcal F(X_j)),~~~\psi_j=\mathcal G(\mathcal F(\pi_j))\circ \rho_j\circ i_j
\]
are actually isomorphisms with inverses given by
\[
 \psi_j^{-1}: \mathcal G(\mathcal F(X_j))\rightarrow  X_j,~~~\psi_j^{-1}=\pi_j\circ\rho_j^{-1}\circ\mathcal G(\mathcal F(i_j));
\]
 this last statement is proved by  writing 
explicitly $\rho_j$ and recalling the decompositions \eqref{dec}, for any object in $\triang^{\infty}_A(A)$.
As $\mathcal G(\mathcal F(\cdot))$ is of the form $\mathcal G(\mathcal F(g))=g\otimes 1$, for any morphism $g$ in $\triang^{\infty}_A(A)$,
 the claim follows.  The thick subcategory $\thick^{\infty}_B(K)$ is studied analogously.


\end{proof}
\section{Proof of thm.~\ref{Thm30} }
We show the proof of thm.~\ref{Thm30} in some detail.
 such proof is analogous to the one of thm.~\ref{Thm29}, modulo technical issue due to the presence of ``roofs''.
Once again, all we need is to prove the commutativity of ``easier'' diagrams in which objects of the form $A_\hbar[i]\langle j\rangle$
 and $K_\hbar[n]\langle m\rangle$ appear.

\begin{proof}
We study the exact\footnote{as usual, ``exact'' is w.r.t the
triangulated structures on the derived categories}
functors $\mathcal F_\hbar$ and $\mathcal G_\hbar$ on
$\DD^{\infty}_{tf}(A_{\hbar})$ and $\DD^{\infty}_{tf}(B_{\hbar})$ to prove
that they induce an equivalence of triangulated categories between
$\triang^{\infty}_{A_{\hbar}}(\tilde{K}_{\hbar})$ and
$\triang^{\infty}_{B_{\hbar}}(B_{\hbar})$.

{\it On objects in $\mathcal S_1$, $\mathcal{S}^{'}_1$.}

We introduce
$\mathcal S_1=\{ A_{\hbar}\langle i\rangle[n], i,n\in\mathbb Z\}$, and
$\mathcal S'_1=\{K_{\hbar}\langle i\rangle[n], i,n\in\mathbb Z\}$. By
definition, objects of $\mathcal S_1$ are (all isomorphism classes
of the) objects $ A_{\hbar}\langle i\rangle[n]$ in
$\DD^{\infty}_{tf}(A_\hbar)$ and similarly for $\mathcal S'_1$.
By definition of the functors $(\mathcal F_\hbar,\mathcal G_\hbar)$ and
proposition~\ref{marroni} we have
\[
\mathcal G_\hbar(\mathcal F_\hbar(A_\hbar\langle i\rangle [n] ))\simeq
A_\hbar\langle i\rangle[n]
\]
in $\DD^{\infty}_{tf}(A_{\hbar})$ and
\[
\mathcal F_\hbar(\mathcal G_\hbar(K_\hbar\langle j\rangle [m] ))\simeq
K_\hbar\langle j\rangle [m]
\]
in $\DD^{\infty}_{tf}(B_{\hbar})$, with
$\mathcal F_\hbar(X)\in \mathcal S'_1$ and $\mathcal
G_\hbar(Y)\in\mathcal S_1$, for
every $X\in\mathcal S_1$, $Y\in\mathcal S'_1$.

\subsubsection{ On commutative diagrams in the
derived categories $\DD^{\infty}_{tf}(A_{\hbar})$ and
$\DD^{\infty}_{tf}(B_{\hbar})$. }
To prove the theorem, we need to consider the following general setting.
Let $X_\hbar$, $Y_\hbar$ be objects in
$\triang^{\infty}_{A_{\hbar}}(A_{\hbar})$
and let $f_\hbar:X_\hbar\rightarrow Y_\hbar$ be a morphism in
$\DD^{\infty}_{tf}(A_{\hbar})$.
We want to prove that there exist isomorphisms
\[
\varphi_{X_\hbar} : X_\hbar\rightarrow \mathcal G_\hbar(\mathcal
F_\hbar(X_\hbar)),~~~\varphi_{Y_\hbar} : Y_\hbar\rightarrow
\mathcal G_\hbar(\mathcal F_\hbar(Y_\hbar))
\]
in $\DD^{\infty}_{tf}(A_{\hbar})$, such that
\begin{eqnarray}
\begin{CD}
X_\hbar @> f_\hbar >>Y_\hbar \\
@V \varphi_{X_\hbar} VV @V \varphi_{Y_\hbar} VV \\
\mathcal G_\hbar(\mathcal F_\hbar(X_\hbar)) @> \mathcal G_\hbar(\mathcal
F_\hbar(f_\hbar)) >>\mathcal G_\hbar(\mathcal F_\hbar(Y_\hbar)) \\
\end{CD}\label{megaborg}
\end{eqnarray}
commutes in $\DD^{\infty}_{tf}(A_{\hbar})$.

Let us represent the morphism $f_\hbar$ by the roof
$(s_\hbar,\bar{f}_\hbar)$, i.e.

\begin{diagram}
& & X'_\hbar & & \\
& \ldTo^{s_\hbar} & & \rdTo{\bar{f}_\hbar} & \\
X_\hbar & & & & Y_\hbar \\
\end{diagram}
for some $X'_\hbar$ in $\DD^{\infty}_{tf}(A_{\hbar})$; by definition of
$\triang^{\infty}_{A_{\hbar}}(A_{\hbar})$ (property $(SO)$) we can infer
that $X'_\hbar$ is an object in $\triang^{\infty}_{A_{\hbar}}(A_{\hbar})$:
in fact $s_\hbar$ is an isomorphism in $\DD^{\infty}_{tf}(A_{\hbar})$.
Then the morphism
$\mathcal G_\hbar(\mathcal F_\hbar(X_\hbar))\rightarrow \mathcal
G_\hbar(\mathcal F_\hbar(Y_\hbar))$ is represented by the roof
$(\mathcal G_\hbar(\mathcal F_\hbar(s_\hbar),\mathcal G_\hbar(\mathcal
F_\hbar(\bar{f}_\hbar)) $, i.e.

\begin{diagram}
& & \mathcal G_\hbar(\mathcal F_\hbar(X'_\hbar)) & & \\
& \ldTo^{\mathcal{G}_\hbar(\mathcal{F}_\hbar(s_\hbar))} & &
\rdTo{\mathcal{G}_\hbar(\mathcal{F}_\hbar(\bar{f}_\hbar))} & \\
\mathcal G_\hbar(\mathcal F_\hbar(X_\hbar)) & & & & \mathcal
G_\hbar(\mathcal F_\hbar(Y_\hbar)) \\
\end{diagram}

We are interested in proving also the commutativity of diagrams
\begin{eqnarray}
\begin{CD}
W_\hbar @> g_\hbar >>Z_\hbar \\
@V \varphi_{W_\hbar} VV @V \varphi_{Z_\hbar} VV \\
\mathcal F_\hbar(\mathcal G_\hbar(W_\hbar)) @> \mathcal F_\hbar(\mathcal
G_\hbar(g_\hbar)) >>\mathcal F_\hbar(\mathcal G_\hbar(Z_\hbar)) \\
\end{CD}\label{megaborg2}
\end{eqnarray}
in $\DD^{\infty}_{tf}(B_{\hbar})$, with $W_\hbar$ and $Z_\hbar$ in
$\triang^{\infty}_{B_{\hbar}}(K_{\hbar})$,
$\varphi_{W_\hbar}$, $\varphi_{Z_\hbar}$ isomorphisms in
$\DD^{\infty}_{tf}(B_{\hbar})$ and morphisms
$g_\hbar$ represented by some roof, like the morphisms $f_\hbar$
introduced above.

We need the following lemma, which reduces the problem of commutativity
in the derived categories to the check of certain relations
involving morphisms in the corresponding homotopy categories. We state
the lemma in the case of diagrams of the form \eqref{megaborg};
the other case is analogous.

\begin{Lem}\label{MM1}
Let $X_\hbar$, $Y_\hbar$ be objects in
$\triang^{\infty}_{A_{\hbar}}(A_{\hbar})$ and let us consider a diagram
of the form \eqref{megaborg},
with $f_\hbar$, $\varphi_{X_\hbar}$, $\varphi_{Y_\hbar}$ and $\mathcal
G_\hbar(\mathcal
F_\hbar(f_\hbar))$ as above.
If there exists a quasi-isomorphism
\[
\bar{\varphi}_{X'_\hbar}: X'_\hbar\rightarrow \mathcal G_\hbar(\mathcal
F_\hbar(X'_\hbar))
\]
in $\mathcal H^{tf}_\infty(A_\hbar)$
s.t. for any morphism $g_\hbar: X'_\hbar\rightarrow Y_\hbar$ in
$\mathcal H^{tf}_\infty(A_\hbar)$ the relation
\begin{eqnarray}
\bar{\varphi}_{Y_\hbar}\circ
g_\hbar=\bar{G}_\hbar(\bar{F}_\hbar(g_\hbar))\circ
\bar{\varphi}_{X'_\hbar}\label{Neon}
\end{eqnarray}
holds true, then, representing the isomorphism
\[
\varphi_{X'_\hbar} : X'_\hbar\rightarrow \mathcal G_\hbar(\mathcal
F_\hbar(X'_\hbar)),
\]
in $\DD^{\infty}_{tf}(A_{\hbar})$ by the roof
\begin{diagram}
& & X'_\hbar & & \\
& \ldTo^{1_\hbar} & & \rdTo{\bar{\varphi}_{X'_\hbar}} & \\
X'_\hbar & & & & \mathcal G_\hbar(\mathcal F_\hbar(X'_\hbar)) \\
\end{diagram}
the diagrams of the form \eqref{megaborg} commute.
\end{Lem}

\begin{proof}
Commutativity of \eqref{megaborg} is equivalent to
\begin{eqnarray}
(1,\bar{\varphi}_{Y_\hbar})\circ (s_\hbar,\bar{f}_\hbar) =
(\bar{G}_\hbar(\bar{F}_\hbar(s_\hbar)),\bar{G}_\hbar(\bar{F}_\hbar(\bar{f}_\hbar))
)\circ (1,\bar{\varphi}_{X'_\hbar})\label{battery}
\end{eqnarray}
in $\DD^{\infty}_{tf}(A_{\hbar})$; the l.h.s. reads
\[
(1,\bar{\varphi}_{Y_\hbar})\circ
(s_\hbar,\bar{f}_\hbar)=(s_\hbar\alpha_\hbar,\bar{\varphi}_{Y_\hbar}\beta_\hbar),
\]
with $\beta=\bar{f}_\hbar\alpha_\hbar$, for some roof
\begin{diagram}
& & Z_\hbar & & \\
& \ldTo^{\alpha_\hbar} & & \rdTo{\beta_\hbar} & \\
X'_\hbar & & & & Y_\hbar \\
\end{diagram}
Then
\[
(1,\bar{\varphi}_{Y_\hbar})\circ (s_\hbar,\bar{f}_\hbar)=
(s_\hbar\alpha_\hbar,\bar{\varphi}_{Y_\hbar}\bar{f}_\hbar\alpha_\hbar)=(s_\hbar,\bar{\varphi}_{Y_\hbar}\bar{f}_\hbar)=
(s_\hbar,\bar{G}_\hbar(\bar{F}_\hbar(\bar{f}_\hbar))\bar{\varphi}_{X'_\hbar}),
\]
where in the last equality we used \eqref{Neon} and the second equality
holds true as $\alpha_\hbar$ is a quasi-isomorphism
\footnote{We recall that equality ``='' between roofs is, by definition,
the equivalence relation betweeen them.}; but
\[
(s_\hbar,\bar{G}_\hbar(\bar{F}_\hbar(\bar{f}_\hbar))
\bar{\varphi}_{X'_\hbar}),
\]
i.e. the roof
\begin{diagram}
& & X'_\hbar & & \\
& \ldTo^{s_\hbar} & & \rdTo{
\bar{G}_\hbar(\bar{F}_\hbar(\bar{f}_\hbar))\varphi_{X'_\hbar}} & \\
X_\hbar & & & & \mathcal G_\hbar(\mathcal F_\hbar(Y_\hbar)) \\
\end{diagram}
is equal to

\begin{diagram}
& & & & X'_\hbar & & & & \\
& & & \ldTo^{s_\hbar} & & \rdTo{ \bar{\varphi}_{X'_\hbar}} & & & \\
& & X_\hbar & & & & \mathcal G_\hbar(\mathcal F_\hbar(X'_\hbar)) & & \\
& \ldTo^{1_\hbar} & & \rdTo{\bar{\varphi}_{X_\hbar}} & &
\ldTo^{\bar{G}_\hbar(\bar{F}_\hbar(s_\hbar))}
& & \rdTo^{\bar{G}_\hbar(\bar{F}_\hbar(\bar{f}_\hbar))} & \\
X_\hbar & & & & \mathcal G_\hbar(\mathcal F_\hbar(X_\hbar)) & & & &
\mathcal G_\hbar(\mathcal F_\hbar(Y_\hbar))
\end{diagram}

i.e. the composition
\[
(\bar{G}_\hbar(\bar{F}_\hbar(s_\hbar)),\bar{G}_\hbar(\bar{F}_\hbar(\bar{f}_\hbar))
)\circ (1,\bar{\varphi}_{X'_\hbar}),
\]

which is the r.h.s. of \eqref{battery}, if and only if
\[
\bar{\varphi}_{X_\hbar}\circ
s_\hbar=\bar{G}_\hbar(\bar{F}_\hbar(s_\hbar))\circ \bar{\varphi}_{X'_\hbar}.
\]
This latter is nothing but \eqref{Neon} applied to the quasi-isomorphism
$s_\hbar$.
\end{proof}

In virtue of the above lemma, we prove that diagrams of the form
\eqref{megaborg} commute, for any $X_\hbar$, $Y_\hbar$ in
$\triang^{\infty}_{A_{\hbar}}(A_{\hbar})$, by choosing a representative
for the morphisms and checking the relations \eqref{battery}.

\subsubsection{ On $\mathcal G_\hbar\circ \mathcal F_\hbar\simeq 1$ on
$\triang^{\infty}_{A_\hbar}(A_\hbar)$}
We begin by proving $\mathcal G_\hbar\circ \mathcal F_\hbar\simeq 1$ on
$\triang^{\infty}_{A_\hbar}(A_\hbar)$
on any pair
\[
X_\hbar:= A_\hbar\langle i' \rangle[n'],~~~~ Y_\hbar:=A_\hbar\langle
j\rangle[m]
\]
of objects in $\mathcal S_1$ and any morphism
 $f_\hbar: A_\hbar\langle i' \rangle[n']\rightarrow A_\hbar\langle
j\rangle[m]$ in $\DD^{\infty}_{tf}(A_{\hbar})$ represented by the roof
\begin{eqnarray}
\begin{diagram}
& & A_\hbar\langle i \rangle[n] & & \\
& \ldTo^{s_\hbar} & & \rdTo{\bar{f}_\hbar} & \\
A_\hbar\langle i' \rangle[n'] & & & & A_\hbar\langle
j\rangle[m] \\
\end{diagram}\label{fbasic}
\end{eqnarray}
where $i,i',n,n',j,m\in\mathbb Z$.

Let $\bar{f}^{(l),r}_\hbar: (A\langle i \rangle[n])[1]\tilde{\otimes}
A[1]^{\tilde{\otimes} r}\rightarrow (A \langle j\rangle[m])[1]$
be the $r$-th Taylor component of $f_\hbar^{(l)}$ for any $l,r\geq 0$.
A quick degree analysis (we recall that $A_\hbar$ is concentrated in
cohomological degree $0$) implies that $\bar{f}^{(l),r}_\hbar\neq 0$ if and only if $r=m-n$, i.e.
there exists one and only one non trivial Taylor component of
$f^{(l)}_\hbar$, if $m-n\geq 0$, for any $l\geq 0$.
To prove the commutativity of the diagram
\begin{eqnarray*}
\begin{CD}
A_\hbar\langle i \rangle[n] @> f_\hbar >>A_\hbar\langle
j\rangle[m] \\
@VV \varphi_{A_\hbar\langle i \rangle[n]} V @VV \varphi_{A_\hbar\langle
j\rangle[m]} V \\
\mathcal G_\hbar(\mathcal F_\hbar(A_\hbar\langle i \rangle[n])) @>
\mathcal G_\hbar(\mathcal
F_\hbar(f_\hbar)) >>\mathcal G_\hbar(\mathcal F_\hbar( A_\hbar\langle
j\rangle[m])) \\
\end{CD}
\end{eqnarray*}
is sufficient, thanks to lemma~\ref{MM1} to prove
\begin{eqnarray}
\bar{\varphi}_{A_\hbar\langle j \rangle[m]} \circ\bar{f}_\hbar=\mathcal
G_\hbar(\mathcal
F_\hbar(\bar{f}_\hbar))\circ \bar{\varphi}_{A_\hbar\langle i
\rangle[n]}\label{pira!}
\end{eqnarray}
representing the isomorphisms
\[
\varphi_{A_\hbar\langle k \rangle[l]}: A_\hbar\langle k
\rangle[l]\rightarrow \mathcal G_\hbar(\mathcal F_\hbar(A_\hbar\langle k
\rangle[l]))
\]
in $\DD^{\infty}_{tf}(A_\hbar)$ by 
\begin{eqnarray}
\begin{diagram}
& & A_\hbar\langle j \rangle[m] & & \\
& \ldTo{1_\hbar} & & \rdTo{\bar{\varphi}_{A_\hbar\langle j \rangle[m]}} & \\
A_\hbar\langle j \rangle[m] & & & & \mathcal G_\hbar(\mathcal
F_\hbar(A_\hbar\langle j \rangle[m])) \\
\end{diagram}\label{FUUU}
\end{eqnarray}
for any $k,l\in\mathbb Z$.
The quasi-isomorphisms $\bar{\varphi}_{A_\hbar\langle j \rangle[m]}$ and
$\bar{\varphi}_{A_\hbar\langle i \rangle[n]}$
can be deduced by the diagram \eqref{diag2}, with due changes. Up to
suspensions and desuspensions w.r.t.\ both the cohomological
and internal degree, we have
\begin{eqnarray}
\bar{\varphi}_{A_\hbar\langle i \rangle[n]}=(1\tilde{\otimes}\mathcal
T_\hbar)\circ(1\tilde{\otimes}\mathrm
L_{A_\hbar})\circ\Phi_{A_\hbar},\label{morfismobase}
\end{eqnarray}
where $\Phi_{A_\hbar}: A_\hbar\rightarrow
A_\hbar\underline{\tilde{\otimes}}_{A_\hbar}A_\hbar$ and
$\mathcal T_\hbar: \underline{\End2}_{B_\hbar}(K_\hbar)\rightarrow
K_\hbar\underline{\tilde{\otimes}}_{B_\hbar}\underline{K}_\hbar$
is described in Corollary 5. To check \eqref{pira!} is immediate, once
we recall that $\mathcal G_\hbar(\mathcal
F_\hbar(\bar{f}_\hbar))=\bar{f}_\hbar\tilde{\otimes}1$.

We finish the proof of the equivalence $\mathcal G_\hbar\circ \mathcal
F_\hbar\simeq 1$ on $\triang^{\infty}_{A_\hbar}(A_\hbar)$
considering the general case.
Denoting by
\[
\mathcal S_r= \underbrace{\mathcal S_1\star\dots\star\mathcal
S_1}_{r-\mbox{times}},
\]
the $r^{th}$
extension of $\mathcal S_1$
for every $r\geq 1$, we note that  $\mathcal F_\hbar(X_\hbar)\in\mathcal S'_r$
for every $X_\hbar\in\mathcal S_r$ and $\mathcal G_\hbar(Y_\hbar)\in\mathcal S_{r'}$
for every $Y_\hbar\in\mathcal S'_{r'}$, and $r,r'\geq 1$.\\

Let $X_\hbar$ and $Y_\hbar$ be objects in
$\triang^{\infty}_{A_{\hbar}}(A_{\hbar})$
and $f_\hbar: X_\hbar\rightarrow Y_\hbar$ be a morphism in
$\DD^{\infty}_{tf}(A_{\hbar})$ represented by the roof
\begin{diagram}
& & X'_\hbar & & \\
& \ldTo^{s_\hbar} & & \rdTo{\bar{f}_\hbar} & \\
X_\hbar & & & & Y_\hbar \\
\end{diagram}
It follows that $X'_\hbar$ is an object in
$\triang^{\infty}_{A_{\hbar}}(A_{\hbar})$ as well.
we show that the diagram
\begin{eqnarray}
\begin{CD}
X_\hbar @> f_\hbar >>Y_\hbar \\
@V \varphi_{X_\hbar} VV @V \varphi_{Y_\hbar} VV \\
\mathcal G_\hbar(\mathcal F_\hbar(X_\hbar)) @> \mathcal G_\hbar(\mathcal
F_\hbar(f_\hbar)) >>\mathcal G_\hbar(\mathcal F_\hbar(Y_\hbar)) \\
\end{CD}\label{megaborg!}
\end{eqnarray}
commutes in $\DD^{\infty}_{tf}(A_{\hbar})$. To do so, we introduce the roofs
{\footnotesize
\[
\begin{diagram}
& & X_\hbar & & \\
& \ldTo^{1_\hbar} & & \rdTo{\bar{\varphi}_{X_\hbar}} & \\
X_\hbar & & & & \mathcal G_\hbar(\mathcal F_\hbar(X_\hbar)), \\
\end{diagram}
\qquad\qquad\begin{diagram}
& & Y_\hbar & & \\
& \ldTo^{1_\hbar} & & \rdTo{\bar{\varphi}_{Y_\hbar}} & \\
Y_\hbar & & & & \mathcal G_\hbar(\mathcal F_\hbar(Y_\hbar)), \\
\end{diagram}
\qquad\qquad
\begin{diagram}
& & \mathcal G_\hbar(\mathcal F_\hbar(X'_\hbar)) & & \\
& \ldTo^{ \mathcal G_\hbar(\mathcal F_\hbar(s_\hbar))} & &
\rdTo{\mathcal G_\hbar(\mathcal F_\hbar(\bar{f}_\hbar))} & \\
\mathcal G_\hbar(\mathcal F_\hbar(X_\hbar)) & & & & \mathcal
G_\hbar(\mathcal F_\hbar(Y_\hbar)), \\
\end{diagram}
\]
}
representing $\varphi_{X_\hbar}$, $\varphi_{Y_\hbar}$ and $ \mathcal
G_\hbar(\mathcal F_\hbar(f_\hbar))$,
for some $\bar{\varphi}_{X_\hbar}$ and $\bar{\varphi}_{Y_\hbar}$ still
to define (see below)
On the other hand, $\mathcal G_\hbar(\mathcal
F_\hbar(\bar{f}_\hbar))=\bar{f}_\hbar\tilde{\otimes}1$.
By definition of the triangulated subcategory
$\triang^{\infty}_{A_{\hbar}}(A_{\hbar})$,
there exist $r',r\geq 0$ s.t. $X'_\hbar\in\mathcal S_{r'}$ and
$Y_\hbar\in \mathcal S_r$, i.e. there exist exact triangles
\begin{eqnarray}
X^1_{\hbar}\rightarrow X'_\hbar \rightarrow
\bar{X}^{r'-1}_\hbar\stackrel{g_\hbar}{\rightarrow} X^1_\hbar[1], &
Y^1_\hbar\rightarrow Y_\hbar\rightarrow
\bar{Y}^{r-1}_\hbar\stackrel{h_\hbar}{\rightarrow} Y^1_\hbar[1],
\label{troll}
\end{eqnarray}
in $\DD^{\infty}(A)$ for some morphisms $g_\hbar$ and $h_\hbar$, with
$X^1_\hbar, Y^1_\hbar\in \mathcal S_1$ and
$\bar{X}^{r'-1}_\hbar\in\mathcal S_{r'-1}$, $
\bar{Y}^{r-1}_\hbar\in\mathcal S_{r-1}$.
Let us focus on the first exact triangle in \eqref{troll}; for the
second one the analysis is analogous.
By definition of exact triangles in $\DD^{\infty}_{tf}(A_\hbar)$ (and
$\triang^{\infty}_{A_{\hbar}}(A_{\hbar})$), such exact triangle
is isomorphic in $\DD^{\infty}_{tf}(A_\hbar)$ to a sequence of the form
\begin{eqnarray}
W_\hbar \stackrel{\alpha_\hbar}{\rightarrow}
Z_\hbar\stackrel{\beta_\hbar}{\rightarrow}
R_\hbar\stackrel{\gamma_\hbar}{\rightarrow} W_\hbar[1],
\label{reb}
\end{eqnarray}
where \eqref{reb} is the image under the canonical functor $\mathcal
Q_{A_\hbar}: \mathcal H^{tf}_\infty
(A_\hbar)\rightarrow\DD^{\infty}(A_\hbar)$ of the exact triangle
$ W_\hbar \stackrel{\bar{\alpha}_\hbar}{\rightarrow}
Z_\hbar\stackrel{\bar{\beta}_\hbar}{\rightarrow}
R_\hbar\stackrel{\bar{\gamma}_\hbar}
{\rightarrow} W_\hbar[1]$. In other words, the morphism $\alpha_\hbar$
is represented by the roof
\[
\begin{diagram}
& & W_\hbar & & \\
& \ldTo^{1_\hbar} & & \rdTo{\bar{\alpha}_{\hbar}} & \\
W_\hbar & & & & Z_\hbar; \\
\end{diagram}
\]
and similarly for $\beta_\hbar$, $\gamma_\hbar$.

But \eqref{reb} is isomorphic in $\DD^{\infty}_{tf}(A_\hbar)$ to the
exact triangle
\begin{eqnarray}
W_\hbar \stackrel{i_\hbar}{\rightarrow}W_\hbar\tilde{\oplus}
R_\hbar\stackrel{p_\hbar}{\rightarrow}
R_\hbar\stackrel{\gamma_\hbar}{\rightarrow} W_\hbar[1],
~~~~\mathrm d_{W_\hbar\tilde{\oplus} R_\hbar}=\left(\begin{array}{cc}
\mathrm d_{W_\hbar} & -\bar{\gamma}_\hbar \\
0 & \mathrm d_{R_\hbar} \end{array}\right),
\label{reb2}
\end{eqnarray}
where $i_\hbar$ and $p_\hbar$ are represented by
\[
\begin{diagram}
& & W_\hbar & & \\
& \ldTo^{1_\hbar} & & \rdTo{\bar{i}_{\hbar}} & \\
W_\hbar & & & & W_\hbar\tilde{\oplus} R_\hbar, \\
\end{diagram}\qquad\qquad
\begin{diagram}
& & W_\hbar\tilde{\oplus} R_\hbar & & \\
& \ldTo^{1_\hbar} & & \rdTo{\bar{p}_{\hbar}} & \\
W_\hbar\tilde{\oplus} R_\hbar & & & & R_\hbar, \\
\end{diagram}
\]
with canonical topological inclusion $\bar{i}_\hbar$ and topological
projection $\bar{p}_\hbar$.
In summary, collecting the isomorphisms of the exact triangles so far,
we arrive at the isomorphisms $X^1_{\hbar}\simeq W_\hbar$
and $\bar{X}^{r'-1}_\hbar\simeq R_\hbar$ in
$\DD^{\infty}_{tf}(A_\hbar)$, implying that $W_\hbar\in\mathcal S_1$ and
$R_\hbar\in\mathcal S_{r'-1}$;
the isomorphism
\[
\rho_{X'_\hbar}: (X'_\hbar,\mathrm d_{X'_\hbar})\rightarrow
(W_{\hbar}\tilde{\oplus}R_\hbar,\mathrm d_{W_{\hbar}\tilde{\oplus}R_\hbar}),
\]
in $\DD^{\infty}_{tf}(A_\hbar)$ follows, as well.
Repeating the same analysis for the exact triangle in which $Y_\hbar$
appears-see \eqref{troll}-we get the isomorphism
\[
\rho_{Y_\hbar}:(Y_\hbar,\mathrm d_{Y_\hbar})\rightarrow
(M_{\hbar}\tilde{\oplus}N_\hbar,\mathrm d_{M_{\hbar}\tilde{\oplus}N_\hbar}),
\]
in $\DD^{\infty}_{tf}(A_\hbar)$ for some $M_\hbar\in\mathcal S_1$ and
$N_\hbar\in\mathcal S_{r-1}$.

In virtue of the above isomorphisms in $\DD^{\infty}_{tf}(A_\hbar)$,
\eqref{megaborg!} commutes if and only if
\begin{eqnarray}
\begin{CD}
W_{\hbar}\tilde{\oplus}R_\hbar @> \tilde{f}_\hbar
 >>M_{\hbar}\tilde{\oplus}N_\hbar \\
@V \varphi_{W_{\hbar}\tilde{\oplus}R_\hbar} VV @VV
\varphi_{M_{\hbar}\tilde{\oplus}N_\hbar} V \\
\mathcal G_\hbar(\mathcal F_\hbar(W_{\hbar}\tilde{\oplus}R_\hbar)) @>
\mathcal G_\hbar(\mathcal
F_\hbar(\tilde{f}_\hbar)) >>\mathcal G_\hbar(\mathcal
F_\hbar(M_{\hbar}\tilde{\oplus}N_\hbar)) \\
\end{CD}\label{megaborg!!}
\end{eqnarray}
does, where
\begin{eqnarray}
\tilde{f}_\hbar&=& \rho_{Y_\hbar}\circ
f_\hbar\circ\rho^{-1}_{X_\hbar},\nonumber\\
\varphi_{W_{\hbar}\tilde{\oplus}R_\hbar}&=& \rho_{\mathcal
G_\hbar(\mathcal
F_\hbar(X_\hbar))}\circ\varphi_{X_{\hbar}}\circ\rho^{-1}_{X_\hbar},\nonumber\\
\varphi_{M_{\hbar}\tilde{\oplus}N_\hbar}&=& \rho_{\mathcal
G_\hbar(\mathcal
F_\hbar(Y_\hbar))}\circ\varphi_{Y_{\hbar}}\circ\rho^{-1}_{Y_\hbar}.\label{toad}
\end{eqnarray}
In the sequel we will explicitly define
$\varphi_{W_{\hbar}\tilde{\oplus}R_\hbar}$ and
$\varphi_{M_{\hbar}\tilde{\oplus}N_\hbar}$; thanks to \eqref{toad}
$\varphi_{X_\hbar}$ and $\varphi_{Y_\hbar}$ will be explicit, as well.
As we did in the proof of thm.~\ref{Thm29} in the subsection ``{\it Induction: $\mathcal S_r$}'', 
we need to distinguish two cases: if $r'=r=2$, i.e.
$W_\hbar, R_\hbar, M_{\hbar},N_\hbar\in\mathcal S_1$, we represent
$\varphi_{W_{\hbar}\tilde{\oplus}R_\hbar}$
and $\varphi_{M_{\hbar}\tilde{\oplus}N_\hbar}$ by the roofs
\[
\begin{diagram}
& & W_{\hbar}\tilde{\oplus}R_\hbar & & \\
& \ldTo^{1_\hbar} & &
\rdTo{\bar{\varphi}_{W_\hbar}\tilde{\oplus}\bar{\varphi}_{R_\hbar}} & \\
W_{\hbar}\tilde{\oplus}R_\hbar & & & & \mathcal G_\hbar(\mathcal
F_\hbar(W_{\hbar}\tilde{\oplus}R_\hbar)), \\
\end{diagram}\qquad\qquad
\begin{diagram}
& & M_{\hbar}\tilde{\oplus}N_\hbar & & \\
& \ldTo^{1_\hbar} & &
\rdTo{\bar{\varphi}_{M_\hbar}\tilde{\oplus}\bar{\varphi}_{N_\hbar}} & \\
M_{\hbar}\tilde{\oplus}N_\hbar & & & & \mathcal G_\hbar(\mathcal
F_\hbar(M_{\hbar}\tilde{\oplus}N_\hbar)), \\
\end{diagram}
\]
where $\bar{\varphi}_{W_\hbar},\dots,\bar{\varphi}_{N_\hbar}$ are given
by \eqref{morfismobase}.
Note that $\bar{\varphi}_{W_\hbar}\tilde{\oplus}\bar{\varphi}_{R_\hbar}$
and
$\bar{\varphi}_{M_\hbar}\tilde{\oplus}\bar{\varphi}_{N_\hbar}$ are
quasi-isomorphisms in $\mathcal H^{tf}_\infty (A_\hbar)$ as,
by definition,
$\bar{\varphi}^{(0)}_{W_\hbar},\dots,\bar{\varphi}^{(0)}_{N_\hbar}$ and
are homotopy equivalences, i.e. isomorphisms, in $\DD^{\infty}(A)$, as
noted in the subsection ``{\it Induction: $\mathcal S_r$}''
in the proof of thm.~\ref{Thm29}.
Commutativity of \eqref{megaborg!!} is easily
proved: we are just following the lines of the proof in thm.~\ref{Thm29}, with due changes.

If $r'\geq 3$ or $r\geq 3$, we need to further decompose
$\bar{X}^{r'-1}_\hbar$ and $\bar{Y}^{r-1}_\hbar$, repeating the above
considerations, a finite number of times. The proof of commutativity of
\eqref{megaborg!!} is conceptually analogous to the
one for the $r'=r=2$ case; we are just considering a
``trivially'' quantized version of the computations which appear at the
very end of the proof of thm.~\ref{Thm29}.

In summary, we have proven the equivalence $\mathcal G_\hbar\circ
\mathcal F_\hbar\simeq 1$ on $\triang^{\infty}_{A_\hbar}(A_\hbar)$.

 \subsubsection{ On $\mathcal F_\hbar\circ \mathcal G_\hbar\simeq 1$ on
$\triang^{\infty}_{B_\hbar}(K_\hbar)$ }

To prove $\mathcal F_\hbar\circ \mathcal G_\hbar\simeq 1$ on
$\triang^{\infty}_{B_\hbar}(K_\hbar)$
we begin by considering pair of objects in $\mathcal S'_1$, following
(once again!)
the lines in the proof of thm.~\ref{Thm29}. We define the derived functor
\[
\bar{\mathcal{G}}_\hbar: \DD^{\infty}_{tf}(B_{\hbar})\rightarrow
\DD^{\infty}_{tf}(A_{\hbar}),
\]
with $M_\hbar\mapsto
\bar{\mathcal{G}}_\hbar(M_\hbar):=M_\hbar\underline{\tilde{\otimes}}_{B_\hbar}\overline{K}_\hbar$

and we consider the following diagram:
{\scriptsize
\begin{eqnarray*}
\begin{diagram}
& & & & M_\hbar\underline{\tilde{\otimes}}_{B_\hbar}(\overline{K}_\hbar
\tilde{\underline{\otimes}}_{A_\hbar}
K_\hbar)\tilde{\underline{\otimes}}_{B_\hbar}\underline{K}_\hbar & & & & \\
& & & \ruTo^{1\tilde{\otimes}\nu_{A_\hbar}} & &
\luTo^{1\tilde{\otimes}\nu_{B_\hbar}\tilde{\otimes} 1} & & & \\
& & (M_\hbar\tilde{\underline{\otimes}}_{B_\hbar}
\overline{K}_\hbar)\bar{\underline{\otimes}}_{A_\hbar} A_\hbar & & & &
M_\hbar\tilde{\underline{\otimes}}_{B_\hbar}(B_\hbar\tilde{\underline{\otimes}}_{B_\hbar}\underline{K}_\hbar)&
& \\
& \ruTo^{\Phi_{M_\hbar\tilde{\underline{\otimes}}_{B_\hbar}
\overline{K}_\hbar} } & & & & & &
\luTo^{1\tilde{\otimes}\Phi_{\underline{K}_\hbar}} & \\
\bar{\mathcal{G}}_\hbar(M_\hbar)=M_\hbar\underline{\tilde{\otimes}}_{B_\hbar}
\overline{K}_\hbar & & & & & & & &
M_\hbar\underline{\tilde{\otimes}}_{B_\hbar}\underline{K}_\hbar=\mathcal
G_\hbar(M_\hbar)
\end{diagram} \nonumber\\
\end{eqnarray*}
}
where the quasi-isomorphisms
$\Phi_{M_\hbar\tilde{\underline{\otimes}}_{B_\hbar} \overline{K}_\hbar}$
and
$\Phi_{\underline{K}_\hbar}$ have been described in prop~\ref{Shuffle},
while the quasi-isomorphisms $\nu_{A_\hbar}$ and $\nu_{B_\hbar}$ appear
in Cor.~\ref{party2}.
All arrows in the above diagram are quasi-isomorphisms of strictly
unital topological $A_\infty$-$A_\hbar$-$A_\hbar$-bimodules.
Introducing the notation
\begin{eqnarray}
T_\hbar(M_\hbar):=M_\hbar\underline{\tilde{\otimes}}_{B_\hbar}(\overline{K}_\hbar
\tilde{\underline{\otimes}}_{A_\hbar}
K_\hbar)\tilde{\underline{\otimes}}_{B_\hbar}\underline{K}_\hbar,\label{Thbar}
\end{eqnarray}
the quasi-isomorphism
\begin{eqnarray*}
\varphi_{M_\hbar}:\bar{\mathcal{G}}_\hbar(M_\hbar)\rightarrow
\mathcal G_\hbar(M_\hbar)
\end{eqnarray*}
in $\DD^{\infty}_{tf}(A_{\hbar})$ is the composition
\begin{eqnarray}
\varphi_{M_\hbar}=\alpha_{M_\hbar} \circ\beta_{M_\hbar}, \label{Post}
\end{eqnarray}
choosing the roofs $(1_\hbar,\bar{\alpha}_{M_\hbar})$ and $(\bar{\beta}_{M_\hbar},1_\hbar)$ 
for $\alpha_{M_\hbar}$ and $\beta_{M_\hbar}$, where
\begin{eqnarray*}
\bar{\alpha}_{M_\hbar}=(1\tilde{\otimes}\nu_{A_\hbar})\circ
\Phi_{M_\hbar\tilde{\underline{\otimes}}_{B_\hbar} \overline{K}_\hbar}, &
\bar{\beta}_{M_\hbar}=(1\tilde{\otimes}\nu_{B_\hbar}\tilde{\otimes}
1)\circ (1\tilde{\otimes}\Phi_{\underline{K}_\hbar}).
\end{eqnarray*}

$\mathcal{F}_\hbar$ and $\mathcal{G}_\hbar$
are exact w.r.t.\ the triangulated structures on
$\DD^{\infty}_{tf}(B_{\hbar})$ and $\DD^{\infty}_{tf}(A_{\hbar})$; the
statement for $\bar{\mathcal{G}}_\hbar$
follows as well: the analysis is similar. Moreover
$\bar{\mathcal{G}}_\hbar(M_\hbar)\in
\triang^{\infty}_{A_\hbar}(A_\hbar)$, for any object $M_\hbar$ in
$\triang^{\infty}_{B_\hbar}(K_\hbar)$;
in fact
\begin{eqnarray}
\bar{\mathcal{G}}_\hbar(K_\hbar)=K_\hbar\underline{\tilde{\otimes}}_{B_\hbar}\overline{K}_\hbar\simeq
A_\hbar
\end{eqnarray}
in $\DD^{\infty}_{tf}(A_{\hbar})$ as it follows by considering the
quasi-isomorphism $\varphi_{K_\hbar}$ and using
$K_\hbar\underline{\tilde{\otimes}}_{B_\hbar}\underline{K}_\hbar\simeq A_\hbar$;
the statement for any object in $\mathcal S'_1$ easily follows; for
$r\geq 2$ we just need to look at exact triangles.
We continue with  the equivalence $\mathcal F_\hbar\circ
\bar{\mathcal{G}}_\hbar\simeq 1$ on $\triang^{\infty}_{B_\hbar}(K_\hbar)$;
 its proof is done by decomposing objects in $\mathcal S'_r$, $r\ge2$, into finite direct sums
of objects in $\mathcal S'_1$, as we did in the preceding subsection.
We need to prove the step $r=1$ explicitly.
Let $K_\hbar\langle i' \rangle[n']$, $K_\hbar\langle
j\rangle[r]$ be two objects in $\mathcal{S}'_1$ and let $g_\hbar:
K_\hbar\langle i' \rangle[n']\rightarrow K_\hbar\langle
j\rangle[r]$ be a morphism in $\DD^{\infty}(B_\hbar)$ represented by the
roof
\begin{eqnarray}
\begin{diagram}
& & K_\hbar\langle i \rangle[l] & & \\
& \ldTo^{s_\hbar} & & \rdTo{\bar{g}_\hbar} & \\
K_\hbar\langle i' \rangle[n'] & & & & K_\hbar\langle
j\rangle[r] \\
\end{diagram}\label{Megusta}
\end{eqnarray}
where $i,i',l,n',j,r\in\mathbb Z$. By degree reasons $\bar{g}_\hbar$ has a unique non trivial Taylor component
${\bar{g}_\hbar}^{(r)}$, for some $r\geq 0$.
The commutativity of the diagram
\begin{eqnarray*}
\begin{CD}
K_\hbar\langle i' \rangle[n'] @> g_\hbar >> K_\hbar\langle
j\rangle[r] \\
@VV \psi_{K_\hbar\langle i' \rangle[n']} V @VV \psi_{K_\hbar\langle
j\rangle[r]} V \\
\mathcal F_\hbar(\mathcal G_\hbar(K_\hbar\langle i' \rangle[n'])) @>
\mathcal F_\hbar(\bar{\mathcal{
G}}_\hbar(g_\hbar)) >>\mathcal F_\hbar(\bar{\mathcal{G}}_\hbar(
K_\hbar\langle
j\rangle[r])) \\
\end{CD}
\end{eqnarray*}
is proven once we show that
\begin{eqnarray}
\begin{CD}
K_\hbar\langle i\rangle[l] @>\bar{g}_\hbar>>
K_\hbar\langle j\rangle[r]\\
@VV \mathcal{Z}^{\hbar}_1 V @VV \mathcal{Z}^{\hbar}_2 V \\
(K_\hbar\underline{\tilde{\otimes}}_{B_\hbar} B_\hbar)\langle i\rangle[l] @.
(K_\hbar\underline{\tilde{\otimes}}_{B_\hbar} B_\hbar)\langle j\rangle[r]\\
@VV \mathcal{Z}^{\hbar}_3 V @VV \mathcal{Z}^{\hbar}_4 V \\
(K_\hbar\underline{\tilde{\otimes}}_{B_\hbar}
\underline{\End2}_{A_\hbar}(K_\hbar)^{op})\langle i\rangle[l] @.
(K_\hbar\underline{\tilde{\otimes}}_{B_\hbar}
\underline{\End2}_{A_\hbar}(K_\hbar)^{op})\langle j\rangle[r]\\
@VV \mathcal{Z}^{\hbar}_5 V @VV \mathcal{Z}^{\hbar}_6 V \\
(K_\hbar\underline{\tilde{\otimes}}_{B_\hbar}
(\overline{K}_\hbar\underline{\tilde{\otimes}}_{A_\hbar} K_\hbar )
)\langle i\rangle[l]
@>\mathcal F_\hbar(\bar{\mathcal{G}}_\hbar(\bar{g}_\hbar))>>
(K_\hbar\underline{\tilde{\otimes}}_{B_\hbar}
(\overline{K}_\hbar\underline{\tilde{\otimes}}_{A_\hbar} K_\hbar )
)\langle j\rangle[r]\\
\end{CD}\label{Kozo2}
\end{eqnarray}
commutes, where the morphisms $\mathcal{Z}^{\hbar}_i$ are those appearing in the
proof of thm.~\ref{Thm29}, section 7, with due changes.
Considering \eqref{Kozo2}, given any $X_\hbar\in \mathcal S'_1$, we
denote by $\bar{\psi}_{X_\hbar}: X_\hbar \rightarrow \mathcal
F_\hbar(\bar{\mathcal{G}}_\hbar(X_\hbar))$ the composition
\begin{eqnarray}
\bar{\psi}_{X_\hbar}=\mathcal{Z}^{\hbar}_5\circ\mathcal{Z}^{\hbar}_3\circ\mathcal{Z}^{\hbar}_1=
(1\tilde{\otimes}\mathcal R_\hbar)\circ(1\tilde{\otimes}\mathrm
R_{B_\hbar})\circ(\Phi_{K_\hbar} ),
\label{Caveza}
\end{eqnarray}
where $\Phi_{K_\hbar}$ is described in prop~\ref{Shuffle}, $\mathrm
R_{B_\hbar}$ is the quantized derived right action and $\mathcal R_\hbar$
appears in cor~\ref{miaumiau}. $\bar{\psi}_{X_\hbar}$ is a
quasi-isomorphism in $\mathcal H^{\infty}_{tf}(A_\hbar)$ as
$\bar{\psi}^{(0)}_{X_\hbar}$ is a quasi-isomorphism, i.e. a homotopy
equivalence in $\DD^{\infty}(A)$.

Using the decomposition of objects in $\mathcal S'_r$ into finite direct
sums of
objects in $\mathcal S'_1$ we finish the proof of the equivalence
$\mathcal F_\hbar \circ \bar{\mathcal G_\hbar}\simeq 1$
on $\triang^{\infty}_{B_\hbar}(K_\hbar)$. We
can repeat the analysis of the previous subsection almost {\it verbatim}.

The equivalence $\mathcal F_\hbar \circ \mathcal
G_\hbar\simeq 1$ on $\triang^{\infty}_{B_\hbar}(K_\hbar)$.
is proved following the same strategy that lead us to
$\mathcal F_\hbar \circ \bar{\mathcal G_\hbar}\simeq 1$ on
$\triang^{\infty}_{B_\hbar}(K_\hbar)$ : all we need is to consider the case
$r=1$ using $\mathcal F_\hbar \circ \bar{\mathcal{G}}_\hbar\simeq 1$ on
$\triang^{\infty}_{B_\hbar}(K_\hbar)$ . All computations and decompositions that appear in the proof of thm.~\ref{Thm29} 
can be repeated here, with due changes.

\subsubsection{Last part of the proof}

Exchanging $A_\hbar$ and $B_\hbar$; i.e. using the topological
$A_{\infty}$-$B_\hbar$-$A_\hbar$-bimodule $(\tilde{K}_\hbar,\mathrm
d_{K_\hbar})$ and
the new functors $\mathcal F_\hbar^{''}=
\cdot\underline{\tilde{\otimes}}_{B_\hbar} \tilde{K}_\hbar$ and
$\mathcal
G_\hbar^{''}=\cdot\underline{\tilde{\otimes}}_{A_\hbar}\underline{\tilde{K}}_\hbar$
we can prove the equivalence of the triangulated categories
$\triang^{\infty}_{B_\hbar}(B_\hbar) $ and
$\triang^{\infty}_{A_\hbar}(\tilde{K}_\hbar)$ with the same techniques
introduced above.

The statement on the thick subcategories follows by additivity of
$\mathcal F_\hbar$
and $\mathcal G_\hbar$ ($\mathcal F_\hbar^{''}$ and $\mathcal
G_\hbar^{''}$ as well), w.r.t.\ the
coproduct in $\DD^{\infty}_{tf}(A_\hbar)$ and
$\DD^{\infty}_{tf}(B_\hbar)$,
i.e. the direct sum of strictly unital topological  $A_{\infty}$-modules.

\end{proof}

\section{On triangulated categories}
In this section we collect some known facts on triangulated categories
and thickness.
We follow the expositions in \cite{GM} and \cite{Krause}.
Let us consider the pair $(\mathcal T,\Sigma)$, where $\mathcal T$ is an
additive category
and $\Sigma: \mathcal T\rightarrow \mathcal T$, $\Sigma(X):=\Sigma X$
an additive autoequivalence.

\begin{Def}
A triangle in $\mathcal T$ is a triple $(\alpha,\beta,\gamma)$ of
morphisms in $\mathcal T$
\[
X\stackrel{\alpha}{\rightarrow}Y\stackrel{\beta}{\rightarrow}Z\stackrel{\gamma}{\rightarrow}\Sigma
Z,
\]
and a morphism bewteen two triangles $(\alpha,\beta,\gamma)$,
$(\alpha',\beta',\gamma')$ is a triple $(\varphi_1,\varphi_2,\varphi_3)$
of morphisms in $\mathcal T$ s.t. the following diagram
\[
\begin{CD}
X @>\alpha >> Y @>\beta >> Z @> \gamma >>
\Sigma X \\
@VV \varphi_1 V @VV \varphi_2 V @VV \varphi_3 V @VV
\Sigma(\varphi_1)V \\
X' @>\alpha' >> Y' @>\beta' >> Z' @> \gamma' >> \Sigma X'
\end{CD}
\]
commutes.

\end{Def}

\begin{Def}
The category $\mathcal T$ is said to be triangulated if it is equipped
with a class of distinguished triangles, called the exact triangles,
satisfying the following axioms.
\begin{itemize}
\item (T1) A triangle isomorphic to an exact triangle is exact. For any
object $X$,
the triangle $0\rightarrow X\stackrel{1}{\rightarrow} X\rightarrow 0$ is
exact. Any morphism $\alpha: X\rightarrow Y$ in $\mathcal T$
can be completed to an exact triangle $X\stackrel{\alpha}{\rightarrow}
Y\stackrel{\beta}{\rightarrow} Z\stackrel{\gamma}{\rightarrow} \Sigma X$.
\item (T2) A triangle $(\alpha,\beta,\gamma)$ is exact if and only if
$(\beta,\gamma,-\Sigma \alpha)$ is exact.
\item (T3) Given two exact triangles $(\alpha,\beta,\gamma)$ and
$(\alpha',\beta',\gamma')$, each pair of morphisms
$\varphi_1$ and $\varphi_2$ satisfying
$\varphi_2\circ \alpha=\alpha'\circ \varphi_1$ can be completed (not
necessarily uniquely) to a morphism
of triangles $(\varphi_1,\varphi_2,\varphi_3)$:
\[
\begin{CD}
X @>\alpha >> Y @>\beta >> Z @> \gamma >>
\Sigma X \\
@VV \varphi_1 V @VV \varphi_2 V @VV \varphi_3 V @VV
\Sigma(\varphi_1)V \\
X' @>\alpha' >> Y' @>\beta' >> Z' @> \gamma' >> \Sigma X'
\end{CD}
\]
\item (T4)\footnote{this is the celebrated $octahedral$ $axiom$.}
Given exact triangles $(\alpha_1,\alpha_2,\alpha_3)$,
$(\beta_1,\beta_2,\beta_3)$ and $(\gamma_1,\gamma_2,\gamma_3)$ with
$\gamma_1=\beta_1\circ\alpha_1$,
there exists an exact triangle $(\delta_1,\delta_2,\delta_3)$ making the
following diagram
\[
\begin{CD}
X @>\alpha_1 >> Y @>\alpha_2 >> U @> \alpha_3
 >> \Sigma X \\
@| @VV \beta_1 V @VV \delta_1 V @|\\
X @>\gamma_1 >> Z @>\gamma_2 >> V @> \gamma_3 >> \Sigma X \\
@. @VV \beta_2 V @VV \delta_2 V @VV \Sigma\alpha_1 V\\
@. W @= W @> \beta_3 >> \Sigma Y \\
@. @VV \beta_3 V @VV \delta_3 V @. \\
@. \Sigma Y @>\Sigma(\alpha_2)>> \Sigma U @.\\
\end{CD}
\]
commutative.
\end{itemize}
\end{Def}
If the category $\mathcal T$ satisfies only the axioms
$(T1)$-$(T2)$-$(T3)$,
then it is said to be a pre-triangulated category.

\begin{Def}
Let $\mathcal T$ be a pre-triangulated category and $\Ab$ be the
category of abelian groups.
A functor $\mathcal F: \mathcal T\rightarrow \mathcal A$, with $\mathcal
A$ abelian, is said to be cohomological if
it sends each exact triangle in $\mathcal T$ to an exact sequence in
$\mathcal A$.
\end{Def}
\begin{Lem}
For each $\in\mathcal T$, the representable functors
\[
\Hom2_{\mathcal T}(X,\cdot):\mathcal T\rightarrow \Ab,
~~~\Hom2_{\mathcal T}(\cdot,X):\mathcal T^{op}\rightarrow \Ab,
\]
are cohomological.
\end{Lem}

 From the above lemma it follows
\begin{Lem}\label{isoTR}
Let $(\varphi_1,\varphi_2,\varphi_3)$ be a morphism between exact
triangles in $\mathcal T$.
If two maps in $\{\varphi_1,\varphi_2,\varphi_3\}$ are isomorphisms,
then also the third.
\end{Lem}

\begin{Def}
Let $(\mathcal T,\Sigma_1)$ and $(\mathcal U,\Sigma_2)$ be triangulated
categories.
An exact functor $\mathcal T\rightarrow \mathcal U$ is a pair $(\mathcal
F,\eta)$ consisting
of a functor $\mathcal F:\mathcal T\rightarrow \mathcal U$ and a natural
isomorphism $\eta: \mathcal F\circ \Sigma_1 \rightarrow \Sigma_2\circ
\mathcal F$ s.t.,
for every exact triangle $X\stackrel{\alpha}{\rightarrow}
Y\stackrel{\beta}{\rightarrow}
Z\stackrel{\gamma}{\rightarrow}\Sigma_1 X$ in $\mathcal T$, the triangle
\[
\mathcal F(X)\stackrel{\mathcal F(\alpha)}{\rightarrow} \mathcal
F(Y)\stackrel{\mathcal F(\beta)
}{\rightarrow} \mathcal F(Z)\stackrel{\mathcal
\eta\circ F(\gamma)}{\rightarrow}\Sigma_2 \mathcal F(X)
\]
is exact in $\mathcal U$.
\end{Def}
The following basic  example is well studied in \cite{Kel4}.
\begin{Example}
Let $(\mathcal T,\Sigma)$ be a triangulated category.
The autoequivalence $(\Sigma,-1)$ is an exact functor w.r.t. the triangulated structure on $\mathcal T$. 
\end{Example}

\subsubsection{On triangulated subcategories}

Let $(\mathcal T,\Sigma)$ be a triangulated category.
\begin{Def}
A non-empty full additive subcategory $\mathcal C$ is a triangulated
subcategory if
\begin{itemize}
\item (S0) $\mathcal S$ is strict; any object isomorphic to an object in
$\mathcal S$ belongs to $\mathcal S$.
\item (S1) $\Sigma^n X\in\mathcal C$ for all $X\in \mathcal C$ and
$n\in\mathbb Z$.
\item (S2) Let $X\rightarrow Y\rightarrow Z\rightarrow \Sigma X$ be any
exact triangle in $\mathcal T$.
If any two objects from $\{X,Y,Z\}$ belong to $\mathcal C$, so also the
third.
\end{itemize}
\end{Def}
A triangulated subcategory $\mathcal C$ inherits a canonical
triangulated structure from $\mathcal T$.
Let $\mathcal U$ and $\mathcal V$ be classes whose objects are
(isomorphism classes of) objects in $\mathcal
T$, where $\mathcal T$ is triangulated.
The class $\mathcal U \star \mathcal V$ is defined as follows:
\[
\mathcal U \star \mathcal V:=\{X\in \mathcal T : U\rightarrow
X\rightarrow V\rightarrow \Sigma U ~\mbox{exact triangle in }
\mathcal T,~ U\in\mathcal U, V\in\mathcal V \}.
\]
The composition $\star$ is associative by the octahedral axiom $(T4)$.
The following notation
\[
\mathcal S_r= \underbrace{S_1\star S_1\star\dots\star S_1}_{r-\mbox{times}}
\]
is unambiguous for $r\geq 1$, for any class $\mathcal S_1$ of objects in
$\mathcal T$.
The objects in $\mathcal S_r$ are called the
extensions of lenght $r$ of objects of $\mathcal S_1$. If $\mathcal T$ is a triangulated category and $M$ is an object in $\mathcal T$,
 the full triangulated  subcategory  generated by $M$ consists of all objects belonging to $\mathcal S_r$  ($r\geq 1$, as above) with 
$\mathcal S_1=\{ M[i],i\in\mathbb Z\}$ (in $\mathcal S_1$ we consider equivalence classes of isomorphic objects) .
Such triangulated subcategory is the smallest full triangulated subcategory in $\mathcal T$ containing $M$.

Its thickening  is the full triangulated subcategory of $\mathcal T$ consisting of all objects $X$ in $\mathcal T$, s.t.
there exist an object $Z$ in the triangulated subcategory generated by $M$ with $Z\simeq X\oplus Y$. 
The thickening is closed under direct summands; actually it is the smallest full triangulated
 subcategory in $\mathcal T$ containing the triangulated subcategory generated by $M$  and closed under direct summands.


\begin{thebibliography}{99}


\bibitem{Q}
F. Bayen, M. Flato, C. Fr\o nsdal, A. Lichnerowicz and D. Sternheimer,
``Deformation theory and quantization I, II'',
Ann. Phys. {\bf 111} (1978).



\bibitem{CFR}
D.~Calaque,G.~Felder, C.A.~Rossi,
``Deformation Quantization with generators and relations'',
available at http://arxiv.org/abs/0911.4377.

 

\bibitem{CFFR}
D.~Calaque,G.~Felder, A.~Ferrario, C.A.~Rossi,
`` Bimodules and branes in deformation quantization``,
Comp. Math. {\bf 147} , 105-160 (2011).


\bibitem{CT1}
A.~Caldararu, J.~Tu,
``Curved A-infinity algebras and Landau-Ginzburg models'',
available at http://arxiv.org/abs/1007.2679.

 
 \bibitem{CF0}
  A.~S.~Cattaneo and G.~Felder,
``Coisotropic submanifolds in Poisson geometry and branes in the Poisson sigma model'',
Lett.Math.Phys. {\bf 69} (2004), 157-175.

 
\bibitem{CF}
  A.~S.~Cattaneo and G.~Felder,
  ``Relative formality theorem and quantisation of coisotropic submanifolds'',
   Adv. Math. {\bf 208}, 521-548 (2007).






 

 


\bibitem{Kad}
T.V.~Kadeishvili,
``The algebraic structure in the homology of an A($\infty$)
algebra''(Russian),
Soobshch. Akad. Nauk. Gruzin, {\bf SSR} 108, 1982.

\bibitem{Keller-mail} 
B.~Keller, private communication.


\bibitem{Kel}
B.~Keller,
`` Introduction to $A_\infty$-algebras and modules'',
Homology Homotopy Appl. {\bf 3}, (2001).

\bibitem{Kel2}
B.~Keller,
``Hochschild cohomology and derived Picard group'',
Journal of Pure and Applied Algebra, {\bf 190} (2004).


\bibitem{Kel4}
B.~Keller,
``Calabi-Yau triangulated categories'',
Trends in representation theory of algebras and related topics,
EMS Ser. Congr. Rep., Eur. Math. Soc., Zürich, 2008.



\bibitem{Kenji}
K.~Levefre-Hasegawa,
``Sur les $A_{\infty}$-categories'',
PhD Thesis, Univ. Paris 7 (2007).

\bibitem{Kont}
M.~Kontsevich,
``Deformation Quantization of Poisson manifolds, I'',
Lett. Math. Phys. 66 (2003).

\bibitem{Kont2}
M.~Kontsevich,
``Homological Algebra of Mirror Symmetry'',
talk at ICM, Zurich (1994).

\bibitem{Kont3}
M.~Kontsevich, J.~Soibelman
``Notes on $A_\infty$-algebras, $A_\infty$-categories and
non-commutative geometry'',
Lecture Notes in Phys., 757, Springer, Berlin, (2009).

\bibitem{Krause}
H.~Krause,
``Derived categories, resolutions, and Brown representability'',
arXiv:math/0511047v3


\bibitem{AF}
A.~Ferrario,
``Poisson sigma model with branes and hyperelliptic Riemann surfaces'',
J. Math. Phys. {\bf 49}, 2008.

\bibitem{ART}
A.~Ferrario, C.A.~Rossi, T.~Willwacher,
``A note on the Koszul complex in deformation quantization'',
Lett. Math. Phys. {\bf 95}, 27-39 (2011).

\bibitem{GJ}
E.~Getzler, J.D.S.~Jones,
``A-infinity algebras and the cyclic bar complex'',
Illinois J. Math. 34 (1990), no. 2, 256–283.

\bibitem{GM}
S.I.~Gelfand, Y.~Manin,
``Methods of Homological Algebra'',
Springer Monographs in Mathematics.

 
\bibitem{Halperin}
Y.~ Félix, S.~ Halperin, J.C.~ Thomas,
``Rational Homotopy Theory'', 
Springer Verlag, (2001).


\bibitem{Neeman}
A.~Neeman,
``Triangulated categories'',
AMS, Princeton University Press.

\bibitem{Carlo}
C.~Rossi,
``The Chevalley-Eilenberg complex quantizes the Koszul complex'',
math.QA/1105.5972.


\bibitem{Shoikhet}
B.~Shoikhet,
``Koszul duality in deformation quantization and Tamarkin's approach to
Kontsevich formality'',
arXiv:0805.0174.



\bibitem{Sta}
J.~Stasheff,
``Homotopy associativity of H-spaces I'',
Trans. Amer. Math. Soc. {\bf 105} (1962)

\bibitem{Tradler}
T.~Tradler,
``Infinity-inner-products on $A_\infty$-algebras'',
J. Homotopy Relat. Struct. 3 (2008), no. 1, 245-271.


\bibitem{Zhang}
D.M.~Lu, J. H.~Palmieri, Q.-S.~Wu, J. J.~Zhang,
``Koszul Equivalences in $A_\infty$-Algebras '',
New York J. Math., {\bf 14}, 2008.


\end{thebibliography}
\end{document}